\documentclass[11pt]{article}
\usepackage{amstext,amssymb,amsmath,amsbsy}
\usepackage{color,soul}
\usepackage{hyperref}
\usepackage{amscd}
\usepackage{amsfonts}
\usepackage{indentfirst}
\usepackage{verbatim}
\usepackage{amsmath}
\usepackage{amsthm}
\usepackage{enumerate}
\usepackage{graphicx}
\usepackage{color}
\usepackage[OT1]{fontenc}
\usepackage[latin1]{inputenc}
\usepackage[english]{babel}
\usepackage{amssymb}
\usepackage{subfig}
\usepackage{algorithm}
\usepackage{algpseudocode}

\newtheorem{Lemma}{Lemma}[section]

\newtheorem{remark}{Remark}[section]

\newtheorem*{Assumption*}{Assumption}

\newtheorem{problem}{Problem}[section]
\newtheorem*{problem*}{Problem}
\numberwithin{equation}{section}

\input{tcilatex}
\begin{document}

\title{Convexification-based globally convergent numerical method for a 1D
coefficient inverse problem with experimental data}
\author{Michael V. Klibanov\thanks{%
Department of Mathematics and Statistics, University of North Carolina at
Charlotte, Charlotte, NC 28223, USA, mklibanv@uncc.edu (corresponding
author), tle55@uncc.edu, loc.nguyen@uncc.edu } \and Thuy T. Le%
\footnotemark[1] \and Loc H. Nguyen\footnotemark[1] \and Anders Sullivan%
\thanks{%
US Army Research Laboratory, 2800 Powder Mill Road, Adelphi, MD 20783-1197,
USA, (anders.j.sullivan.civ@mail.mil, lam.h.nguyen2civ@mail.mil)} \and Lam
Nguyen\footnotemark[2]}
\date{}
\maketitle

\begin{abstract}
To compute the spatially distributed dielectric constant from the
backscattering data, we study a coefficient inverse problem for a 1D
hyperbolic equation. To solve the inverse problem, we establish a new
version of Carleman estimate and then employ this estimate to construct a
cost functional which is strictly convex on a convex bounded set with an
arbitrary diameter in a Hilbert space. The strict convexity property is
rigorously proved. This result is called the convexification theorem and is
considered as the central analytical result of this paper. Minimizing this
convex functional by the gradient descent method, we obtain the desired
numerical solution to the coefficient inverse problems. We prove that the
gradient descent method generates a sequence converging to the minimizer and
we also establish a theorem confirming that the minimizer converges to the
true solution as the noise in the measured data and the regularization
parameter tend to zero. Unlike the methods that are based on optimization,
our convexification method converges globally in the sense that it delivers
a good approximation of the exact solution without requiring any initial
guess. Results of numerical studies of both computationally simulated and
experimental data are presented.
\end{abstract}

%\affil[]{Department of Mathematics and Statistics, University of North Carolina at
%Charlotte, Charlotte, NC 28223, USA}



\section{Introduction}

\label{sec intro}

We develop a new version of the convexification method to numerically solve
a highly nonlinear and severely ill-posed inverse problem for a 1D
hyperbolic equation. Applications of this technique are in detection and
identification of explosives, see some details in Section 7. This paper
belongs to a series of works that establish a variety of versions of the
convexification method to solve nonlinear inverse problems for many partial
differential equations \cite{VoKlibanovNguyen:IP2020,
KhoaKlibanovLoc:SIAMImaging2020, Klibanov:jiip2017,
KlibanovIoussoupova:SMA1995, KlibanovKolesov:cma2019,
KlibanovKolNguyen:SISC2019, KlibanovKolesov:ip2018, KlibanovLiZhang:ip2019,
Klibanov:ip2020, KlibanovLiZhang:SIAM2019, SmirnovKlibanovNguyen:IPI2020}.
The key point of the convexification method for each inverse problem in
those publications is to use a suitable weight function to construct a
globally strictly convex weighted Tikhonov-like functional. The weight is
the Carleman Weight Function (CWF), i.e. the function which is involved as
the weight in the Carleman estimate for the corresponding Partial
Differential Operator. The unique minimizer of such a functional directly
yields the desired numerical solution of that nonlinear inverse problem. The
above mentioned global strict convexity guarantees that we can solve that
nonlinear inverse problem without any advanced knowledge of the true
solutions. Therefore, we say that our convexification method is globally
convergent. By \textquotedblleft globally convergent", we mean:

\begin{enumerate}
\item There exists a theorem rigorously confirming that our method delivers
at least one point in a sufficiently small neighborhood of the exact
solution without requiring a good initial guess of the true solution;

\item This theorem is verified numerically.
\end{enumerate}

In this paper, item 1 is reached by establishing a new Carleman estimate and
applying it to prove a new version of the convexification theorem. Item 2 is
reached for both computationally simulated data and experimental data.

We also refer here to some recent publications of another research group 
\cite{BAUDOUIN:SIAMNumAna:2017,Baudouin:2020,Boulakia:esaim2021} and the
publication \cite{LeNguyen:2020} by members of our group. These papers work
on two different versions of the convexification. CWFs still play a crucial
role in both versions. We now describe the main difference between our above
cited works and ones of this group. Papers \cite%
{BAUDOUIN:SIAMNumAna:2017,Baudouin:2020,Boulakia:esaim2021,LeNguyen:2020}
work for the case when at least one of initial conditions is not vanishing.
Unlike this, we consider in almost all above cited publications the case
when the initial condition in a hyperbolic equation is the Dirac $\delta -$%
function and similar conditions for CIPs for the Helmholtz equation. In this
regard, the only exception is the publication \cite{Klibanov:ip2020}, which
also works for the case when a sort of an initial condition is not vanishing.

Let $n(x)$ be the refractive index and let $c(x)=n^{2}(x)$ for all $x\in 
\mathbb{R}$ be the spatially distributed dielectric constant. If we scale
the speed of light traveling in the air or vacuum to be $1$, then $1/\sqrt{%
c(x)}$ is the speed of light in the medium. Let $\epsilon $ and $M$ be two
fixed numbers with $0<\epsilon \ll 1<M<\infty $. Assume that the spatially
distributed dielectric constant $c$ belongs to $C^{3}(\mathbb{R})$ and that 
\begin{equation}
c(x)=\left\{ 
\begin{array}{ll}
\in \lbrack \underline{c},\overline{c}] & \mbox{if }x\in \lbrack \epsilon
,M], \\ 
1 & \mbox{otherwise}%
\end{array}%
\right.  \label{1.1}
\end{equation}%
for some known constants $0<\underline{c}<\overline{c}<\infty .$ The
smoothness condition $c\in C^{3}(\mathbb{R})$ is imposed only for the
theoretical part while it can be relaxed in the numerical study. Let $%
u=u(x,t)$, $(x,t)\in \mathbb{R}\times \lbrack 0,\infty )$, be the solution
to the following initial value problem 
\begin{equation}
\left\{ 
\begin{array}{rcll}
c(x)u_{tt}(x,t) & = & u_{xx}(x,t) & (x,t)\in \mathbb{R}\times (0,\infty ),
\\ 
u(x,0) & = & 0 & x\in \mathbb{R}, \\ 
u_{t}(x,0) & = & \delta _{0}(x) & x\in \mathbb{R},%
\end{array}%
\right.  \label{main eqn}
\end{equation}%
where $\delta _{0}$ is the Dirac function with its support $\left\{
0\right\} $. 
%In this paper, we propose a globally convergent numerical method, named as the convexification, to solve the following coefficient inverse problem. 
The problem of our interest is formulated as follows.

\begin{problem}[Coefficient inverse problem]
Let $T>0$ be the length of the time interval. Measuring the functions $%
g_{0}(t)$ and $g_{1}(t),$ 
\begin{equation}
g_{0}(t)=u(\epsilon ,t)\quad \mbox{and }\quad g_{1}(t)=u_{x}(\epsilon ,t)
\label{1.3}
\end{equation}%
for $t\in \lbrack 0,T+\epsilon ],$ determine the function $c(x)$ for all $%
x\in \lbrack \epsilon ,M].$ \label{cip}
\end{problem}

We refer the reader to \cite{deHoopetal2018, Katchalovetal:2001,
Korpela,Montalto:2014} for some uniqueness and stability results for
coefficient inverse problems that are similar to Problem \ref{cip} to
identify $c,$ given the Dirichlet-to-Neumann map data. The uniqueness and
stability results for Problem \ref{cip} follow from \cite{Romanov:VNU1986}
(chapter 2) as well as directly from our computational method in this paper.
In \cite{Kab1} the Gelfand-Levitan method \cite{GL} was numerically
implemented for a similar CIP. Next, this method was extended to the 2D case 
\cite{Kab1,Kab2}.

Problem \ref{cip} arises from the following experiment. Let an emitter sends
an electric wave into the inspected area, in which the target we want to
identify is hidden. Then, we measure the back scattering wave using a
detector located near the emitter. An example of this device is the Forward
Looking Radar built by the US Army Research Laboratory (ARL) \cite%
{NguyenWong:pspie2007}. This radar device is placed on the top of a moving
vehicle during the data collecting process, see \cite{NguyenWong:pspie2007}
for more details. 
Since the given data has one dimension for each target, reconstructing $d-$%
dimensional function with $d>1$ is impossible. Thus, we have no choice but
to model the wave propagation by a 1D hyperbolic equation. This 1D model was
verified numerically multiple times in the past in the sense that it can be
used to successfully compute the dielectric constants of explosive-like
targets from experimental data provided by ARL, see \cite%
{Karchevskyetal:2013, KlibanovAlex:SIAMjam:2017,
KlibanovKolesov:ip2018,KlibanovLoc:ipi2016, Kuzhuget:IP2012, Smirnov:ip2020}.

\begin{figure}[h!]
\begin{center}
\subfloat[The target is placed in the
air]{\includegraphics[width=.4\textwidth]{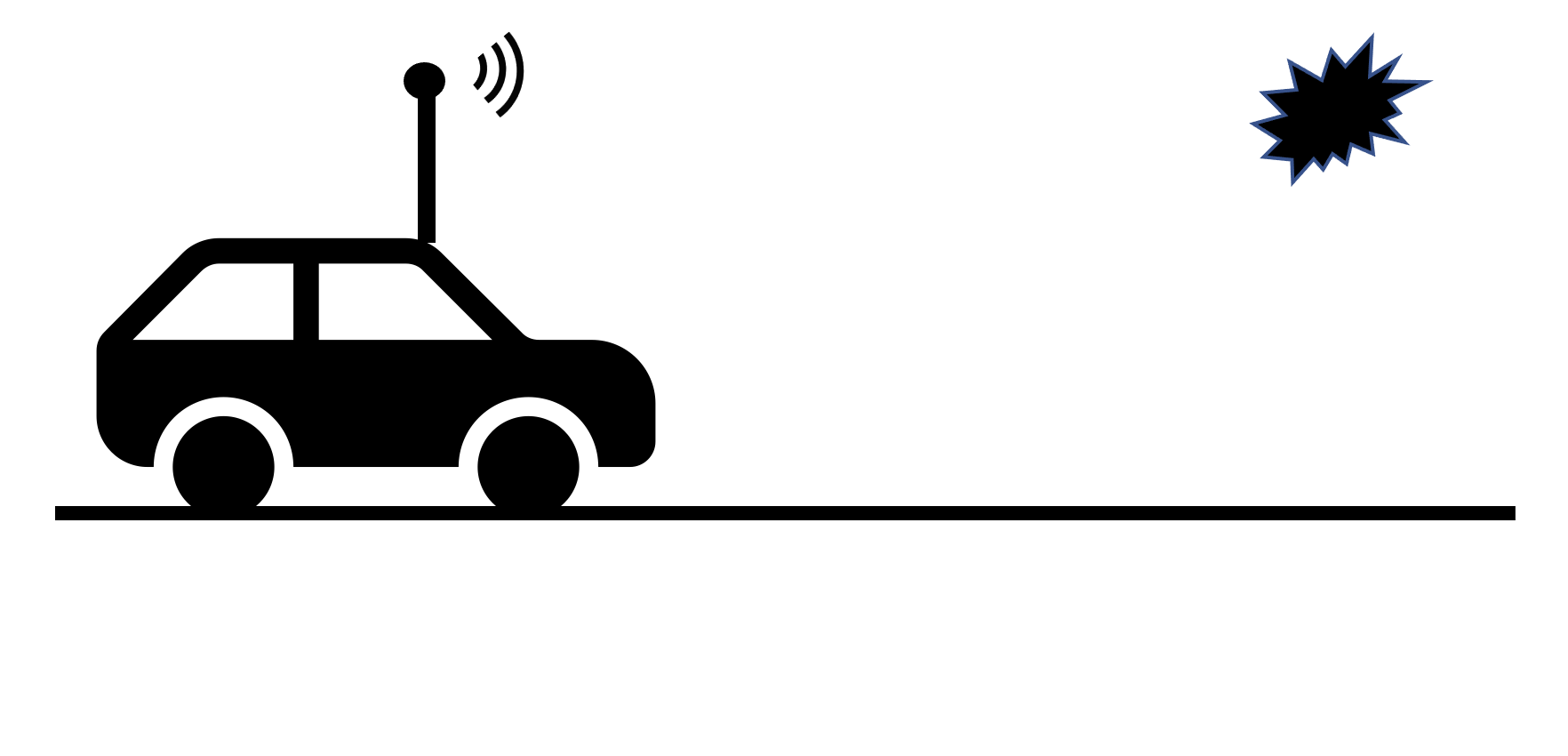}} \quad 
\subfloat[The
target is buried under the
ground]{\includegraphics[width=.4\textwidth]{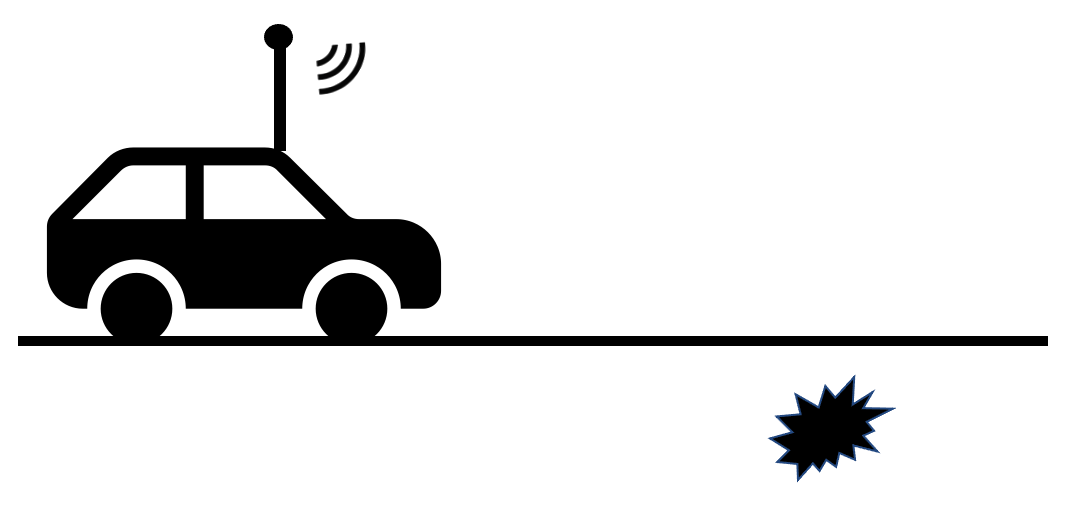}}
\end{center}
\caption{\textit{The schematic for the data generating and collecting
device. A device, called radar, emits an acoustic source and then collect
the time-dependent backscattering wave. In the physical experiment, we
consider two cases: (a) the target is placed in the air and (b) the target
is buried a few centimeters under the ground. }}
\label{fig sch}
\end{figure}

By solving Problem \ref{cip}, we obtain the spatially distributed dielectric
constant. This computed dielectric constant provides the location and some
information about the constituent material of the target. This problem has
applications in detecting antipersonnel explosive devices. The latter is one
of important Army's interests. In this experiment, only $g_{0}$ is measured
while the data $g_{1}$ is missing. However, in Section \ref{sec exp}, we
explain how to approximate the missing function $g_{1}$. The schematic
diagram of the data collection is displayed on Figure \ref{fig sch}.We have
also discovered recently that Problem \ref{cip} plays the key role in the
nonlinear synthetic-aperture radar (SAR) imaging, including SAR experimental
data \cite{Klibanov2021ThroughtheWallNS,Klibanov:2ndSAR2021}.

Natural approaches for solutions of nonlinear inverse problems, that are
widely used in the scientific community, are based on the least-squares
optimization. However, the use of optimization-based methods is limited to
the case when a good initial guess for the true solution of Problem \ref{cip}
is known. However, it is rarely available in the reality. This requirement
is due to the fact that the those least-squares functionals are non convex
and typically have multiple local minima and ravines, see, e.g. \cite[Figure
1]{ScalesSmithFischerLjcp1992} for a convincing example of this well-known
challenge. Hence, the least-squares optimization method is not applicable to
solve Problem \ref{cip}. Another approach to solve Problem \ref{cip} is the
use of the Born approximation or Born series. This approach is effective if
the true dielectric constant is a sufficiently small perturbation of a known
background function. Hence, the methods based on Born approximation or/and
Born series work for the case when the size of the target is small and the
contrast $c_{\mathrm{target}}/c_{\mathrm{bckgr}}\ll 1$ where $c_{\mathrm{%
target}}$ is the dielectric constant of the target and $c_{\mathrm{bckgr}}$
is the dielectric constant of the background (or the environment around the
target). For example, it was demonstrated numerically on \cite%
{Klibanov2021ThroughtheWallNS} that the Born approximation is not capable to
deliver accurate values of dielectric constants of targets for SAR-like data
in the case of high target/background contrasts.

To overcome these limitations, Klibanov and Ioussoupova introduced the
convexification method \cite{KlibanovIoussoupova:SMA1995}. Since then, it
has been intensively applied to solve nonlinear coefficient inverse problems 
\cite{KlibanovNik:ra2017, VoKlibanovNguyen:IP2020,
KhoaKlibanovLoc:SIAMImaging2020, Klibanov:sjma1997, Klibanov:nw1997,
Klibanov:ip2015, KlibanovKolesov:cma2019, KlibanovKolNguyen:SISC2019,
KlibanovLiZhang:ip2019, KlibanovLiZhang:SIAM2019,
SmirnovKlibanovNguyen:IPI2020, Smirnov:ip2020}. The reconstructions due to
the convexification are successful even for the challenging case of
experimental data \cite{VoKlibanovNguyen:IP2020,
KlibanovKolNguyen:SISC2019,Smirnov:ip2020}. The main idea of the
convexification is to employ CWFs and Carleman estimates to convexify the
least-squares functional. In other words, when we employ a suitable CWF in
the mismatch functional, the resulting functional is strictly convex on a
convex bounded set of an arbitrary diameter in an appropriate Hilbert space.
The minimizer of this strictly convex functional, which can be found without
an initial guess, is an approximation of the desired solution. Hence, we
claim and then prove that our numerical method is \textit{globally convergent%
}, see the first paragraph of this section for the definition of globally
convergent.

The original idea of applying Carleman estimates to coefficient inverse
problems was first published in \cite{BukhgeimKlibanov:smd1981} back in 1981
by Bukhgeim and Klibanov to prove uniqueness theorems for a wide class of
coefficient inverse problems. Some follow up publications can be found in,
e.g. \cite{Isakov:bookSpringer2017, Klibanov:ip1992,
KlibanovTimonov:u2004,LeNguyen:2020, LocNguyen:ip2019,
NguyenLiKlibanov:IPI2019, TriggianiandYao:amo2002}. Surveys on the
Bukhgeim-Klibanov method can be found in \cite%
{Klibanov:jiipp2013,Yamamoto:ip2009}, also, see section 1.10 of the book 
\cite[Chapter 1]{BeilinaKlibanovBook}. It was discovered later in \cite%
{KlibanovIoussoupova:SMA1995} that the idea of \cite%
{BukhgeimKlibanov:smd1981} can be used to develop globally convergent
numerical methods for coefficient inverse problems using the convexification.

The inverse problem in this paper, Problem \ref{cip}, is identical with the
inverse problem in \cite{SmirnovKlibanovNguyen:IPI2020,Smirnov:ip2020}. The
convexification method in \cite{SmirnovKlibanovNguyen:IPI2020,Smirnov:ip2020}
is effective but there are some rooms to improve. The method of \cite%
{SmirnovKlibanovNguyen:IPI2020,Smirnov:ip2020} has two stages. On stage 1,
the authors used the well-known change of variables as in \cite[Chapter 2, §7%
]{Romanov:VNU1986} to reduce the original inverse problem to the inverse
problem of computing the potential of a 1D hyperbolic equation. This
resulting inverse problem to compute the potential can be solved using one
of two versions of the convexificaton method of either \cite%
{SmirnovKlibanovNguyen:IPI2020} or \cite{Smirnov:ip2020}. On stage 2, the
authors computed the dielectric constant from the knowledge of the
reconstructed potential of stage 1. This second stage is a quite complicated
one, due to that change of variables, and caused many difficulties in its
numerical implementation. This motivates us to propose a different version
of the convexification method in this paper. In this paper, we solve Problem %
\ref{cip} directly, i.e. without the change of variable of \cite[Chapter 2, §%
7]{Romanov:VNU1986}. By this, the numerical implementation is significantly
simplified.

The key points that guarantee the success of our method involve:

\begin{enumerate}
\item The derivation of a nonlinear and non local partial differential
equation without the presence of the unknown coefficient.

\item Two new Carleman estimates for this equation.

\item A new version of the convexification method to solve the above
equation.

\item The theorem ensuring the global strict convexity of the cost
functional constructed by the convexification method of item 3.

\item The theorem, which guarantees the global convergence of the gradient
descent method of the minimization of the strictly convex functional
resulting from the convexification method.
\end{enumerate}

%(1) the derivation of a nonlinear and nonlocal differential equation without the presence of the dielectric constant, (2) a new Carleman estimate for this equation, (3) a new version of the  convexification method to solve the equation above, (4) the convergence of the convexification method and (5) the convergence of the gradient descent method to minimize the convex functional constructed the the convexification method.

Since Problem \ref{cip} is a coefficient inverse problem for a 1D hyperbolic
equation, we mention here works \cite%
{KlibanovLiZhang:SIAM2019,KlibanovLiZhang:arxiv2020} in which the authors
solved 3D versions of Problem \ref{cip}. In \cite{KlibanovLiZhang:SIAM2019}
a coefficient inverse problem for the 3D analog of equation (\ref{main eqn})
with a single location of the source was solved numerically via a version of
the convexification method. In \cite{KlibanovLiZhang:arxiv2020} a
coefficient inverse problem of the 3D analog of equation (\ref{main eqn})
with the point source running along a straight line was solved via the
linear integral equation invented by M.M. Lavrent'ev in 1964 \cite%
{Lavrentev:1964}; also, see formula (7.18) of the book \cite%
{Lavrentiev:AMS1986} for this equation. A new numerical method for the
solution of the Lavrent'ev equation was proposed in \cite%
{KlibanovLiZhang:arxiv2020}.

This paper is organized as follows. In Section \ref{sec q}, we derive an
important equation. Solution of this equation can be directly used to
compute the desired dielectric constant. In Section \ref{sec Car}, we prove
a new Carleman estimate. This estimate is an important generalization of the
one in \cite{SmirnovKlibanovNguyen:IPI2020}. In Section \ref{sec conv}, we
prove the convexification theorem. Also in Section \ref{sec conv}, we prove
the global convergence of the gradient descent method of the minimization of
the globally strictly convex cost functional constructed by the
convexification method. In Section \ref{sec sim}, we present the numerical
results obtained for computationally simulated data. In Section \ref{sec exp}%
, we present the numerical results obtained from experimental data.
Concluding remarks are made in Section \ref{sec rem}.

\section{A Partial Differential Equation in Which the Unknown Coefficient is
Not Present}

\label{sec q}

For each $x\in \mathbb{R}$, define 
\begin{equation}
\tau (x)=\int_{0}^{x}\sqrt{c(s)}ds.  \label{tau}
\end{equation}%
The function $\tau (x)$ is the travel time. This is the time the wave needs
to propagate from the source position $\left\{ x=0\right\} $ to the point $x$%
, see \cite[Chapter 2, §7]{Romanov:VNU1986}. It is well-known that $\tau $
satisfies the \textit{eikonal} equation 
\begin{equation}
|\tau ^{\prime}(x)|^2=c(x), \quad \mbox{for all }x>\epsilon .  \label{eik}
\end{equation}%
Since $c(x)=1$ for all $x<\epsilon $, see condition (\ref{1.1}), then 
\begin{equation}
\tau (x)=x, \quad \mbox{for all }x\leq \epsilon .  \label{2.444444}
\end{equation}%
In particular, 
\begin{equation}
\tau (\epsilon )=\epsilon .  \label{tauepsilon}
\end{equation}%
The following lemma is important for derivation of the numerical method in
this paper.

\begin{Lemma}
The function $u(x, t)$ has the form 
\begin{equation}
u(x, t) = \frac{H(t - |\tau(x)|)}{2 c^{1/4}(x)} + \widehat u(x, t), \quad
(x, t) \in \mathbb{R} \times (0, \infty)  \label{2.3}
\end{equation}
where $\widehat u$ is a function in $C^2(t \geq \tau(x))$ and $\widehat u(x,
\tau(x)) = 0$. As a result, 
\begin{equation}
\lim_{t \to \tau(x)^+} u(x, t) = \frac{1}{2 c^{1/4}(x)}, \quad 
\mbox{for all
} x > \epsilon.  \label{u2c}
\end{equation}
\label{lem 2.1}
\end{Lemma}

\textbf{Proof.} By (\ref{tau}) $\tau \left( x\right) $ is an increasing
function and has an inverse. Recall the well-known change of variable, see
formulas (8) and (9) in \cite{Smirnov:ip2020}, 
\begin{align}
v(x,t)& =\displaystyle u(\tau ^{-1}(x),t)c^{1/4}(\tau ^{-1}(x)), & & %
\mbox{for all }x\in \mathbb{R},t\in (0,\infty ),  \label{2.4444} \\
\ S(x)& =\displaystyle c^{-1/4}(\tau ^{-1}(x)), & & \mbox{for all }x\in 
\mathbb{R},  \label{2.5} \\
r(x)& =\displaystyle\frac{S^{\prime \prime }(x)}{S(x)}-2\Big[\frac{S^{\prime
}(x)}{S(x)}\Big]^{2}, & & \mbox{for all }x\in \mathbb{R}.  \label{2.6}
\end{align}%
By a straightforward computation, we deduce from (\ref{main eqn}), %
(\ref{2.4444}), (\ref{2.5}) and (\ref{2.6}) that 
\begin{equation*}
\left\{ 
\begin{array}{rcll}
v_{tt}(x,t) & = & v_{xx}(x,t)+r(x)v(x,t) & x\in \mathbb{R},t\in (0,\infty ),
\\ 
v(x,0) & = & 0 & x\in \mathbb{R}, \\ 
v_{t}(x,0) & = & \delta _{0}(x) & x\in \mathbb{R}.%
\end{array}%
\right.
\end{equation*}%
Hence, $v$ has the form, see \cite[Chapter 2, §3]{Romanov:VNU1986}, 
\begin{equation}
v(x,t)=\frac{H(t-|x|)}{2}+\frac{1}{2}\int_{(x-t)/2}^{(x+t)/2}r(\xi )\Big(%
\int_{|\xi |}^{t-|x-\xi |}v(\xi ,\tau )d\tau \Big)d\xi .  \label{2.7}
\end{equation}%
Using (\ref{2.4444}) and (\ref{2.7}, we obtain for all $x\in \mathbb{R}$, $%
t\in (0,\infty )$ 
\begin{equation*}
u(x,t)c^{1/4}(x)=v(\tau (x),t)=\frac{H(t-|\tau (x)|)}{2}+\frac{1}{2}%
\int_{(\tau (x)-t)/2}^{(\tau (x)+t)/2}r(\xi )\Big(\int_{|\xi |}^{t-|\tau
(x)-\xi |}v(\xi ,\tau )d\tau \Big)d\xi .
\end{equation*}%
%
%\begin{equation*}
%u(x,t)c^{1/4}(x)=v(\tau (x),t)
%\end{equation*}%
%\begin{equation}
%=\frac{H(t-|\tau (x)|)}{2}+\frac{1}{2}\int_{(\tau (x)-t)/2}^{(\tau
%(x)+t)/2}r(\xi )\Big(\int_{|\xi |}^{t-|\tau (x)-\xi |}v(\xi ,\tau )d\tau %
%\Big)d\xi .  \label{2.70}
%\end{equation}%
We obtain (\ref{2.3}). Formula (\ref{u2c}) is a direct consequence of %
(\ref{2.3}). \qed

We also need Lemma 2.2 below, which is a direct consequence of one of
results of \cite[Lemma 2.2]{SmirnovKlibanovNguyen:IPI2020}. Using this
lemma, we establish a boundary condition at $x=M$ where we cannot measure
any information of the function $u$.

\textbf{Lemma 2.2.} (Absorbing boundary condition) \emph{Assume that }$%
M>\max \{c(x):x\in R\}$\emph{. Then }%
\begin{equation}
u_{x}(M,t)+u_{t}(M,t)=0,\quad \mbox{for all }t>0  \label{abso}
\end{equation}%
\emph{and }%
\begin{equation}
u_{x}(x,t)-u_{t}(x,t)=0,\quad \mbox{for all }x<0,t>0.  \label{abso1}
\end{equation}

Let $\Omega _{T}=(\epsilon ,M)\times (0,T).$ Introduce the change of
variables, 
\begin{equation}
q(x,t)=u(x,t+\tau (x)),\quad \mbox{for all }(x,t)\in \Omega _{T}
\label{q def}
\end{equation}%
where $u$ is the solution of (\ref{main eqn}). We derive a partial
differential equation for the function $q$, in which the function $c$ is not
present. The first step is to show that $q(x,0)>0$ for all $x>\epsilon $. In
fact, using (\ref{q def}) with $t=0$ and (\ref{u2c}) with $t=\tau
(x)^{+}=\tau (x)$, we obtain 
\begin{equation}
q(x,0)=u(x,\tau (x))=\frac{1}{2c^{1/4}(x)}\geq \underline{q}:=\frac{1}{2%
\overline{c}^{1/4}}>0,\quad \mbox{for all }x>\epsilon  \label{2.14}
\end{equation}%
where $\overline{c}$ is a known upper bound of the function $c$, see %
(\ref{1.1}). The estimate (\ref{2.14}) is important because we will see
later in this section that $q(x,0)$ is the denominator of some components
for the desired differential equation that governs the function $q$. We next
differentiate (\ref{q def}) with respect to $x$ to obtain 
\begin{equation}
q_{x}(x,t)=u_{x}(x,t+\tau (x))+u_{t}(x,t+\tau (x))\tau ^{\prime }(x)
\label{qx}
\end{equation}%
for all $x>\epsilon $, $t>0$. Thus, for all $x>\epsilon $, $t>0$, we obtain 
\begin{equation}
q_{xt}(x,t)=u_{xt}(x,t+\tau (x))+u_{tt}(x,t+\tau (x))\tau ^{\prime }(x)
\label{2.2222}
\end{equation}%
and 
\begin{multline}
q_{xx}(x,t)=u_{xx}(x,t+\tau (x))+2u_{xt}(x,t+\tau (x))\tau ^{\prime
}(x)+u_{tt}(x,t+\tau (x))|\tau ^{\prime 2} \\
+u_{t}(x,t+\tau (x))\tau ^{\prime \prime }(x).  \label{2.2}
\end{multline}%
It follows from the governing equation (\ref{main eqn}), the \textit{eikonal}
equation (\ref{eik}), (\ref{2.2222}), and (\ref{2.2}) that 
\begin{align}
q_{xx}(x,t)& =2u_{xt}(x,t+\tau (x))\tau ^{\prime }(x)+2u_{tt}(x,t+\tau
(x))|\tau ^{\prime 2}+u_{t}(x,t+\tau (x))\tau ^{\prime \prime }(x)  \notag \\
& =2q_{xt}(x,t)\tau ^{\prime }(x)+u_{t}(x,t+\tau (x))\tau ^{\prime \prime
}(x)  \notag \\
& =2q_{xt}(x,t)\tau ^{\prime }(x)+q_{t}(x,t)\tau ^{\prime \prime }(x)
\label{2.4}
\end{align}%
for all $x>\epsilon ,t>0.$ We next eliminate $\tau $ from (\ref{2.4}). Using %
(\ref{tau}) and (\ref{2.14}), we have 
\begin{equation}
\tau ^{\prime }(\mathbf{x})=\sqrt{c(x)}=\frac{1}{4q^{2}(x,0)}.  \label{2.11}
\end{equation}%
Therefore, for all $x>\epsilon $ 
\begin{equation}
\tau ^{\prime \prime }(\mathbf{x})=-\frac{q_{x}(x,0)}{2q^{3}(x,0)}.
\label{2.12}
\end{equation}%
Combining (\ref{2.4}), (\ref{2.11}) and (\ref{2.12}), we arrive at the
following equation for the function $q$ 
\begin{equation}
q_{xx}(x,t)-\frac{q_{xt}(x,t)}{2q^{2}(x,0)}+\frac{q_{t}(x,t)q_{x}(x,0)}{%
2q^{3}(x,0)}=0,\quad \mbox{for all }(x,t)\in (\epsilon ,M)\times (0,T).
\label{q eqn}
\end{equation}%
We need to solve equation (\ref{q eqn}) for the function $q$. In the next
step, we find the boundary values $q(\epsilon ,t)$, $q_{x}(\epsilon ,t)$ and 
$q_{x}(M,t)$ for $t\in \lbrack 0,T].$ It follows from (\ref{q def}) that 
\begin{equation*}
q(\epsilon ,t)=u(\epsilon ,t+\tau (\epsilon ))=u(\epsilon ,t+\epsilon
),\quad \mbox{for all }t>0.
\end{equation*}%
Therefore, using (\ref{1.3}) and (\ref{tauepsilon}), we have 
\begin{equation*}
q(\epsilon ,t)=g_{0}(t+\epsilon ),\quad \mbox{for all }t>0.
\end{equation*}%
By (\ref{qx}), 
\begin{equation*}
q_{x}(\epsilon ,t)=u_{x}(\epsilon ,t+\tau (\epsilon ))+u_{t}(\epsilon
,t+\tau (\epsilon ))\tau ^{\prime }(\epsilon ),\quad \mbox{for all }t>0.
\end{equation*}%
By (\ref{2.444444}) $\tau ^{\prime }(\epsilon )=1$. This, together with %
(\ref{1.3}), implies 
\begin{equation*}
q_{x}(\epsilon ,t)=g_{1}(t+\epsilon )+g_{0}^{\prime }(t+\epsilon ),\quad %
\mbox{for all }t>0.
\end{equation*}%
On the other hand, assume that $M>\Vert c\Vert _{L^{\infty }(\mathbb{R})}.$
Due to (\ref{abso}) and the fact that by (\ref{1.1}) 
\begin{equation*}
\tau ^{\prime }(M)=\sqrt{c(M)}=1,
\end{equation*}%
we obtain for all $t\in \lbrack 0,T],$ 
\begin{equation*}
q_{x}(M,t)=u_{x}(M,t+\tau (M))+u_{t}(M,t+\tau (M))\tau ^{\prime }(M)=0.
\end{equation*}%
In summary, we have proved the following proposition.

\begin{proposition}
The function $q$ defined in (\ref{q def}) satisfies 
\begin{equation}
\left\{ 
\begin{array}{ll}
q_{xx}(x,t)-\frac{q_{xt}(x,t)}{2q^{2}(x,0)}+\frac{q_{t}(x,t)q_{x}(x,0)}{%
2q^{3}(x,0)}=0 & (x,t)\in (\epsilon ,M)\times (0,T), \\ 
q(x,0)\geq \underline{q}>0 & x\in \lbrack \epsilon ,M], \\ 
q(\epsilon ,t)=g_{0}(t+\epsilon ) & t\in \lbrack 0,T], \\ 
q_{x}(\epsilon ,t)=g_{1}(t+\epsilon )+g_{0}^{\prime }(t+\epsilon ) & t\in
\lbrack 0,T], \\ 
q_{x}(M,t)=0 & t\in \lbrack 0,T],%
\end{array}%
\right.   \label{q sys}
\end{equation}%
where the number $\underline{q}$ is defined in (\ref{2.14}).
\end{proposition}

\textbf{Remark 2.1}. \emph{Recall that in the statement of Problem \ref{cip}%
, we need the data }$g_{0}$\emph{\ and }$g_{1}$\emph{\ are known in }$%
[0,T+\epsilon ].$\emph{\ This is because we need the boundary conditions for
the function }$q$\emph{\ in (\ref{q sys}) to be well-defined.}

\textbf{Remark 2.2}.\emph{\ Solving numerically nonlinear partial
differential equations like the equation in (\ref{q sys}), in which the non
local term }$q(x,0)$\emph{\ is involved, is interesting not only in the area
of inverse problems but, more generally, in the area of scientific
computations.}

\textbf{Remark 2.3}. \emph{We use the space }$H^{5}(\Omega _{T})$\emph{\
below because, by embedding theorem, }%
\begin{equation}
H^{5}(\Omega _{T})\subset C^{3}\left( \overline{\Omega _{T}}\right) ,
\label{2}
\end{equation}%
\emph{and }%
\begin{equation}
\left\Vert f\right\Vert _{C^{3}\left( \overline{\Omega }_{T}\right) }\leq
C_{1}\left[ f\right] ,\quad \mbox{for all }f\in H^{5}(\Omega _{T}).
\label{3}
\end{equation}%
\emph{This helps us to prove the convexification theorem. Here the constant }%
$C_{1}>0$\emph{\ depends only on the domain }$\Omega _{T}.$\emph{\ Below }$%
\left[ \cdot ,\cdot \right] $\emph{\ and }$\left[ \cdot \right] $\emph{\
denote the scalar product and the norm respectively in the space }$%
H^{5}(\Omega _{T})$\emph{\ of real valued functions.}

Problem \ref{cip} is reduced to the problem of computing the function $q$
that satisfies (\ref{q sys}). In fact, having $q$, we can compute $c$ via
the formula 
\begin{equation}
c(x)=\frac{1}{16q^{4}(x,0)},\quad \mbox{for all }x\in \lbrack \epsilon ,M].
\label{q2c}
\end{equation}

We solve (\ref{q sys}) for $q$ using the convexification method. By
convexification, we mean that we use the Carleman weigh function $%
e^{-\lambda (x+\alpha t)}$ to convexify the mismatch functional for some
suitably chosen parameters $\lambda $ and $\alpha $. Let the operator $%
F:H\rightarrow L^{2}(\Omega _{T})$ be given by 
\begin{equation}
F(q)=q_{xx}-\frac{q_{xt}}{2q^{2}(x,0)}+\frac{q_{t}q_{x}(x,0)}{2q^{3}(x,0)}%
,\quad \mbox{for all }q\in H^{5}(\Omega _{T}).  \label{3.4}
\end{equation}%
We define the weighted Tikhonov-like functional $J_{\lambda ,\alpha ,\beta }$
as: 
\begin{equation}
J_{\lambda ,\alpha ,\beta }(q)=\int_{\Omega }e^{-2\lambda (x+\alpha t)}\big|%
F(q)\big|^{2}dxdt+\beta \left[ q\right] ^{2}.  \label{2.25}
\end{equation}%
We claim, given a certain closed convex set in the space $H^{5}(\Omega _{T})$
of an arbitrary diameter $d>0,$ there exist constants $\lambda _{0}\geq 1$
and $\alpha _{0}>0,$ depending on $d$ and some other parameters, such that
whenever $\lambda \geq \lambda _{0}$, $0<\alpha <\alpha _{0}$ and $\beta \in %
\left[ 2e^{-2\lambda \alpha T},1\right) ,$ functional (\ref{2.25}) is
strictly convex on that set. Furthermore, $J_{\lambda ,\alpha ,\beta }(q)$
has the unique minimizer on that set. We define that set below. Here, $\beta
\Vert q\Vert _{H^{4}(\Omega _{T})}^{2}$ is the regularization term and $%
\beta $ is the regularization parameter. This claim is one of the main
results in this paper, see Section \ref{sec conv}. Other important results
behind our numerical algorithm of Section \ref{sec conv}, include:

\begin{enumerate}
\item The confirmation that the well-known gradient descent method delivers
a sequence converging to the unique minimizer of $J_{\lambda ,\alpha ,\beta
} $ if starting at an arbitrary point of the above mentioned set.

\item The convergence of the minimizers of $J_{\lambda ,\alpha ,\beta }$ to
the true solution of (\ref{q sys}) as the noise contained in the measured
data tends to zero.
\end{enumerate}

Thus, this is \emph{global convergence: }see items 1 and 2 in Section 1.
Those important results are proved based on a new Carleman estimate of the
next section.

\section{Two New Carleman Estimates}

\label{sec Car}

Let the function $m\left( x\right) \in C^{1}([\epsilon ,M])$. For two
numbers $m_{0},m_{1},0<m_{0}<m_{1},$ we assume that 
\begin{equation}
0<m_{0}<m(x)<m_{1},\quad \mbox{ for all }x\in \lbrack \epsilon ,M]
\label{mbound}
\end{equation}%
Let%
\begin{equation}
m_{2}=\max_{[\epsilon ,M]}\left\vert m^{\prime }\left( x\right) \right\vert .
\label{1}
\end{equation}

\bigskip For $v\in H^{2}(\Omega _{T}),$ define the operator 
\begin{equation*}
L_{0}v(x,t)=v_{xx}(x,t)-m(x)v_{xt}(x,t),\quad \mbox{for all }(x,t)\in \Omega
_{T}.
\end{equation*}

\textbf{Theorem 3.1}. \emph{1. There exist a number }$\alpha _{0}>0$\emph{\
and a sufficiently large number }$\lambda _{0}\geq 1,$\emph{\ both depending
only on }$m_{0},$\emph{\ }$m_{1}$\emph{, }$m_{2},\epsilon ,$\emph{\ }$T,$%
\emph{\ such that for all }$\alpha \in (0,\alpha _{0}),$\emph{\ all }$%
\lambda \geq \lambda _{0}$\emph{\ and for all functions }$v\in H^{2}(\Omega
_{T})$\emph{\ the following Carleman estimate is valid: }%
\begin{multline}
\int_{\Omega _{T}}e^{-2\lambda (x+\alpha t)}|L_{0}v|^{2}dxdt\geq C\lambda
\int_{\Omega _{T}}e^{-2\lambda (x+\alpha t)}\big(|v_{x}|^{2}+|v_{t}|^{2}+%
\lambda ^{2}|v|^{2}\big)dtdx \\
-\int_{0}^{T}C\lambda e^{-2\lambda (\epsilon +\alpha t)}\big(|v_{x}(\epsilon
,t)|^{2}+\lambda ^{2}|v(\epsilon ,t)|^{2}+|v_{t}(\epsilon ,t)|^{2}\big)dt \\
+C\lambda \int_{\epsilon }^{M}e^{-2\lambda x}\big(|v_{x}(x,0)|^{2}+\lambda
^{2}|v(x,0)|^{2}\big)dx \\
-C\lambda \int_{\epsilon }^{M}e^{-2\lambda (x+\alpha T)}\big(%
|v_{x}(x,T)|^{2}+\lambda ^{2}|v(x,T)|^{2}\big)dx,
\label{100}
\end{multline}%
\emph{where }$C=C\left( m_{0},m_{1},m_{2},\epsilon ,M\right) >0$\emph{\ is a
generic constant depending only on listed parameters.}

\emph{2. Let the function }$g\in C^{1}(\overline{\Omega _{T}}).$\emph{\
There exist a number }$\alpha _{0}>0$\emph{\ and a sufficiently large number 
}$\lambda _{0}\geq 1,$\emph{\ both depending only on }$m_{0},$\emph{\ }$%
m_{1} $\emph{, }$m_{2},\epsilon ,$\emph{\ }$M,\Vert g\Vert _{C^{1}(\overline{%
\Omega _{T}})}$\emph{, }$T,$\emph{\ such that for all }$\alpha \in (0,\alpha
_{0}),$\emph{\ all }$\lambda \geq \lambda _{0}$\emph{\ and for all functions 
}$v\in H^{2}(\Omega _{T})$\emph{\ the following Carleman estimate is valid:}%
\begin{multline}
\int_{\Omega _{T}}e^{-2\lambda (x+\alpha t)}|L_{0}v|^{2}dxdt+\int_{\Omega
_{T}}e^{-2\lambda (x+\alpha t)}gv_{xt}\left( x,t\right) v\left( x,0\right)
dxdt \\
\geq C\lambda \int_{\Omega _{T}}e^{-2\lambda (x+\alpha t)}\big(%
|v_{x}|^{2}+|v_{t}|^{2}+\lambda ^{2}|v|^{2}\big)dtdx-\int_{0}^{T}C\lambda
e^{-2\lambda (\epsilon +\alpha t)}\big(|v_{x}(\epsilon ,t)|^{2}+\lambda
^{2}|v(\epsilon ,t)|^{2}+|v_{t}(\epsilon ,t)|^{2}\big)dt \\
+C\lambda \int_{\epsilon }^{M}e^{-2\lambda x}\big(|v_{x}(x,0)|^{2}+\lambda
^{2}|v(x,0)|^{2}\big)dx-C\lambda \int_{\epsilon }^{M}e^{-2\lambda (x+\alpha
T)}\big(|v_{x}(x,T)|^{2}+\lambda ^{2}|v(x,T)|^{2}\big)dx,  \label{101}
\end{multline}%
\emph{where }$C=C\left( m_{0},m_{1},m_{2},\epsilon ,M,\left\Vert
g\right\Vert _{C^{1}\left( \overline{\Omega }_{T}\right) }\right) >0$\emph{\
is a generic constant depending only on listed parameters. }

\textbf{Remark 3.1}. \emph{Theorem 3.1 is a significant
generalization of the Carleman estimate  in \cite%
{SmirnovKlibanovNguyen:IPI2020}. In fact, the function }$m$\emph{\ in that
theorem is a constant. On the other hand, the function }$m(x)$\emph{\ in
this paper is }$m(x)=1/\left( 2q^{2}(x,0)\right) $\emph{\ for }$x\in \lbrack
\epsilon ,M].$

\textbf{Remark 3.2}. \emph{An unusual element of the second Carleman
estimate (\ref{101}) of Theorem 3.1 is the presence of the nonlinear term }$%
v_{xt}\left( x,t\right) v\left( x,0\right) ,$\emph{\ which contains the
derivative }$v_{xt}\left( x,t\right) $\emph{\ involved in the operator }$%
L_{0}$\emph{\ as well as the non local term }$v\left( x,0\right) .$\emph{\
To the best of our knowledge a similar term was involved only in the
Carleman estimate of the paper \cite{BeilinaKlibanov:narwa2015}. However,
that term was different from the one in (\ref{101}) and its treatment was
different as well.}

\textbf{Remark 3.3}. \emph{For brevity, we treat below the constant }$C>0$%
\emph{\ of both parts of Theorem 3.1 as one depending on }$%
m_{0},m_{1},m_{2},\epsilon ,M,\left\Vert g\right\Vert _{C^{1}\left( 
\overline{\Omega }_{T}\right) },T:$\emph{\ as in the second part of this
theorem.}

\noindent \textbf{Proof of Theorem 3.1.}\textit{\ }We prove this theorem
only for functions\textit{\ }$v\in C^{2}\left( \overline{\Omega }_{T}\right)
.$ The case $v\in H^{2}\left( \Omega _{T}\right) $ follows from density
arguments. We split the proof in several steps. First, we prove (\ref{100}).

\noindent \textit{Step 1.} Define the function $w$, 
\begin{equation}
w(x,t)=e^{-\lambda (x+\alpha t)}v(x,t), \quad \mbox{for all }(x,t)\in \Omega
_{T}.  \label{4.2222}
\end{equation}%
We have in $\Omega _{T}$, 
\begin{equation}
\begin{array}{rcl}
v & = & e^{\lambda (x+\alpha t)}w, \\ 
v_{x} & = & e^{\lambda (x+\alpha t)}\big(w_{x}+\lambda w\big), \\ 
v_{xx} & = & e^{\lambda (x+\alpha t)}\big(w_{xx}+2\lambda w_{x}(x,t)+\lambda
^{2}w\big), \\ 
v_{xt} & = & e^{\lambda (x+\alpha t)}\big(w_{xt}+\lambda \alpha
w_{x}+\lambda w_{t}+\lambda ^{2}\alpha w\big).%
\end{array}%
\end{equation}%
Therefore, 
\begin{equation}
e^{-\lambda (x+\alpha t)}L_{0}v=w_{xx}-mw_{xt}+\lambda \big(2-\alpha m\big)%
w_{x}-\lambda mw_{t}+\lambda ^{2}\big(1-\alpha m\big)w.  \label{4.2}
\end{equation}%
Let $s$ be a positive number, which we will choose later. It follows from %
(\ref{4.2}) that 
\begin{equation}
e^{-sx}e^{-2\lambda (x+\alpha t)}[L_{0}v]^{2}=e^{-sx}\Big[\big(%
w_{xx}-mw_{xt}+\lambda ^{2}\big(1-\alpha m\big)w\big)+\left( \lambda \big(%
2-\alpha m\big)w_{x}-\lambda mw_{t}\right) \Big].  \label{4.3}
\end{equation}%
Since $\left( a+b\right) ^{2}\geq 2ab,\forall a,b\in \mathbb{R},$ then 
\begin{equation*}
|L_{0}v|^{2}\geq \left( \lambda \big(2-\alpha m\big)w_{x}-\lambda
mw_{t}\right) \big(w_{xx}-mw_{xt}+\lambda ^{2}\big(1-\alpha m\big)w\big),
\end{equation*}%
Hence, 
\begin{equation}
e^{-sx}e^{-2\lambda (x+\alpha t)}|L_{0}v|^{2}\geq I_{1}+I_{2},  \label{4.55}
\end{equation}%
where 
\begin{align}
I_{1}& =2\lambda e^{-sx}\big(2-\alpha m\big)w_{x}\Big[w_{xx}-mw_{xt}+\lambda
^{2}\big(1-\alpha m\big)w\Big],  \label{4.5} \\
I_{2}& =-2\lambda e^{-sx}mw_{t}\Big[w_{xx}-mw_{xt}+\lambda ^{2}\big(1-\alpha
m\big)w\Big].  \label{4.6}
\end{align}

\noindent \textit{Step 2.} In this step, we estimate $I_{1}.$ By (\ref{4.5}), 
\begin{equation*}
I_{1}=2\lambda \big(2-\alpha m\big)e^{-sx}w_{x}w_{xx}-2\lambda m\big(%
2-\alpha m\big)e^{-sx}w_{x}w_{xt}+2\lambda ^{3}\big(2-\alpha m\big)\big(%
1-\alpha m\big)e^{-sx}w_{x}w
\end{equation*}%
\begin{equation*}
=\lambda \left( 2-\alpha m\right) e^{-sx}\left( |w_{x}|^{2}\right)
_{x}-\lambda m\left( 2-\alpha m\right) e^{-sx}\left( |w_{x}|^{2}\right)
_{t}+\lambda ^{3}\left( 2-\alpha m\right) \left( 1-\alpha m\right)
e^{-sx}\left( |w|^{2}\right) _{x}.
\end{equation*}%
Thus,%
\begin{equation*}
I_{1}=\lambda \left( \left( 2-\alpha m\right) e^{-sx}|w_{x}|^{2}\right)
_{x}+\lambda \left( \left( s\left( 2-\alpha m\right) +\alpha m^{\prime
}\right) e^{-sx}\right) |w_{x}|^{2}+\left( -\lambda m\left( 2-\alpha
m\right) e^{-sx}|w_{x}|^{2}\right) _{t}
\end{equation*}%
\begin{equation*}
+\left( \lambda ^{3}\left( 2-\alpha m\right) \left( 1-\alpha m\right)
e^{-sx}|w|^{2}\right) _{x}+\lambda ^{3}\left( \left( s\left( 2-\alpha
m\right) \left( 1-\alpha m\right) -3\alpha m^{\prime }+2\alpha
^{2}mm^{\prime }\right) e^{-sx}\right) |w|^{2}.
\end{equation*}
This formula is equivalent to 
\begin{multline}
I_{1}=\lambda \left( \left( s\left( 2-\alpha m\right) +\alpha m^{\prime
}\right) e^{-sx}\right) |w_{x}|^{2} +\lambda^{3}\left( \left( s\left(
2-\alpha m\right) \left( 1-\alpha m\right) -3\alpha m^{\prime }+2\alpha
^{2}mm^{\prime }\right) e^{-sx}\right) |w|^{2} \\
+\left( \lambda \left( 2-\alpha m\right) e^{-sx}|w_{x}|^{2} +\lambda
^{3}\left( 2-\alpha m\right) \left( 1-\alpha m\right) e^{-sx}|w|^{2}\right)
_{x} \\
+\left( -\lambda m\left( 2-\alpha m\right) e^{-sx}|w_{x}|^{2}\right) _{t}.
\label{4.7}
\end{multline}

\noindent \textit{Step 3.} We now estimate $I_{2}$. By (\ref{4.6}), we have%
\begin{align*}
I_{2}&=-2\lambda e^{-sx}mw_{t}\left[ w_{xx}-mw_{xt}+\lambda ^{2}\left(
1-\alpha m\right) w\right] \\
&=\left( -2\lambda e^{-sx}mw_{t}w_{x}\right) _{x}+2\lambda
e^{-sx}mw_{xt}w_{x}-2\lambda \left( sm-m^{\prime }\right) w_{t}w_{x}+\left(
\lambda e^{-sx}m^{2}w_{t}^{2}\right) _{x} \\
&\hspace{6cm}+\lambda e^{-sx}\left( s-2mm^{\prime }\right) w_{t}^{2}+\left(
-\lambda ^{3}e^{-sx}m\left( 1-\alpha m\right) w^{2}\right) _{t} \\
&=\lambda e^{-sx}\left( sm^{2}-2mm^{\prime }\right) w_{t}^{2}-2\lambda
\left( sm-m^{\prime }\right) w_{t}w_{x} +\left( \lambda
e^{-sx}mw_{x}^{2}-\lambda ^{3}e^{-sx}m\left( 1-\alpha m\right) w^{2}\right)
_{t} \\
&\hspace{8cm} +\left( \lambda e^{-sx}m^{2}w_{t}^{2}-2\lambda
e^{-sx}mw_{t}w_{x}\right) _{x}.
\end{align*}%
Thus,%
\begin{multline}
I_{2}=\lambda e^{-sx}\left( s-2mm^{\prime }\right) w_{t}^{2}-2\lambda \left(
sm-m^{\prime }\right) w_{t}w_{x} +\left( \lambda e^{-sx}mw_{x}^{2}-\lambda
^{3}e^{-sx}m\left( 1-\alpha m\right) w^{2}\right) _{t} \\
+\left( \lambda e^{-sx}m^{2}w_{t}^{2}-2\lambda e^{-sx}mw_{t}w_{x}+\lambda
e^{-sx}m^{2}w_{t}^{2}\right) _{x}.  \label{4.8}
\end{multline}

\noindent \textit{Step 4.} In this step, we estimate $I_{1}+I_{2}$ as $%
s\rightarrow \infty $. Below, the notation $O(1/s)$ indicates the quantity
satisfying 
\begin{equation*}
|O(1/s)|\leq \frac{C}{s}\quad \mbox{as }s\rightarrow \infty ,
\end{equation*}%
where $C>0$ is independent on $x,t,s$. Adding (\ref{4.7}) and (\ref{4.8}),
we obtain 
\begin{multline}
I_{1}+I_{2}\geq \lambda s\Big((2-\alpha
m)|w_{x}|^{2}-2mw_{t}w_{x}+m^{2}|w_{t}|^{2}\Big)e^{-sx}-C\lambda \left(
|w_{x}|^{2}+|w_{t}|^{2}\right) e^{-sx} \\
+\lambda ^{3}s\Big((2-\alpha m)(1-\alpha m)+O(1/s)\Big)e^{-sx}|w|^{2}+\Big(%
\lambda \big(2-\alpha m\big)e^{-sx}|w_{x}|^{2}+\lambda ^{3}\big(2-\alpha m%
\big)\big(1-\alpha m\big)e^{-sx}|w|^{2}\Big)_{x} \\
+\Big(-2\lambda e^{-sx}mw_{t}w_{x}+\lambda e^{-sx}m^{2}|w_{t}|^{2}\Big)_{x}+%
\Big(-\lambda m(1-\alpha m)(|w_{x}|^{2}+\lambda ^{2}|w|^{2})e^{-sx}\Big)_{t}.
\label{4.9}
\end{multline}%
We estimate the the first term in the right hand side of (\ref{4.9}). Using
the inequality 
\begin{equation*}
-2mw_{t}w_{x}\geq -\frac{3}{4}m^{2}|w_{t}|^{2}-\frac{4}{3}|w_{x}|^{2},
\end{equation*}%
we obtain 
\begin{equation}
\lambda s\Big((2-\alpha m)|w_{x}|^{2}-2mw_{t}w_{x}+m^{2}|w_{t}|^{2}\Big)%
e^{-sx}\geq \lambda s\Big[\Big(\frac{2}{3}-\alpha m\Big)|w_{x}|^{2}+\frac{1}{%
4}m^{2}|w_{t}|^{2}\Big].  \label{4.11}
\end{equation}%
Let $\alpha _{0}=2/\left( 3m_{1}\right) ,$ where $m_{1}$ and $m_{2}$ are
defined in (\ref{mbound}). Then, it follows from (\ref{4.11}) that for all $%
\alpha \in (0,\alpha _{0})$ 
\begin{equation}
\lambda s\Big((2-\alpha m)|w_{x}|^{2}-2mw_{t}w_{x}+m^{2}|w_{t}|^{2}\Big)%
e^{-sx}\geq C\lambda s\big(|w_{x}|^{2}+|w_{t}|^{2}\big)e^{-sx}.  \label{4.12}
\end{equation}%
Hence, by (\ref{4.9}) and (\ref{4.12}), we can find a number $%
s_{0}=s_{0}(m_{1},m_{2})>0$ such that for all $s\geq s_{0}$ 
\begin{equation}
I_{1}+I_{2}\geq C\lambda s\big(|w_{x}|^{2}+|w_{t}|^{2}\big)e^{-sx}+C\lambda
^{3}s|w|^{2}e^{-sx}+U_{x}+V_{t},  \label{4.13}
\end{equation}%
where 
\begin{equation*}
U=\lambda \big(2-\alpha m\big)e^{-sx}|w_{x}|^{2}+\lambda ^{3}\big(2-\alpha m%
\big)\big(1-\alpha m\big)e^{-sx}|w|^{2}-2\lambda e^{-sx}mw_{t}w_{x}+\lambda
e^{-sx}m^{2}|w_{t}|^{2}
\end{equation*}%
and 
\begin{equation*}
V=-\lambda m(1-\alpha m)(|w_{x}|^{2}+\lambda ^{2}|w|^{2})e^{-sx}.
\end{equation*}

\noindent \textit{Step 5.} In this step, we estimate the integrals of $U_{x}$
and $V_{t}$. It is obvious that 
\begin{equation*}
U\geq C\lambda \big(|w_{x}|^{2}+\lambda ^{2}|w|^{2}+|w_{t}|^{2}\big).
\end{equation*}%
and the opposite inequality is also true with a different constant $C$.
Hence, 
\begin{equation}
\int_{\Omega _{T}}U_{x}dxdt\geq -\int_{0}^{T}C\lambda \big(|w_{x}(\epsilon
,t)|^{2}+\lambda ^{2}|w(\epsilon ,t)|^{2}+|w_{t}(\epsilon ,t)|^{2}\big)%
e^{-s\epsilon }dt.  \label{4.16}
\end{equation}

On the other hand, recalling numbers $m_{0}$ and $m_{1}$ in (\ref{mbound})
and that $\alpha \in \left( 0,2/\left( 3m_{1}\right) \right) $, we obtain 
\begin{equation*}
\lambda m(1-\alpha m)\geq \lambda m_{0}(1-\alpha m_{1})\geq C\lambda ,\quad %
\mbox{for all }\alpha \in (0,\alpha _{0}).
\end{equation*}%
Hence,%
\begin{align}
\int_{\Omega _{T}}V_{t}dtdx&=\int_{\epsilon }^{M}\left( -V\left( x,0\right)
+V\left( x,T\right) \right) dx  \notag \\
&\geq C\lambda \int_{\epsilon }^{M}\big(|w_{x}(x,0)|^{2}+\lambda
^{2}|w(x,0)|^{2}\big)e^{-sx}dx  \notag \\
&\hspace{4cm} -C\lambda \int_{\epsilon }^{M}\big(|w_{x}(x,T)|^{2}+\lambda
^{2}|w(x,T)|^{2}\big)e^{-sx}dx.  \label{4.17}
\end{align}

\noindent \textit{Step 6.} Combining (\ref{4.55}), (\ref{4.13}), (\ref{4.16})
and (\ref{4.17}), setting $s=s_{0}$ and regarding $e^{-s_{0}\left(
M-\epsilon \right) }$ as a part of the constant $C$, we obtain 
\begin{multline}
\int_{\Omega _{T}}e^{-2\lambda (x+\alpha t)}[L_{0}v(x,t)]^{2}dxdt\geq
C\lambda \int_{\Omega _{T}}\big(|w_{x}(x,t)|^{2}+|w_{t}(x,t)|^{2}+\lambda
^{2}|w(x,t)|^{2}\big)dtdx \\
-\int_{0}^{T}C\lambda \big(|w_{x}(\epsilon ,t)|^{2}+\lambda ^{2}|w(\epsilon
,t)|^{2}+|w_{t}(\epsilon ,t)|^{2}\big)dt+C\lambda \int_{\epsilon }^{M}\big(%
|w_{x}(x,0)|^{2}+\lambda ^{2}|w(x,0)|^{2}\big)dx \\
-C\lambda \int_{\epsilon }^{M}\big(|w_{x}(x,T)|^{2}+\lambda ^{2}|w(x,T)|^{2}%
\big)dx.  \label{4.19}
\end{multline}%
It follows from (\ref{4.2222}) that $w=e^{-\lambda (x+\alpha t)}v$. Hence, 
\begin{align}
|w_{x}|^{2}+\frac{1}{2}\lambda ^{2}|w|^{2}& =e^{-2\lambda (x+\alpha t)}\big|%
-\lambda v+v_{x}\big|^{2}+\frac{1}{2}\lambda ^{2}|w|^{2}  \notag \\
& =e^{-2\lambda (x+\alpha t)}\big(\frac{3}{2}\lambda ^{2}|v|^{2}-2\lambda
vv_{x}+|v_{x}|^{2}\big)\geq \Big(\frac{1}{6}\lambda ^{2}|v|^{2}+\frac{1}{4}%
|v_{x}|^{2}\Big)  \label{4.20}
\end{align}%
and 
\begin{align}
|w_{t}|^{2}+\frac{1}{2}\lambda ^{2}|w|^{2}& =e^{-2\lambda (x+\alpha t)}\big|%
-\lambda \alpha v+v_{t}\big|^{2}+\frac{1}{2}\lambda ^{2}|w|^{2}  \notag \\
& \geq e^{-2\lambda (x+\alpha t)}\Big(\big(\alpha ^{2}+\frac{1}{2}\big)%
\lambda ^{2}|v|^{2}-2\lambda \alpha vv_{t}+|v_{t}|^{2}\Big)  \notag \\
& \geq e^{-2\lambda (x+\alpha t)}\Big(\big(\alpha ^{2}+\frac{1}{2}\big)%
\lambda ^{2}|v|^{2}-\frac{1+4\alpha ^{2}}{4\alpha ^{2}}\lambda ^{2}\alpha
^{2}|v|^{2}-\frac{4\alpha ^{2}}{1+4\alpha ^{2}}|v_{t}|^{2}+|v_{t}|^{2}\Big).
\notag \\
& =e^{-2\lambda (x+\alpha t)}\Big(\frac{1}{4}\lambda ^{2}|v|^{2}+\frac{1}{%
1+4\alpha ^{2}}|v_{t}|^{2}\Big).  \label{4.21}
\end{align}%
Adding (\ref{4.20}) and (\ref{4.21}), we have 
\begin{equation}
|w_{x}|^{2}+|w_{t}|^{2}+\lambda ^{2}|w|^{2}\geq e^{-2\lambda (x+\alpha
t)}C(|v_{x}|^{2}+|v_{t}|^{2}+\lambda ^{2}|v|^{2}).  \label{4.22}
\end{equation}

On the other hand, for all $s \geq s_0$ and $(x, t) \in [\epsilon, M] \times
[0, T],$ 
\begin{equation}
|w_x|^2 = e^{-2\lambda (x + \alpha t)}\big|-\lambda v + v_x\big|^2 \leq
e^{-2\lambda (x + \alpha t)}\big(2\lambda^2 |v|^2 + 2|v_x|^2\big)
\label{4.23}
\end{equation}
and 
\begin{equation}
|w_t|^2 = e^{-2\lambda (x + \alpha t)}\big|-\lambda \alpha v + v_t\big|^2
\leq e^{-2\lambda (x + \alpha t)}\big(2\lambda^2 \alpha^2 |v|^2 + 2|v_t|^2%
\big).  \label{4.24}
\end{equation}

Combining (\ref{4.19}), (\ref{4.22}), (\ref{4.23}) and (\ref{4.24}), we
obtain \ref{100}.

\emph{Step 7}. We now prove the second Carleman estimate (\ref{101}).
Consider the expression which is a part of the second term in the first line
of (\ref{101})
\begin{align*}
& e^{-2\lambda (x+\alpha t)}gv_{xt}\left( x,t\right) v\left( x,0\right) \\
& =v\left( e^{-2\lambda (x+\alpha t)}gv_{x}\left( x,t\right) v\left(
x,0\right) \right) _{t}-g_{t}v_{x}\left( x,t\right) v\left( x,0\right)
e^{-2\lambda (x+\alpha t)}-2\lambda \alpha gv_{x}\left( x,t\right) v\left(
x,0\right) e^{-2\lambda (x+\alpha t)} \\
& \geq \left( e^{-2\lambda (x+\alpha t)}gv_{x}\left( x,t\right) v\left(
x,0\right) \right) _{t}-C\left\vert v_{x}\left( x,t\right) \right\vert
^{2}e^{-2\lambda (x+\alpha t)}-C\lambda ^{2}v^{2}\left( x,0\right)
e^{-2\lambda (x+\alpha t)}.
\end{align*}%
Hence, 
\begin{multline}
\int_{\Omega _{T}}e^{-2\lambda (x+\alpha t)}gv_{xt}\left( x,t\right) v\left(
x,0\right) dxdt\geq -C\int_{\Omega }\left[ \left\vert v\left( x,T\right)
\right\vert ^{2}+\left\vert v\left( x,0\right) \right\vert ^{2}\right]
e^{-2\lambda (x+\alpha T)}dx \\
-C\int_{\Omega }\left[ \left\vert v_{x}\left( x,0\right) \right\vert
^{2}+\left\vert v\left( x,0\right) \right\vert ^{2}\right] e^{-2\lambda
x}dx-C\int_{\Omega _{T}}\left\vert v_{x}\left( x,t\right) \right\vert
^{2}e^{-2\lambda (x+\alpha t)}dx \\
-C\lambda ^{2}\int_{\Omega }\left\vert v\left( x,0\right) \right\vert
^{2}e^{-2\lambda x}\left( \int_{0}^{T}e^{-2\lambda \alpha t}dt\right) dx.
\label{240}
\end{multline}%
Since 
\begin{equation*}
\int_{0}^{T}e^{-2\lambda \alpha t}dt=\frac{1}{2\lambda \alpha }\left(
1-e^{-2\lambda \alpha T}\right) \leq \frac{1}{2\lambda \alpha },
\end{equation*}%
(\ref{240}) implies%
\begin{multline}
\int_{\Omega _{T}}e^{-2\lambda (x+\alpha t)}gv_{xt}\left( x,t\right) v\left(
x,0\right) dxdt\geq -C\int_{\Omega }\left\vert v\left( x,T\right)
\right\vert ^{2}e^{-2\lambda (x+\alpha T)}dx \\
-C\int_{\Omega }\left[ \left\vert v_{x}\left( x,0\right) \right\vert
^{2}+\lambda \left\vert v\left( x,0\right) \right\vert ^{2}\right]
e^{-2\lambda x}dx-C\int_{\Omega _{T}}\left\vert v_{x}\left( x,t\right)
\right\vert ^{2}e^{-2\lambda (x+\alpha t)}dx.  \label{241}
\end{multline}%
Summing up (\ref{241}) with (\ref{100}), we obtain%
\begin{multline}
\int_{\Omega _{T}}e^{-2\lambda (x+\alpha t)}|L_{0}v|^{2}dxdt+\int_{\Omega
_{T}}e^{-2\lambda (x+\alpha t)}gv_{xt}\left( x,t\right) v\left( x,0\right)
dxdt \\
\geq C\int_{\Omega _{T}}e^{-2\lambda (x+\alpha t)}\big(\lambda \left( 1-%
\frac{1}{\lambda }\right) |v_{x}|^{2}+\lambda |v_{t}|^{2}+\lambda ^{3}|v|^{2}%
\big)dtdx \\
-C\lambda \int_{0}^{T}e^{-2\lambda (\epsilon +\alpha t)}\big(|v_{x}(\epsilon
,t)|^{2}+|v_{t}(\epsilon ,t)|^{2}+\lambda ^{2}|v(\epsilon ,t)|^{2}\big)dt \\
+C\int_{\epsilon }^{M}e^{-2\lambda x}\big(\lambda \left( 1-\frac{1}{\lambda }%
\right) |v_{x}(x,0)|^{2}+\lambda ^{3}\left( 1-\frac{1}{\lambda ^{2}}\right)
|v(x,0)|^{2}\big)dx \\
-C\lambda \int_{\epsilon }^{M}e^{-2\lambda (x+\alpha T)}\big(%
|v_{x}(x,T)|^{2}+\lambda ^{2}|v(x,T)|^{2}\big)dx.  \label{242}
\end{multline}%
Since $1-1/\lambda \geq 1/2$ for $\lambda \geq 2,$ (\ref{242}) implies %
(\ref{101}) for sufficiently large values of $\lambda \geq \lambda _{0}.$ 
\qed

\section{Some Results of Convex Analysis}

\label{sec 4}

We need results of this section for the proof of the existence and
uniqueness of the minimizer of our convexification functional $J_{\lambda
,\alpha ,\beta }$ on the closed convex set $\overline{B(R)}\cap H.$ Both
this functional and this set are introduced in the next section 5. A close
analog of Theorem 4.1 of this section is Theorem 2.1 of \cite%
{KlibanovNik:ra2017}. However, that theorem was proven only for the case
when the considered functional is strictly convex on a ball with the center
at $\left\{ 0\right\} $ in a Hilbert space. On the other hand, since the
closed convex set $\overline{B(R)}\cap H$ is not such a ball, then we need
an analog of that theorem for an arbitrary closed convex set in a Hilbert
space. This is exactly what is done in the current section (Theorem 4.1),
and our proof is similar with the proof of Theorem 2.1 of \cite%
{KlibanovNik:ra2017}.

Let $H$ be a Hilbert space of real valued functions. In this section, we
denote $\left\Vert \cdot \right\Vert $ and $\left\langle \cdot ,\cdot
\right\rangle $ \ respectively the norm and the scalar product $H$. Let $%
G\subset H$ be a closed convex set. Let $I:G\rightarrow \mathbb{R}$ be a
functional. We assume the existence of the Fr\'echet derivative $I^{\prime
}\left( x\right) ,\forall x\in G$ of the functional $I.$ The Fréchet
derivative $I^{\prime }\left( x\right) \in H$ at a point $x\in G$ is
understood as%
\begin{equation*}
I\left( y\right) -I\left( x\right) =\left\langle I^{\prime }\left( x\right)
,y-x\right\rangle +o\left( \left\Vert x-y\right\Vert \right) ,\left\Vert
x-y\right\Vert \rightarrow 0,y\in G,
\end{equation*}%
We denote the action of $I^{\prime }\left( x\right) $ on the vector $h\in H$
as $I^{\prime }\left( x\right) \left( h\right) ,$ where $I^{\prime }\left(
x\right) \left( h\right) =\left\langle I^{\prime }\left( x\right)
,h\right\rangle .$ We assume that $I^{\prime }\left( x\right) $ is Lipschitz
continuous, i.e.%
\begin{equation}
\left\Vert I^{\prime }\left( x\right) -I^{\prime }\left( y\right)
\right\Vert \leq D\left\Vert x-y\right\Vert ,\forall x,y\in G,  \label{8.1}
\end{equation}%
$D=const.>0$. We assume the strict convexity of the functional $I\left(
x\right) $ on the set $G,$%
\begin{equation}
I\left( y\right) -I\left( x\right) -I^{\prime }\left( x\right) \left(
y-x\right) \geq \varkappa \left\Vert x-y\right\Vert ^{2},\forall x,y\in G,
\label{800}
\end{equation}%
where $\varkappa =const.>0.$ We have along with (\ref{800}):%
\begin{equation}
I\left( x\right) -I\left( y\right) -I^{\prime }\left( y\right) \left(
x-y\right) \geq \varkappa \left\Vert x-y\right\Vert ^{2},\forall x,y\in G.
\label{801}
\end{equation}%
Summing up (\ref{800}) and (\ref{801}), we obtain 
\begin{equation}
\left( I^{\prime }\left( x\right) -I^{\prime }\left( y\right) ,x-y\right)
\geq 2\varkappa \left\Vert x-y\right\Vert ^{2},\forall x,y\in G.  \label{803}
\end{equation}

\textbf{Lemma 4.1 }\cite{KlibanovNik:ra2017}, \cite[Chapter 10, section 3]%
{Minoux}. \emph{Assume that conditions (\ref{8.1}) and (\ref{800}) are in
place. A point }$x_{\min }\in G$\emph{\ is a point of a relative minimum of
the functional }$I\left( x\right) $\emph{\ on the set }$G$ \emph{if and only
if the following inequality is true:} \emph{\ \ }%
\begin{equation}
\left( I^{\prime }\left( x_{\min }\right) ,x_{\min }-y\right) \leq 0,\text{ }%
\forall y\in G.  \label{8.3}
\end{equation}%
\emph{If a point of a relative minimum of the functional }$I\left( x\right) $%
\emph{\ on the set }$G$\emph{\ exists, then then this point is unique. Thus,
this point is the point of the global minimum of }$I\left( x\right) $ \emph{%
on the set }$G.$

Choose an arbitrary point $x\in H$. The point $x_{pr}$ is called projection
of the point $x$ on the set $G$ if%
\begin{equation*}
\left\Vert x-x_{pr}\right\Vert =\inf_{y\in G}\left\Vert x-y\right\Vert \text{
and }x_{pr}\in G.
\end{equation*}

\textbf{Lemma 4.2 }\cite[Chapter 10, section 3]{Minoux}. \emph{Each point }$%
x\in H$\emph{\ has the unique projection }$x_{pr}$\emph{\ on the set }$G.$%
\emph{\ Furthermore, }%
\begin{equation*}
\left( x_{pr}-x,z-x_{pr}\right) \geq 0,\forall z\in G.
\end{equation*}%
\emph{Define the projection operator }$P_{G}:H\rightarrow G$\emph{\ as }$%
P_{G}\left( x\right) =x_{pr}\in G.$\emph{\ Then} 
\begin{equation}
\left\Vert P_{G}\left( x\right) -P_{G}\left( y\right) \right\Vert \leq
\left\Vert x-y\right\Vert ,\forall x,y\in H.  \label{8.8}
\end{equation}

\textbf{Lemma 4.3 }\cite{KlibanovNik:ra2017}. \emph{The functional }$I\left(
x\right) $\emph{\ achieves its global minimal value at the point }$x_{\min
}\in G$\emph{\ on the set }$G$\emph{\ if and only if there exists a number }$%
\mu >0$ \emph{such that} 
\begin{equation}
x_{\min }=P_{G}\left( x_{\min }-\mu I^{\prime }\left( x_{\min }\right)
\right) .  \label{8.9}
\end{equation}%
\emph{If (\ref{8.9}) is valid for one number }$\mu ,$\emph{\ then it is also
valid for all numbers }$\mu >0.$

We now construct the gradient projection method of the minimization of the
functional $I\left( x\right) $ on the set $G.$ Choose an arbitrary point $%
x_{0}\in G$ \ and let 
\begin{equation}
x_{n+1}=P_{G}\left( x_{n}-\mu I^{\prime }\left( x_{n}\right) \right)
,n=0,1,2,...  \label{8.11}
\end{equation}

\textbf{Theorem 4.1}. \emph{Assume that conditions (\ref{8.1}) and (\ref{800}%
) are in place. Then there exists unique point of the relative minimum }$%
x_{\min }$\emph{\ of the functional }$I\left( x\right) $ \emph{on the set }$%
G.$\emph{\ In fact, }$x_{\min }$ \emph{is the unique point of the global
minimum of }$I\left( x\right) $\emph{\ on the set }$G.$\emph{\ Let }$D$\emph{%
\ and }$\varkappa $\emph{\ be the numbers in (\ref{8.1}) and (\ref{800})
respectively and let the number }$\varkappa \in \left( 0,D\right) .$\emph{\
Assume that the number }$\mu $\emph{\ in (\ref{8.11}) is so small that }%
\begin{equation}
0<\mu <\frac{2\varkappa }{D^{2}}.  \label{8.12}
\end{equation}%
\emph{Let} $q\left( \mu \right) =\left( 1-2\mu \varkappa +\mu
^{2}D^{2}\right) ^{1/2}.$ \emph{Then sequence (\ref{8.11}) converges to the
point of the global minimum }$x_{\min }$\emph{\ and }%
\begin{equation}
\left\Vert x_{n}-x_{\min }\right\Vert \leq q^{n}\left( \mu \right)
\left\Vert x_{0}-x_{\min }\right\Vert .  \label{8.13}
\end{equation}%
\emph{Furthermore, (\ref{8.3}) holds}.

\textbf{Proof}. First, we observe that since $\varkappa \in \left(
0,D\right) ,$ then (\ref{8.12}) implies that the number $q\left( \mu \right)
\in \left( 0,1\right) .$ Consider the operator%
\begin{equation*}
K:G\rightarrow G,
\end{equation*}%
\begin{equation*}
K\left( x\right) =P_{G}\left( x-\mu I^{\prime }\left( x\right) \right) .
\end{equation*}%
We now show that the operator $K$ is a contractual mapping operator. By (\ref%
{8.8}) we have for all $x,y\in G:$%
\begin{equation}
\begin{array}{c}
\left\Vert K\left( x\right) -K\left( y\right) \right\Vert ^{2}\leq
\left\Vert \left( x-\mu I^{\prime }\left( x\right) \right) -\left( y-\mu
I^{\prime }\left( y\right) \right) \right\Vert ^{2} \\ 
=\left\Vert \left( x-y\right) -\mu \left( I^{\prime }\left( x\right)
-I^{\prime }\left( y\right) \right) \right\Vert ^{2} \\ 
=\left\Vert x-y\right\Vert ^{2}+\mu ^{2}\left\Vert I^{\prime }\left(
x\right) -I^{\prime }\left( y\right) \right\Vert ^{2}-2\mu \left( I^{\prime
}\left( x\right) -I^{\prime }\left( y\right) ,x-y\right) .%
\end{array}
\label{8.14}
\end{equation}%
By (\ref{8.1}) $\mu ^{2}\left\Vert I^{\prime }\left( x\right) -I^{\prime
}\left( y\right) \right\Vert ^{2}\leq \mu ^{2}D^{2}\left\Vert x-y\right\Vert
^{2}.$ Next, by (\ref{803}) 
\begin{equation*}
-2\mu \left( I^{\prime }\left( x\right) -I^{\prime }\left( y\right)
,x-y\right) \leq -2\mu \varkappa \left\Vert x-y\right\Vert ^{2}.
\end{equation*}%
Hence, (\ref{8.14}) implies:%
\begin{equation*}
\left\Vert K\left( x\right) -K\left( y\right) \right\Vert ^{2}\leq \left(
1-2\mu \varkappa +\mu ^{2}D^{2}\right) \left\Vert x-y\right\Vert
^{2}=q^{2}\left( \mu \right) \left\Vert x-y\right\Vert ^{2}.
\end{equation*}%
Hence, the operator $K$ is \ a contraction mapping of the set $G.$ Hence,
there exists unique point%
\begin{equation}
x_{\min }=P_{G}\left( x_{\min }-\mu I^{\prime }\left( x_{\min }\right)
\right) ,x_{\min }\in G  \label{8.15}
\end{equation}%
and the convergence rate (\ref{8.13}) holds. Lemma 4.3 and (\ref{8.15})
imply that 
\begin{equation}
I\left( x_{\min }\right) =\min_{x\in G}I\left( x\right) .  \label{8.16}
\end{equation}%
Finally, Lemma 4.1 and (\ref{8.16}) imply that (\ref{8.3}) holds. \ \qed

\section{The Convexification Theorem}

\label{sec conv}

Using (\ref{1.3}) and (\ref{q sys}), define the set of admissible solutions $%
H$ as 
\begin{equation}
H=\left\{ 
\begin{array}{c}
q\in H^{5}(\Omega _{T}):q\left( \epsilon ,t\right) =g_{0}\left( t+\epsilon
\right) ,q_{x}\left( \epsilon ,t\right) =g_{1}\left( t+\epsilon \right)
+g_{0}^{\prime }\left( t+\epsilon \right) , \\ 
q_{x}\left( M,t\right) =0,\quad \mbox{for all }t\in \left[ 0,T\right] , \\ 
q\left( x,0\right) \geq \underline{q}=1/\left( 2\overline{c}^{1/4}\right)
>0,\quad \mbox{for all }x\in \left[ \epsilon ,M\right]%
\end{array}%
\right\} ,  \label{20}
\end{equation}%
see (\ref{1.1}) and (\ref{2.14}) for $\underline{q}$. We also define the
subspace $H_{0}$ of the space $H^{5}(\Omega _{T})$ as 
\begin{equation}
H_{0}=\big\{h\in H^{5}(\Omega _{T}):h(\epsilon ,t)=0,h_{x}(\epsilon
,t)=0,h_{x}(M,t)=0\mbox{ for all }t\in \lbrack 0,T]\big\}.  \label{21}
\end{equation}

Throughout this paper, we assume that the set $H$ is non empty. Let $R$ be
an arbitrary positive number. We define the set $B(R)$ as%
\begin{equation}
B(R)=\left\{ q\in H^{5}(\Omega _{T}):\left[ q\right] <R\right\} .
\label{1000}
\end{equation}%
The aim of this section is to prove that for all sufficiently large $\lambda 
$ and under some conditions imposed on $\beta $, the functional $J_{\lambda
,\alpha ,\beta }$ is strictly convex on the set $\overline{B(R)}\cap H.$
Theorem 5.1 below is our main result in this paper.

\subsection{The convexification theorem}

\label{sec 5.1}

\textbf{Theorem 5.1}. \emph{1. For any }$q\in H$\emph{\ and for any set of
parameters }$\lambda ,\alpha ,\beta $\emph{\ the functional }$J_{\lambda
,\alpha ,\beta }$\emph{\ has the Fr\'echet derivative }$J_{\lambda ,\alpha
,\beta }^{\prime }(q)\in H_{0}.$\emph{\ This derivative is Lipschitz
continuous in }$\overline{B\left( R\right) }\cap H$\emph{, i.e. there exists
a constant }$D>0$\emph{\ such that }%
\begin{equation}
\left[ J_{\lambda ,\alpha ,\beta }^{\prime }(q_{2})-J_{\lambda ,\alpha
,\beta }^{\prime }(q_{1})\right] \leq D\left[ q_{2}-q_{1}\right] ,\quad %
\mbox{for all }q_{1},q_{2}\in \overline{B(R)}\cap H.  \label{200}
\end{equation}

\emph{2. Let }$\overline{c}>0$\emph{\ be the number in (\ref{1.1}), let }$%
\alpha \in \left( 0,2\sqrt{\overline{c}}\right) $\emph{\ and let }$\lambda
_{0}\geq 1$\emph{\ be the same as in Theorem 3.1}$.$\emph{\ Then
there exists a constant }%
\begin{equation}
\lambda _{1}=\lambda _{1}\left( R,T,\epsilon ,\underline{c},\overline{c}%
,\alpha ,M\right) \geq \lambda _{0}  \label{22}
\end{equation}%
\emph{depending only on listed parameters such that for all }$\lambda \geq
\lambda _{1},$\emph{\ }$\beta \in \left[ 2e^{-\lambda \alpha T},1\right) $%
\emph{\ the functional }$J_{\lambda ,\alpha ,\beta }(q)$\emph{\ is strictly
convex on the set }$\overline{B\left( R\right) }\cap H.$\emph{\ More
precisely, the following inequality holds for an arbitrary pair of functions 
}$q,q+h\in \overline{B(R)}\cap H,$%
\begin{multline}
J_{\lambda ,\alpha ,\beta }(q+h)-J_{\lambda ,\alpha ,\beta }(q)-J_{\lambda
,\alpha ,\beta }^{\prime }(q)(h)\geq C\int_{\Omega _{T}}e^{-2\lambda
(x+\alpha t)}\Big[|h_{x}|^{2}+|h_{t}|^{2}+|h|^{2}\Big]dxdt \\
+C\int_{\epsilon }^{M}e^{-2\lambda x}\big(|(h_{x}(x,0)|^{2}+|h(x,0)|^{2}\big)%
dx+\frac{\beta }{2}\left[ h\right] ^{2},\quad \mbox{for
all }\lambda \geq \lambda _{1},  \label{3.11}
\end{multline}%
\emph{where the constant }$C=C\left( R,T,\epsilon ,\underline{c},\overline{c}%
,\alpha ,M\right) >0$\emph{\ depends only on listed parameters.}

\emph{3. There exists unique minimizer }$q_{\min }\in \overline{B(R)}\cap H$%
\emph{\ of the functional }$J_{\lambda ,\alpha ,\beta }(q)$\emph{\ on the
set }$\overline{B(R)}\cap H$\emph{\ and the following inequality holds:}%
\begin{equation}
\left[ J_{\lambda ,\alpha ,\beta }^{\prime }(q_{\min }),q_{\min }-q\right]
\leq 0,\text{ }\forall q\in \overline{B(R)}\cap H.  \label{900}
\end{equation}

\textbf{Proof.} Since both functions $q,q+h\in \overline{B\left( R\right) }%
\cap H$, then 
\begin{equation}
h\in \overline{B\left( 2R\right) }\cap H_{0}=\left\{ h\in H^{5}\left( \Omega
_{T}\right) :\left[ h\right] \leq 2R\right\} \cap H_{0}.  \label{202}
\end{equation}%
For every function $h$ satisfying (\ref{202}), denote $O(|h(x,0)|^{2})$ all
functions satisfying the inequality%
\begin{equation}
\left\vert O(|h(x,0)|^{2})\right\vert \leq C|h(x,0)|^{2}  \label{3.10}
\end{equation}%
and similarly for all other quantities in which the function $h$ and its
first order derivatives are involved. We also note that by (\ref{3.4}), %
(\ref{20}), (\ref{2}), (\ref{3}), (\ref{1000}) and (\ref{202})
\begin{equation}
\left\Vert h\right\Vert _{C^{2}\left( \overline{\Omega _{T}}\right)
},\left\Vert F\left( q\right) \right\Vert _{C^{1}\left( \overline{\Omega _{T}%
}\right) }\leq C,\quad \mbox{for all }h\in \overline{B\left( 2R\right) }\cap
H_{0},q\in \overline{B\left( R\right) }\cap H.  \label{4}
\end{equation}

We have%
\begin{equation}
F(q+h)=(q+h)_{xx}-\frac{(q+h)_{xt}}{2(q+h)^{2}(x,0)}+\frac{%
(q+h)_{t}(q+h)_{x}(x,0)}{2(q+h)^{3}(x,0)}.  \label{204}
\end{equation}%
It follows immediately from the Taylor formula that 
\begin{equation}
\frac{1}{2(q+h)^{2}(x,0)}=\frac{1}{2q^{2}(x,0)}-\frac{h\left( x,0\right) }{%
q^{3}(x,0)}+O\left( |h(x,0)|^{2}\right) ,  \label{205}
\end{equation}%
and 
\begin{equation}
\frac{1}{2(q+h)^{3}(x,0)}=\frac{1}{2q^{3}(x,0)}-\frac{3h\left( x,0\right) }{%
2q^{4}(x,0)}+O\left( |h(x,0)|^{2}\right) .  \label{206}
\end{equation}%
Using (\ref{3.4}) and (\ref{204})-(\ref{206}), we obtain%
\begin{multline}
F(q+h)=F\left( q\right) +\Big[\Big(h_{xx}-\frac{h_{xt}}{2q^{2}(x,0)}\Big)%
+h\left( x,0\right) \frac{q_{xt}}{q^{3}(x,0)}+h_{t}\frac{q_{x}\left(
x,0\right) }{2q^{3}(x,0)}+h_{x}\frac{q_{t}}{2q^{3}(x,0)} \\
-h\left( x,0\right) \frac{3q_{t}q_{x}\left( x,0\right) }{2q^{4}(x,0)}\Big]+%
\frac{h_{xt}h\left( x,0\right) }{2q^{3}(x,0)} \\
+O\left( \left\vert h\left( x,0\right) \right\vert ^{2}\right) +O\left(
\left\vert h_{t}\right\vert \left\vert h_{x}\left( x,0\right) \right\vert
\right) +O\left( \left\vert h_{t}\right\vert \left\vert h\left( x,0\right)
\right\vert \right) +O\left( \left\vert h_{x}\left( x,0\right) \right\vert
\left\vert h\left( x,0\right) \right\vert \right) .  \label{207}
\end{multline}%
%
%
%
%
%\textcolor{red}{\hl{AAAAA}}
Denote 
\begin{multline}
L_{\text{lin}}\left( h\right) =\Big[\Big(h_{xx}-\frac{h_{xt}}{2q^{2}(x,0)}%
\Big)+h\left( x,0\right) \frac{q_{xt}}{q^{3}(x,0)}+h_{t}\frac{q_{x}\left(
x,0\right) }{2q^{3}(x,0)}+h_{x}\frac{q_{t}}{2q^{3}(x,0)} \\
-h\left( x,0\right) \frac{3q_{t}q_{x}\left( x,0\right) }{2q^{4}(x,0)}\Big],
\label{208}
\end{multline}%
\begin{equation}
L_{\text{nonlin}}^{\left( 1\right) }\left( h\right) =\frac{h_{xt}h\left(
x,0\right) }{2q^{3}(x,0)},  \label{209}
\end{equation}%
and 
\begin{multline}
L_{\text{nonlin}}^{\left( 2\right) }\left( h\right) =O\left( \left\vert
h\left( x,0\right) \right\vert ^{2}\right) +O\left( \left\vert
h_{t}\right\vert \left\vert h_{x}\left( x,0\right) \right\vert \right)
+O\left( \left\vert h_{t}\right\vert \left\vert h\left( x,0\right)
\right\vert \right) \\
+O\left( \left\vert h_{x}\left( x,0\right) \right\vert \left\vert h\left(
x,0\right) \right\vert \right) .  \label{210}
\end{multline}%
Clearly, the operator $L_{\text{lin}}\left( h\right) $ depends linearly on $%
h $ and operators $L_{\text{nonlin}}^{\left( 1\right) }\left( h\right) ,L_{%
\text{nonlin}}^{\left( 2\right) }\left( h\right) $ depend nonlinearly. Using %
(\ref{207})-(\ref{210}) we obtain 
\begin{equation*}
F(q+h)=F\left( q\right) +L_{\text{lin}}\left( h\right) +L_{\text{nonlin}%
}^{\left( 1\right) }\left( h\right) +L_{\text{nonlin}}^{\left( 2\right)
}\left( h\right) .
\end{equation*}%
Hence,%
\begin{multline}
\left\vert F(q+h)\right\vert ^{2}-\left\vert F\left( q\right) \right\vert
^{2}=2F\left( q\right) L_{\text{lin}}\left( h\right) +2F\left( q\right) L_{%
\text{nonlin}}^{\left( 1\right) }\left( h\right) \\
+2F\left( q\right) L_{\text{nonlin}}^{\left( 2\right) }\left( h\right)
+\left\vert L_{\text{lin}}\left( h\right) +L_{\text{nonlin}}^{\left(
1\right) }\left( h\right) +L_{\text{nonlin}}^{\left( 2\right) }\left(
h\right) \right\vert ^{2}.  \label{3.5555}
\end{multline}%
Hence, (\ref{3.5555}) implies that 
\begin{multline}
J_{\lambda ,\alpha ,\beta }(q+h)-J_{\lambda ,\alpha ,\beta
}(q)-2\int_{\Omega _{T}}e^{-2\lambda (x+\alpha t)}F\left( q\right) L_{\text{%
lin}}\left( h\right) dxdt+2\beta \left[ q,h\right] \\
=2\int_{\Omega _{T}}e^{-2\lambda (x+\alpha t)}\left\vert L_{\text{lin}%
}\left( h\right) +L_{\text{nonlin}}^{\left( 1\right) }\left( h\right) +L_{%
\text{nonlin}}^{\left( 2\right) }\left( h\right) \right\vert ^{2}dxdt+\beta
\Vert h\Vert _{H^{5}(\Omega _{T})}^{2} \\
+2\int_{\Omega _{T}}e^{-2\lambda (x+\alpha t)}F\left( q\right) L_{\text{%
nonlin}}^{\left( 1\right) }\left( h\right) dxdt+2\int_{\Omega
_{T}}e^{-2\lambda (x+\alpha t)}F\left( q\right) L_{\text{nonlin}}^{\left(
2\right) }\left( h\right) dxdt.  \label{211}
\end{multline}%
In particular, it follows from (\ref{209})-(\ref{211}) that%
\begin{multline}
\left\vert J_{\lambda ,\alpha ,\beta }(q+h)-J_{\lambda ,\alpha ,\beta
}(q)-2\int_{\Omega _{T}}e^{-2\lambda (x+\alpha t)}F\left( q\right) L_{\text{%
lin}}\left( h\right) dxdt+2\beta \left[ q,h\right] \right\vert \\
=o\left( \left[ h\right] \right) \text{ as }\left[ h\right] \rightarrow 0.
\label{212}
\end{multline}%
Denote%
\begin{equation}
\widetilde{J}_{\lambda ,\alpha ,\beta ,q}\left( h\right) =-2\int_{\Omega
_{T}}e^{-2\lambda (x+\alpha t)}F\left( q\right) L_{\text{lin}}\left(
h\right) dxdt+2\beta \left[ q,h\right] .  \label{213}
\end{equation}%
Then $Z_{\lambda ,\alpha ,\beta ,q}:H_{0}\rightarrow \mathbb{R}$ is a
bounded linear functional for every $q\in H.$ Hence, by Riesz theorem, there
exists unique function $J_{\lambda ,\alpha ,\beta }^{\prime }(q)\in H_{0}$
such that 
\begin{equation}
Z_{\lambda ,\alpha ,\beta ,q}\left( h\right) =\left[ J_{\lambda ,\alpha
,\beta }^{\prime }(q),h\right] ,\quad \mbox{ for all }h\in H_{0}.
\label{214}
\end{equation}%
Therefore, it follows from (\ref{212})-(\ref{214}) that 
\begin{equation}
J_{\lambda ,\alpha ,\beta }^{\prime }(q)\in H_{0}  \label{215}
\end{equation}%
is the Fréchet derivative of the functional $J_{\lambda ,\alpha ,\beta }$ at
the point $q\in H.$ The Lipschitz continuity property (\ref{200}) of $%
J_{\lambda ,\alpha ,\beta }^{\prime }(q)$ can be proven completely similarly
with the proof of Theorem 3.1 in \cite{KlibanovNik:ra2017}. Hence, we omit
the proof of (\ref{200}).

Thus, (\ref{211})-(\ref{215}) imply that 
\begin{multline}
J_{\lambda ,\alpha ,\beta }(q+h)-J_{\lambda ,\alpha ,\beta }(q)-J_{\lambda
,\alpha ,\beta }^{\prime }(q)\left( h\right) \\
=2\int_{\Omega _{T}}e^{-2\lambda (x+\alpha t)}\left\vert L_{\text{lin}%
}\left( h\right) +L_{\text{nonlin}}^{\left( 1\right) }\left( h\right) +L_{%
\text{nonlin}}^{\left( 2\right) }\left( h\right) \right\vert ^{2}dxdt+\beta %
\left[ h\right] ^{2} \\
+2\int_{\Omega _{T}}e^{-2\lambda (x+\alpha t)}F\left( q\right) L_{\text{%
nonlin}}^{\left( 1\right) }\left( h\right) dxdt+2\int_{\Omega
_{T}}e^{-2\lambda (x+\alpha t)}F\left( q\right) L_{\text{nonlin}}^{\left(
2\right) }\left( h\right) dxdt.  \label{216}
\end{multline}%
We now estimate the right hand side of (\ref{216}) from the below. First, we
rewrite formula (\ref{208}) as%
\begin{equation}
L_{\text{lin}}\left( h\right) =L_{\text{lin}}^{\left( 1\right) }\left(
h\right) +L_{\text{lin}}^{\left( 2\right) }\left( h\right) ,  \label{217}
\end{equation}%
where 
\begin{align}
L_{\text{lin}}^{\left( 1\right) }\left( h\right) & =h_{xx}-\frac{h_{xt}}{%
2q^{2}(x,0)},  \label{218} \\
L_{\text{lin}}^{\left( 2\right) }\left( h\right) & =h\left( x,0\right) \frac{%
q_{xt}}{q^{3}(x,0)}+h_{t}\frac{q_{x}\left( x,0\right) }{2q^{3}(x,0)}+h_{x}%
\frac{q_{t}}{2q^{3}(x,0)}-h\left( x,0\right) \frac{3q_{t}q_{x}\left(
x,0\right) }{2q^{4}(x,0)}.  \label{219}
\end{align}%
Hence, Cauchy-Schwarz inequality, (\ref{3.10}) and (\ref{4}) imply that 
\begin{equation*}
\left\vert L_{\text{lin}}\left( h\right) \right\vert ^{2}\geq \left\vert L_{%
\text{lin}}^{\left( 1\right) }\left( h\right) \right\vert ^{2}-2L_{\text{lin}%
}^{\left( 1\right) }\left( h\right) L_{\text{lin}}^{\left( 2\right) }\left(
h\right) \geq \frac{1}{2}\left\vert L_{\text{lin}}^{\left( 1\right) }\left(
h\right) \right\vert ^{2}-2\left\vert L_{\text{lin}}^{\left( 2\right)
}\left( h\right) \right\vert ^{2}
\end{equation*}%
\begin{equation*}
=\frac{1}{2}\left\vert L_{\text{lin}}^{\left( 1\right) }\left( h\right)
\right\vert ^{2}+O\left( \left\vert h\left( x,0\right) \right\vert
^{2}\right) +O\left( h_{t}^{2}\right) +O\left( h_{x}^{2}\right) .
\end{equation*}%
Thus, 
\begin{equation}
\left\vert L_{\text{lin}}\left( h\right) \right\vert ^{2}\geq \frac{1}{2}%
\left\vert L_{\text{lin}}^{\left( 1\right) }\left( h\right) \right\vert
^{2}+O\left( \left\vert h\left( x,0\right) \right\vert ^{2}\right) +O\left(
h_{t}^{2}\right) +O\left( h_{x}^{2}\right) .  \label{220}
\end{equation}%
Next, (\ref{3.10}), (\ref{4}), (\ref{209}), (\ref{210}) and (\ref{220})
imply that 
\begin{align*}
\big|L_{\text{lin}}\left( h\right) & +L_{\text{nonlin}}^{\left( 1\right)
}\left( h\right) +L_{\text{nonlin}}^{\left( 2\right) }\left( h\right) \big|%
^{2}\geq \left\vert L_{\text{lin}}\left( h\right) \right\vert ^{2}+2L_{\text{%
lin}}\left( h\right) \left( L_{\text{nonlin}}^{\left( 1\right) }\left(
h\right) +L_{\text{nonlin}}^{\left( 2\right) }\left( h\right) \right) \\
& \geq \frac{1}{2}\left\vert L_{\text{lin}}\left( h\right) \right\vert
^{2}-2\left( L_{\text{nonlin}}^{\left( 1\right) }\left( h\right) +L_{\text{%
nonlin}}^{\left( 2\right) }\left( h\right) \right) ^{2} \\
& \geq \frac{1}{4}\left\vert L_{\text{lin}}^{\left( 1\right) }\left(
h\right) \right\vert ^{2}+O\left( \left\vert h\left( x,0\right) \right\vert
^{2}\right) +O\left( \left\vert h_{x}\left( x,0\right) \right\vert
^{2}\right) +O\left( h_{t}^{2}\right) +O\left( h_{x}^{2}\right) .
\end{align*}%
Thus, 
\begin{multline}
\left\vert L_{\text{lin}}\left( h\right) +L_{\text{nonlin}}^{\left( 1\right)
}\left( h\right) +L_{\text{nonlin}}^{\left( 2\right) }\left( h\right)
\right\vert ^{2} \\
\geq \frac{1}{4}\left\vert L_{\text{lin}}^{\left( 1\right) }\left( h\right)
\right\vert ^{2}+O\left( \left\vert h\left( x,0\right) \right\vert
^{2}\right) +O\left( \left\vert h_{x}\left( x,0\right) \right\vert
^{2}\right) +O\left( h_{t}^{2}\right) +O\left( h_{x}^{2}\right) .
\label{221}
\end{multline}%
Next, by (\ref{3.10}), (\ref{4}) and (\ref{210})
\begin{equation}
F\left( q\right) L_{\text{nonlin}}^{\left( 2\right) }\left( h\right) \geq
O\left( \left\vert h\left( x,0\right) \right\vert ^{2}\right) +O\left(
\left\vert h_{x}\left( x,0\right) \right\vert ^{2}\right) +O\left(
h_{t}^{2}\right) .  \label{222}
\end{equation}%
Summing up (\ref{221}) and (\ref{222}) and substituting then in (\ref{216}),
we obtain 
\begin{multline}
J_{\lambda ,\alpha ,\beta }(q+h)-J_{\lambda ,\alpha ,\beta }(q)-J_{\lambda
,\alpha ,\beta }^{\prime }(q)\left( h\right) \geq \frac{1}{2}\int_{\Omega
_{T}}e^{-2\lambda (x+\alpha t)}\left[ \left\vert L_{\text{lin}}^{\left(
1\right) }\left( h\right) \right\vert ^{2}+4F\left( q\right) L_{\text{nonlin}%
}^{\left( 1\right) }\left( h\right) \right] dxdt \\
+\beta \Vert h\Vert _{H^{5}(\Omega _{T})}^{2}-C\int_{\Omega
_{T}}e^{-2\lambda (x+\alpha t)}\left[ \left\vert h\left( x,0\right)
\right\vert ^{2}+\left\vert h_{x}\left( x,0\right) \right\vert
^{2}+h_{t}^{2}+h_{x}^{2}\right] dxdt.  \label{223}
\end{multline}

We now apply the second Carleman estimate (\ref{101}) of Theorem 3.1, recalling that $h\left( \epsilon ,t\right) =h_{t}\left( \epsilon
,t\right) =h_{x}\left( \epsilon ,t\right) =0$ and using (\ref{218}), 
\begin{multline}
\frac{1}{2}\int_{\Omega _{T}}e^{-2\lambda (x+\alpha t)}\left[ \left\vert L_{%
\text{lin}}^{\left( 1\right) }\left( h\right) \right\vert ^{2}+4F\left(
q\right) L_{\text{nonlin}}^{\left( 1\right) }\left( h\right) \right] dxdt \\
\geq C\lambda \int_{\Omega _{T}}e^{-2\lambda (x+\alpha t)}\big(%
|h_{x}(x,t)|^{2}+|h_{t}(x,t)|^{2}+\lambda ^{2}|h(x,t)|^{2}\big)dtdx \\
+C\lambda \int_{\epsilon }^{M}e^{-2\lambda x}\big(|h_{x}(x,0)|^{2}+\lambda
^{2}|h(x,0)|^{2}\big)dx \\
-C\lambda \int_{\epsilon }^{M}e^{-2\lambda (x+\alpha T)}\big(%
|h_{x}(x,T)|^{2}+\lambda ^{2}|h(x,T)|^{2}\big)dx,\quad \mbox{for all }%
\lambda \geq \lambda _{0}\geq 1,  \label{224}
\end{multline}%
where $\lambda _{0}$ is was chosen in Theorem 3.1. Since constants 
$C$ are different in Theorem 3.1 and (\ref{224}), we denote $C$ in %
(\ref{224}) as $\widetilde{C}.$ Substituting (\ref{224}) in (\ref{223}), we
obtain%
\begin{multline}
J_{\lambda ,\alpha ,\beta }(q+h)-J_{\lambda ,\alpha ,\beta }(q)-J_{\lambda
,\alpha ,\beta }^{\prime }(q)\left( h\right) \\
\geq C\lambda \int_{\Omega _{T}}e^{-2\lambda (x+\alpha t)}\big(%
|h_{x}(x,t)|^{2}+|h_{t}(x,t)|^{2}+\lambda ^{2}|h(x,t)|^{2}\big)dtdx \\
+C\lambda \int_{\epsilon }^{M}e^{-2\lambda x}\big(|h_{x}(x,0)|^{2}+\lambda
^{2}|h(x,0)|^{2}\big)dx \\
-\widetilde{C}\int_{\Omega _{T}}e^{-2\lambda (x+\alpha t)}\left[ \left\vert
h\left( x,0\right) \right\vert ^{2}+\left\vert h_{x}\left( x,0\right)
\right\vert ^{2}+h_{t}^{2}+h_{x}^{2}\right] dxdt \\
-C\lambda \int_{\epsilon }^{M}e^{-2\lambda (x+\alpha T)}\big(%
|h_{x}(x,T)|^{2}+\lambda ^{2}|h(x,T)|^{2}\big)dx+\beta \Vert h\Vert
_{H^{5}(\Omega _{T})}^{2}  \label{225}
\end{multline}%
for all $\lambda \geq \lambda _{0}\geq 1.$ Choose $\lambda _{1}\geq \lambda
_{0}$ so large that $C\lambda _{1}/2>\widetilde{C}.$ Then, replacing for
convenience in (\ref{225}) $C/2$ with $C$ again, we obtain%
\begin{multline}
J_{\lambda ,\alpha ,\beta }(q+h)-J_{\lambda ,\alpha ,\beta }(q)-J_{\lambda
,\alpha ,\beta }^{\prime }(q)\left( h\right) \\
\geq C\lambda \int_{\Omega _{T}}e^{-2\lambda (x+\alpha t)}\big(%
|h_{x}(x,t)|^{2}+|h_{t}(x,t)|^{2}+\lambda ^{2}|h(x,t)|^{2}\big)dtdx \\
+C\lambda \int_{\epsilon }^{M}e^{-2\lambda x}\big(|h_{x}(x,0)|^{2}+\lambda
^{2}|h(x,0)|^{2}\big)dx \\
-C\lambda \int_{\epsilon }^{M}e^{-2\lambda (x+\alpha T)}\big(%
|h_{x}(x,T)|^{2}+\lambda ^{2}|h(x,T)|^{2}\big)dx+\beta \left[ h\right] ^{2},
\label{226}
\end{multline}%
for all $\lambda \geq \lambda _{0}\geq 1.$ Now, since $\lambda _{1}$ is
sufficiently large, then $2e^{-\lambda \alpha T}>C\lambda ^{3}e^{-2\lambda
(\epsilon +\alpha T)},$ for all $\lambda \geq \lambda _{1}.$ Hence, setting
in (\ref{226}) $\beta \in \left[ 2e^{-\lambda \alpha T},1\right) $ and
recalling that by the trace theorem $\left\Vert h(x,T)\right\Vert
_{H^{1}\left( \epsilon ,M\right) }\leq C\left[ h\right] ,$ we obtain desired
estimate (\ref{3.11}).

The existence and uniqueness of the minimizer $q_{\min }\in \overline{%
B\left( R\right) }\cap H$ of the functional $J_{\lambda ,\alpha ,\beta }(q)$
on the set $\overline{B\left( R\right) }\cap H$ as well as inequality %
(\ref{900}) follow from Theorem 4.1, (\ref{200}) and the strict convexity
property (\ref{3.11}). \qed

\subsection{The accuracy of the minimizer}

\label{sec 5.2}

Following the concept of Tikhonov for ill-posed problems \cite%
{Tihkonov:kapg1995}, we assume that there exists an \textquotedblleft ideal"
solution $q^{\ast }$ of problem (\ref{q sys}) with \textquotedblleft ideal",
i.e. noiseless data $g_{0}^{\ast }$ and $g_{1}^{\ast }.$ More precisely, let 
\begin{equation*}
H^{\ast }=\left\{ 
\begin{array}{c}
q\in H^{5}(\Omega _{T}):q\left( \epsilon ,t\right) =g_{0}^{\ast }\left(
t+\epsilon \right) ,q_{x}\left( \epsilon ,t\right) =g_{1}^{\ast }\left(
t+\epsilon \right) +\left( g_{0}^{\ast }\right) ^{\prime }\left( t+\epsilon
\right) , \\ 
q_{x}\left( M,t\right) =0,\quad \mbox{for all }t\in \left[ 0,T\right] , \\ 
q\left( x,0\right) \geq \underline{q}=1/\left( 2\overline{c}^{1/4}\right)
,\quad \mbox{for all }x\in \left[ \epsilon ,M\right]%
\end{array}%
\right\} .
\end{equation*}%
Keeping in mind the result of subsection 5.3, we assume that 
\begin{equation}
q^{\ast }\in B\left( R/3\right) \cap H^{\ast }.  \label{603}
\end{equation}%
Also, let a sufficiently small number $\delta >0$ be the noise level in the
data. Introduce the set $H^{\delta }$ as%
\begin{equation}
H^{\delta }=\left\{ 
\begin{array}{c}
q\in H^{5}(\Omega _{T}):q\left( \epsilon ,t\right) =g_{0}^{\delta }\left(
t+\epsilon \right) ,q_{x}\left( \epsilon ,t\right) =g_{1}^{\delta }\left(
t+\epsilon \right) +\left( g_{0}^{\delta }\right) ^{\prime }\left(
t+\epsilon \right) , \\ 
q_{x}\left( M,t\right) =0,\quad \mbox{for all }t\in \left[ 0,T\right] , \\ 
q\left( x,0\right) \geq \underline{q}=1/\left( 2\overline{c}^{1/4}\right)
,\quad \mbox{for all }x\in \left[ \epsilon ,M\right]%
\end{array}%
\right\} ,  \label{600}
\end{equation}%
where $g_{0}^{\delta }$ and $g_{1}^{\delta }$ are noisy data $g_{0}$ and $%
g_{1}.$ Suppose that there exists a function $G^{\delta }\left( x,t\right)
\in H^{5}(\Omega _{T})\cap H^{\delta }.$ There exists a function $G^{\ast
}\left( x,t\right) \in H^{5}(\Omega _{T})\cap H^{\ast }.$ We assume that 
\begin{equation}
\left[ G^{\delta }-G^{\ast }\right] <\delta .  \label{604}
\end{equation}

Let $c^{\ast }\left( x\right) $ be the exact solution of Problem \ref{cip}.
Similarly with (\ref{q2c}), we define 
\begin{equation}
c^{\ast }\left( x\right) =\frac{1}{16\left( q^{\ast }\left( x,0\right)
\right) ^{4}}\leq \overline{c},\text{ }x\in \left[ \epsilon ,M\right] .
\label{605}
\end{equation}%
Let $q_{\min }\in \overline{B\left( R\right) }\cap H^{\delta }$ be the
minimizer of the functional $J_{\lambda ,\alpha ,\beta }\left( q\right) ,$
which is found in Theorem 5.1. To indicate the dependence of $q_{\min }$\ on
the noise level $\delta ,$\ we denote $q_{\min }$ as $q_{\min }^{\delta }.$\
Following (\ref{q2c}), define the function $c_{\min }^{\delta }\left(
x\right) ,$\ which corresponds to the function $q_{\min }^{\delta },$\ as 
\begin{equation}
c_{\min }^{\delta }\left( x\right) =\frac{1}{16\left( q_{\min }^{\delta
}\left( x,0\right) \right) ^{4}},\text{ }x\in \left[ \epsilon ,M\right] .
\label{6050}
\end{equation}

\textbf{Theorem 5.2 (}stability of minimizers in the presence of noise%
\textbf{). }\emph{Assume that }$\overline{B\left( R/3\right) }\cap H^{\delta
}\neq \varnothing $\emph{. Let the function }$q^{\ast }$\emph{be the exact
solution of problem (\ref{q sys}). Suppose that conditions (\ref{603}) and (%
\ref{604}) are in place. Let the set }$H$\emph{\ in (\ref{20}) be replaced
with the set }$H^{\delta }$\emph{\ in (\ref{600}). Let }$c^{\ast }\left(
x\right) $\emph{\ be the exact target coefficient, as in (\ref{605}). Let
parameters }$\lambda _{1}$\emph{\ and }$\alpha $\emph{\ be the same as in
Theorem 5.1. Also, similarly with Theorem 5.1, let }$\lambda \geq \lambda
_{1}$\emph{\ and }$\beta =2e^{-\lambda \alpha T}$\emph{. Assume that} 
\begin{equation}
\alpha T>2M.  \label{606}
\end{equation}%
\emph{Choose a number }$T_{0}\in \left( 0,T\right) $\emph{\ such that }%
\begin{equation}
\alpha T>2M+2\alpha T_{0},  \label{607}
\end{equation}%
\emph{which is possible by (\ref{606}). Choose the number }$\delta
_{0}=\delta _{0}\left( R,T,\epsilon ,\underline{c},\overline{c},\alpha
,M\right) >0$\emph{\ depending only on listed parameters so small that }%
\begin{equation}
\delta _{0}<2R/3,  \label{699}
\end{equation}%
\emph{\ }%
\begin{equation}
\ln \left( \delta _{0}^{-2/\left( \alpha T\right) }\right) \geq \lambda
_{1}\left( R,T,\epsilon ,\underline{c},\overline{c},\alpha ,M\right) .
\label{700}
\end{equation}%
\emph{Let }$\delta \in \left( 0,\delta _{0}\right) .$\emph{\ Introduce two
numbers }$\rho _{1},\rho _{2}\in \left( 0,1\right) :$%
\begin{equation}
\text{ }\rho _{1}=\frac{\alpha T-2\left( M+\alpha T_{0}\right) }{\alpha T},%
\text{ }\rho _{2}=\frac{\alpha T-2M}{\alpha T};\text{ }\rho _{1},\rho
_{2}\in \left( 0,1\right) .  \label{611}
\end{equation}%
\emph{Then the following stability estimates are valid for all }$\delta \in
\left( 0,\delta _{0}\right) :$%
\begin{equation}
\left\Vert q^{\ast }-q_{\min }^{\delta }\right\Vert _{H^{1}\left( \Omega
_{T_{0}}\right) }\leq C\delta ^{\rho _{1}},\text{ }  \label{608}
\end{equation}%
\begin{equation}
\left\Vert q^{\ast }\left( x,0\right) -q_{\min }^{\delta }\left( x,0\right)
\right\Vert _{H^{1}\left( \epsilon ,M\right) }\leq C\delta ^{\rho _{2}},
\label{609}
\end{equation}%
\begin{equation}
\left\Vert c^{\ast }\left( x\right) -c_{\min }^{\delta }\left( x,0\right)
\right\Vert _{H^{1}\left( \epsilon ,M\right) }\leq C\delta ^{\rho _{2}},
\label{610}
\end{equation}%
\emph{where the constant }$C=C\left( R,T,\epsilon ,\underline{c},\overline{c}%
,\alpha ,M\right) >0$\emph{\ depends only on listed parameters.}

\textbf{Proof.} Denote $p^{\ast }=q^{\ast }-G^{\ast }.$ Consider the
function $p^{\ast }+G.$ Then $p^{\ast }+G\in H^{\delta }.$ In addition, $%
p^{\ast }+G\in B\left( R\right) .$ Indeed, using (\ref{603}), (\ref{604}), (%
\ref{699}) and triangle inequality, we obtain 
\begin{equation*}
\left[ p^{\ast }+G\right] =\left[ q^{\ast }+\left( G-G^{\ast }\right) \right]
\leq \left[ q^{\ast }\right] +\left[ G-G^{\ast }\right] <\frac{R}{3}+\delta
<R.
\end{equation*}

Thus, $p^{\ast }+G\in B\left( R\right) \cap H^{\delta }.$ Hence, we can use
Theorem 5.1, for all $\lambda \geq \lambda _{1}.$ Thus, by (\ref{3.11}) 
\begin{equation*}
J_{\lambda ,\alpha ,\beta }(p^{\ast }+G)-J_{\lambda ,\alpha ,\beta }(q_{\min
}^{\delta })-J_{\lambda ,\alpha ,\beta }^{\prime }(q_{\min }^{\delta
})\left( p^{\ast }+G-q_{\min }^{\delta }\right)
\end{equation*}%
\begin{equation*}
\geq C\int_{\Omega _{T}}e^{-2\lambda (x+\alpha t)}\Big[|\left( p^{\ast
}+G-q_{\min }^{\delta }\right) _{x}|^{2}+|\left( p^{\ast }+G-q_{\min
}^{\delta }\right) _{t}|^{2}\Big]dxdt
\end{equation*}%
\begin{equation}
+C\int_{\Omega _{T}}e^{-2\lambda (x+\alpha t)}|\left( p^{\ast }+G-q_{\min
}^{\delta }\right) |^{2}dxdt  \label{616}
\end{equation}%
\begin{equation*}
+C\int_{\epsilon }^{M}e^{-2\lambda x}\big(|(\left( p^{\ast }+G-q_{\min
}^{\delta }\right) _{x}(x,0)|^{2}+|\left( p^{\ast }+G-p_{\min }^{\delta
}\right) (x,0)|^{2}\big)dx.
\end{equation*}%
Next, since by (\ref{900}) $-J_{\lambda ,\alpha ,\beta }^{\prime }(q_{\min
}^{\delta })\left( p^{\ast }+G-q_{\min }^{\delta }\right) \leq 0$ and also $%
-J_{\lambda ,\alpha ,\beta }(q_{\min }^{\delta })\leq 0,$ then (\ref{616})
implies:%
\begin{equation*}
J_{\lambda ,\alpha ,\beta }(p^{\ast }+G)\geq
\end{equation*}%
\begin{equation*}
\geq C\int_{\Omega _{T}}e^{-2\lambda (x+\alpha t)}\Big[|\left( p^{\ast
}+G-q_{\min }^{\delta }\right) _{x}|^{2}+|\left( p^{\ast }+G-q_{\min
}^{\delta }\right) _{t}|^{2}\Big]dxdt
\end{equation*}%
\begin{equation*}
+C\int_{\Omega _{T}}e^{-2\lambda (x+\alpha t)}|\left( p^{\ast }+G-q_{\min
}^{\delta }\right) |^{2}dxdt
\end{equation*}%
\begin{equation}
+C\int_{\epsilon }^{M}e^{-2\lambda x}\big(|(\left( p^{\ast }+G-q_{\min
}^{\delta }\right) _{x}(x,0)|^{2}+|\left( p^{\ast }+G-p_{\min }^{\delta
}\right) (x,0)|^{2}\big)dx.  \label{617}
\end{equation}%
By (\ref{2.25}) 
\begin{equation*}
J_{\lambda ,\alpha ,\beta }(p^{\ast }+G^{\delta })=\int_{\Omega
_{T}}e^{-2\lambda (x+\alpha t)}\big|F(p^{\ast }+G^{\delta })\big|^{2}dxdt
\end{equation*}%
\begin{equation}
+\beta \left[ p^{\ast }+G^{\delta }\right] ^{2}  \label{618}
\end{equation}%
\begin{equation*}
=J_{\lambda ,\alpha ,\beta }^{0}(q^{\ast }+\left( G^{\delta }-G^{\ast
}\right) )+2e^{-\lambda T}\left[ q^{\ast }+\left( G^{\delta }-G^{\ast
}\right) \right] ^{2}.
\end{equation*}%
We have $F(p^{\ast }+G^{\ast })=F\left( q^{\ast }\right) =0.$ Hence, (\ref%
{604}) implies 
\begin{equation}
\big|F(q^{\ast }+\left( G^{\delta }-G^{\ast }\right) )\big|^{2}\leq C\delta
^{2}.  \label{6180}
\end{equation}%
Hence, using (\ref{618}) and (\ref{6180}), we obtain 
\begin{equation*}
J_{\lambda ,\alpha ,\beta }^{0}(q^{\ast }+\left( G^{\delta }-G^{\ast
}\right) )\leq C\delta ^{2}.
\end{equation*}%
Hence, by (\ref{618}) 
\begin{equation}
J_{\lambda ,\alpha ,\beta }(p^{\ast }+G^{\delta })\leq C\left( \delta
^{2}+\beta \right) =C\left( \delta ^{2}+2e^{-\lambda \alpha T}\right) .
\label{619}
\end{equation}

Now, since $T_{0}\in \left( 0,T\right) ,$ then $\Omega _{T_{0}}\subset
\Omega _{T}.$ Also,%
\begin{equation*}
e^{-2\lambda (x+\alpha t)}\geq e^{-2\lambda (M+\alpha T_{0})}\text{ in }%
\Omega _{T_{0}}.
\end{equation*}%
Replacing in the first integral of (\ref{617}) $\Omega _{T}$ with $\Omega
_{T_{0}},$ we only make inequality (\ref{617}) stronger. Hence, (\ref{617})
and (\ref{619}) lead to%
\begin{equation*}
\left\Vert q^{\ast }-q_{\min }^{\delta }+\left( G^{\delta }-G^{\ast }\right)
\right\Vert _{H^{1}\left( \Omega _{T_{0}}\right) }\leq Ce^{\lambda (M+\alpha
T_{0})}\left( \delta +\sqrt{\beta /2}\right)
\end{equation*}%
\begin{equation*}
=Ce^{\lambda (M+\alpha T_{0})}\left( \delta +e^{-\lambda \alpha T/2}\right) ,
\end{equation*}%
\begin{equation*}
\left\Vert q^{\ast }-q_{\min }^{\delta }+\left( G^{\delta }-G^{\ast }\right)
\right\Vert _{H^{1}\left( \epsilon ,M\right) }\leq Ce^{\lambda M}\left(
\delta +e^{-\lambda \alpha T/2}\right) .
\end{equation*}%
The triangle inequality, (\ref{604}) and the last two estimates lead to:%
\begin{equation}
\left\Vert q^{\ast }-q_{\min }^{\delta }\right\Vert _{H^{1}\left( \Omega
_{T_{0}}\right) }\leq Ce^{\lambda (M+\alpha T_{0})}\left( \delta +\sqrt{%
\beta /2}\right)  \label{622}
\end{equation}%
\begin{equation*}
=Ce^{\lambda (M+\alpha T_{0})}\left( \delta +e^{-\lambda \alpha T/2}\right) ,
\end{equation*}%
\begin{equation}
\left\Vert q^{\ast }-q_{\min }^{\delta }\right\Vert _{H^{1}\left( \epsilon
,M\right) }\leq Ce^{\lambda M}\left( \delta +e^{-\lambda \alpha T/2}\right) .
\label{623}
\end{equation}%
Choose the number $\delta _{0}=\delta _{0}\left( R,T,\epsilon ,\underline{c},%
\overline{c},\alpha ,M\right) >0$ as in (\ref{700}). Let $\delta \in \left(
0,\delta _{0}\right) .$ In (\ref{622}) and (\ref{623}) choose $\lambda
=\lambda \left( \delta \right) >\lambda _{1}$ such that $e^{-\lambda \alpha
T/2}=\delta ,$ i.e. $\lambda =\ln \left( \delta ^{-2/\left( \alpha T\right)
}\right) .$ Then in (\ref{622})%
\begin{equation}
Ce^{\lambda (M+\alpha T_{0})}\left( \delta +e^{-\lambda \alpha T/2}\right)
=2C\delta ^{\rho _{1}},\text{ \ }\rho _{1}=\frac{\alpha T-2\left( M+\alpha
T_{0}\right) }{\alpha T}\in \left( 0,1\right) ,  \label{624}
\end{equation}%
and in (\ref{623})%
\begin{equation}
Ce^{\lambda M}\left( \delta +e^{-\lambda \alpha T/2}\right) =2C\delta ^{\rho
_{2}},\text{ }\rho _{2}=\frac{\alpha T-2M}{\alpha T}\in \left( 0,1\right) .
\label{625}
\end{equation}%
The fact that numbers $\rho _{1},\rho _{2}\in \left( 0,1\right) $ follows
from (\ref{607}). Estimates (\ref{611})-(\ref{609}) follow immediately from (%
\ref{622})-(\ref{625}). Estimate (\ref{610}) follows from (\ref{605}), (\ref%
{6050}), (\ref{609}) and the requirement $q\left( x,0\right) \geq \underline{%
q}=1/\left( 2\overline{c}^{1/4}\right) $ in both sets $H^{\delta }$ and $%
H^{\ast }.$ \qed

\subsection{The global convergence of the gradient descent method}

\label{sec 5.3}

Starting from the work \cite{KlibanovNik:ra2017}, in all above cited works
on the convexification, the global convergence of the gradient projection
method was proven, see (\ref{8.11}) for this method. However, it is hard to
practically implement projection operators. For this reason, a simpler
gradient descent method was used in those works and results were successful.
In two recent publications \cite{Klibanov:2ndSAR2021,LeNguyen:preprint2021}
the global convergence of the gradient descent method, being applied to some
analogs of the functional $J_{\lambda ,\alpha ,\beta }$, was proven, which
has justified those numerical results.

In this section, we first formulate an analog of that theorem of \cite%
{Klibanov:2ndSAR2021}, which is applicable to our case. Suppose that
assumptions of Theorem 5.2 are in place. Let 
\begin{equation}
q_{0}^{\delta }\in B\left( R/3\right) \cap H^{\delta }  \label{5.0}
\end{equation}%
be the starting point of the minimizing sequence of the gradient descent
method, 
\begin{equation}
q_{n}^{\delta }=q_{n-1}^{\delta }-\eta J_{\lambda ,\alpha ,\beta }^{\prime
}(q_{n-1}^{\delta }),\quad n=1,2,...,  \label{5.1}
\end{equation}%
where $\eta >0$ is a small step size, which we will choose later. Along with
the functions $q_{\min }^{\delta }$ and $q_{n}^{\delta },$ we also introduce
corresponding coefficients and $c_{n}^{\delta }\left( x\right) ,$ which are
calculated by formula (\ref{q2c}),%
\begin{equation}
\text{ }c_{n}^{\delta }\left( x\right) =\frac{1}{16\left( q_{n}^{\delta
}\left( x,0\right) \right) ^{4}},\quad x\in \left[ \epsilon ,M\right] .
\label{227}
\end{equation}%
Below functions $c_{\min }^{\delta }\left( x\right) $ are as in (\ref{6050}%
). Since we make sure below that our functions $q_{n}^{\delta },q_{\min
}^{\delta }\in B\left( R\right) \cap H^{\delta },$ then by the third line of
(\ref{600}) 
\begin{equation}
\frac{1}{16\left( q_{\min }^{\delta }\left( x,0\right) \right) ^{4}},\frac{1%
}{16\left( q_{n}^{\delta }\left( x,0\right) \right) ^{4}}\leq \overline{c}.
\label{228}
\end{equation}

\begin{remark}
Since by Theorem 5.1 $J_{\lambda ,\alpha ,\beta }^{\prime }(q_{n-1}^{\delta
})\in H_{0},$ then in (\ref{5.1}) boundary conditions of (\ref{600}) are the
same for all functions $q^{(n)},$ $n=1,2,\dots .$
\end{remark}

The following theorem follows immediately from either Theorem 4.6 of \cite%
{Klibanov:2ndSAR2021} or Theorem 2.2 of \cite{LeNguyen:preprint2021} as well
as from the trace theorem, Theorem 5.2, (\ref{227}) and (\ref{228}).

\textbf{Theorem 5.3.} \emph{Assume that conditions of Theorem 5.2 as well as
conditions (\ref{5.0}) and (\ref{5.1}) hold. Then, there exists a
sufficiently small number }$\eta _{0}>0$\emph{\ such that for any }$\eta \in
(0,\eta _{0})$\emph{\ all functions }$q_{n}^{\delta }\in B\left( R\right)
\cap H^{\delta }$\emph{\ and there exists a number }$\theta =\theta \left(
\eta \right) \in (0,1)$\emph{\ such that the following convergence estimates
are valid }%
\begin{equation*}
\left[ q_{n}^{\delta }-q_{\mathrm{min}}^{\delta }\right] \leq \theta ^{n}%
\left[ q_{0}^{\delta }-q_{\mathrm{min}}^{\delta }\right] ,
\end{equation*}%
\begin{equation*}
\left\Vert c_{n}^{\delta }-c_{\min }^{\delta }\right\Vert _{H^{4}\left(
\epsilon ,M\right) }\leq C_{1}\theta ^{n}\left[ q_{0}^{\delta }-q_{\mathrm{%
min}}^{\delta }\right] ,
\end{equation*}%
\begin{equation*}
\left\Vert q^{\ast }-q_{n}^{\delta }\right\Vert _{H^{1}\left( \Omega
_{T_{0}}\right) }\leq C\delta ^{\rho _{1}}+\theta ^{n}\left[ q_{0}^{\delta
}-q_{\mathrm{min}}^{\delta }\right] ,\text{ }
\end{equation*}%
\begin{equation*}
\left\Vert q^{\ast }\left( x,0\right) -q_{n}^{\delta }\left( x,0\right)
\right\Vert _{H^{1}\left( \epsilon ,M\right) }\leq C\delta ^{\rho
_{2}}+\theta ^{n}\left[ q_{0}^{\delta }-q_{\mathrm{min}}^{\delta }\right] ,
\end{equation*}%
\begin{equation*}
\left\Vert c^{\ast }\left( x\right) -c_{\min }^{\delta }\left( x,0\right)
\right\Vert _{H^{1}\left( \epsilon ,M\right) }\leq C\delta ^{\rho
_{2}}+C_{1}\theta ^{n}\left[ q_{0}^{\delta }-q_{\mathrm{min}}^{\delta }%
\right] ,
\end{equation*}%
\emph{where the constant }$C=C\left( R,T,\epsilon ,\underline{c},\overline{c}%
,\alpha ,M\right) >0$\emph{\ depends only on listed parameters.}

\subsection{The algorithm}

\label{sec 5.4}

Theorems 5.1-5.3 suggest the Algorithm \ref{alg} to solve Problem \ref{cip}.
These theorems rigorously guarantee that Algorithm \ref{alg} globally
converges to a good approximation of the exact solution $c^{\ast }\left(
x\right) $ of Problem \ref{cip}, as long as the level of noise in the data
is sufficiently small. For brevity, we drop the symbol $\delta $ in the
description of this algorithm. In our numerical studies, we choose
parameters $\lambda ,\alpha $ and $\eta $ by a trial and error procedure
only for one test, which we call \textquotedblleft reference test". Next, we
use the same values of these parameters for all other tests.

\begin{algorithm}
	\caption{\label{alg} A numerical method to solve Problem \ref{cip}}
	\begin{algorithmic}[1]
	%\State \label{choice} Choose $\lambda$, $\alpha$ and $\eta$ by a trial and error process. 
	\State  \label{choose q0} Set $n = 0$ and choose a function $q^{0}$ in $B(B/3) \cap H$. 
	\State \label{step min} Minimize
the functional $J_{\lambda ,\alpha ,\beta }$ subject to boundary constraints
  (\ref{q sys}) using the gradient descent method. 
	Denote the obtained minimizer $q_{\rm comp}.$
	\State \label{step 4} Set $c_{\rm comp}(x) = \frac{1}{16 q_{\rm comp}^4(x, 0)}$ for $x \in [\epsilon, M].$
 	\end{algorithmic}
\end{algorithm}

\section{Numerical Studies with Computationally Simulated Data}

\label{sec sim}

In this section, we describe our numerical implementation of the above
Algorithm \ref{alg}, including our strategy to choose the initial solution $%
q^{(0)}$ in Step \ref{choose q0} of Algorithm \ref{alg}. We also present
some details of finding the minimizer of $J_{\lambda ,\alpha ,\beta }$ in
Step \ref{step min}. In addition, to illustrate the efficiency of our
method, we describe some numerical results for computationally simulated
data.

\subsection{Data generation}

\label{sec 6.1}

To generate the data for the forward problem, we use absorbing boundary
conditions (\ref{abso}), (\ref{abso1}) and, therefore, replace problem %
(\ref{main eqn}) with the following one, which we solve numerically: 
\begin{equation}
\left\{ 
\begin{array}{rcll}
c(x)u_{tt}(x,t) & = & u_{xx}(x,t) & (x,t)\in (-a,a)\times (0,T), \\ 
u(-a,t)-u_{x}(-a,t) & = & 0 & t\in (0,T), \\ 
u(a,t)+u_{x}(a,t) & = & 0 & t\in (0,T), \\ 
u(x,0) & = & 0 & x\in \mathbb{R}, \\ 
u_{t}(x,0) & = & \widetilde{\delta }_{0}(x) & x\in \mathbb{R},%
\end{array}%
\right.  \label{6.1}
\end{equation}%
where $a=5$, $T=6$ and $\widetilde{\delta }_{0}(x)=\frac{30}{\sqrt{2\pi }}%
e^{-\frac{(30x)^{2}}{2}}$ is a smooth approximation of the Dirac function $%
\delta _{0}$. We solve problem (\ref{6.1}) by the implicit finite difference
scheme. We choose the implicit scheme because it is much more stable than
the explicit method. In the finite differences, we arrange a uniform
partition for the interval $[-a,a]$ as $\{y_{0}=-a,y_{1},\dots
,y_{N}=a\}\subset \lbrack -a,a]$ with $y_{i}=a+2ia/N_{x}$, $i=0,\dots ,N_{x}$%
, where $N_{x}$ is a large number. In the time domain, we split the interval 
$[0,T]$ into $N_{t}+1$ uniform sub-intervals $[t_{j},t_{j+1}]$, $j=0,\dots
,N_{t}$, with $t_{j}=jT/N_{t},$ where $N_{t}$ is a large number. In our
computational setting, $N_{x}=3000$ and $N_{t}=300$.

By matching the absorbing boundary conditions in (\ref{6.1}) and the
absorption conditions (\ref{abso})--(\ref{abso1}) in Lemma 2.2, we see that
the solution of problem (\ref{main eqn}) can be approximated on $%
[-a,a]\times \lbrack 0,T]$ by the solution of problem (\ref{6.1}). However,
since the Dirac function is replaced by the function $\widetilde{\delta }%
_{0} $, there is a computational error in the computed function $u$ near $%
(x=0,t=0)$. It follows from the presentation (\ref{2.3}) that when $x$ is in
a small neighborhood of $\left\{ x=0\right\} $, where $c(x)=1$, the function 
$u(x,t)=1/2$ if $t<|\tau (x)|$. One can see in Figure \ref{fig 1a} that $%
u(0,t)<\frac{1}{2}$ when $t$ small. Therefore, we simply correct this error
by reassigning $u(x,t)=\frac{1}{2}$ when $x$ and $t$ are near $0$. The
function $u(0,t)$ after this data correction process is displayed in Figure %
\ref{fig 1b}. In our computational program, we set $u(x,t)=\frac{1}{2}$ when 
$(x,t)\in \lbrack 0,0.0067]\times \lbrack 0,0.26].$ 
\begin{figure}[h]
\begin{center}
\subfloat[\label{fig
1a}]{\includegraphics[width=.3\textwidth]{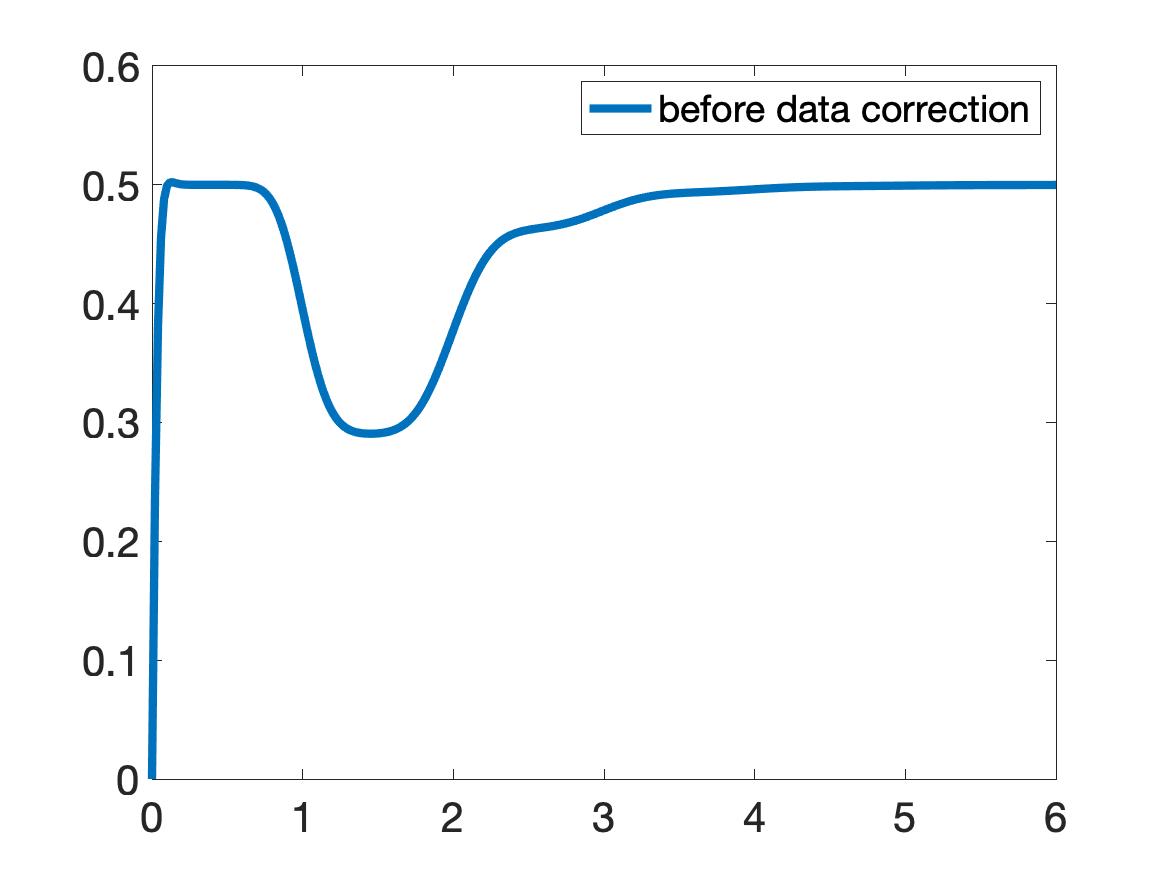}} \quad 
\subfloat[\label{fig
1b}]{\includegraphics[width=.3\textwidth]{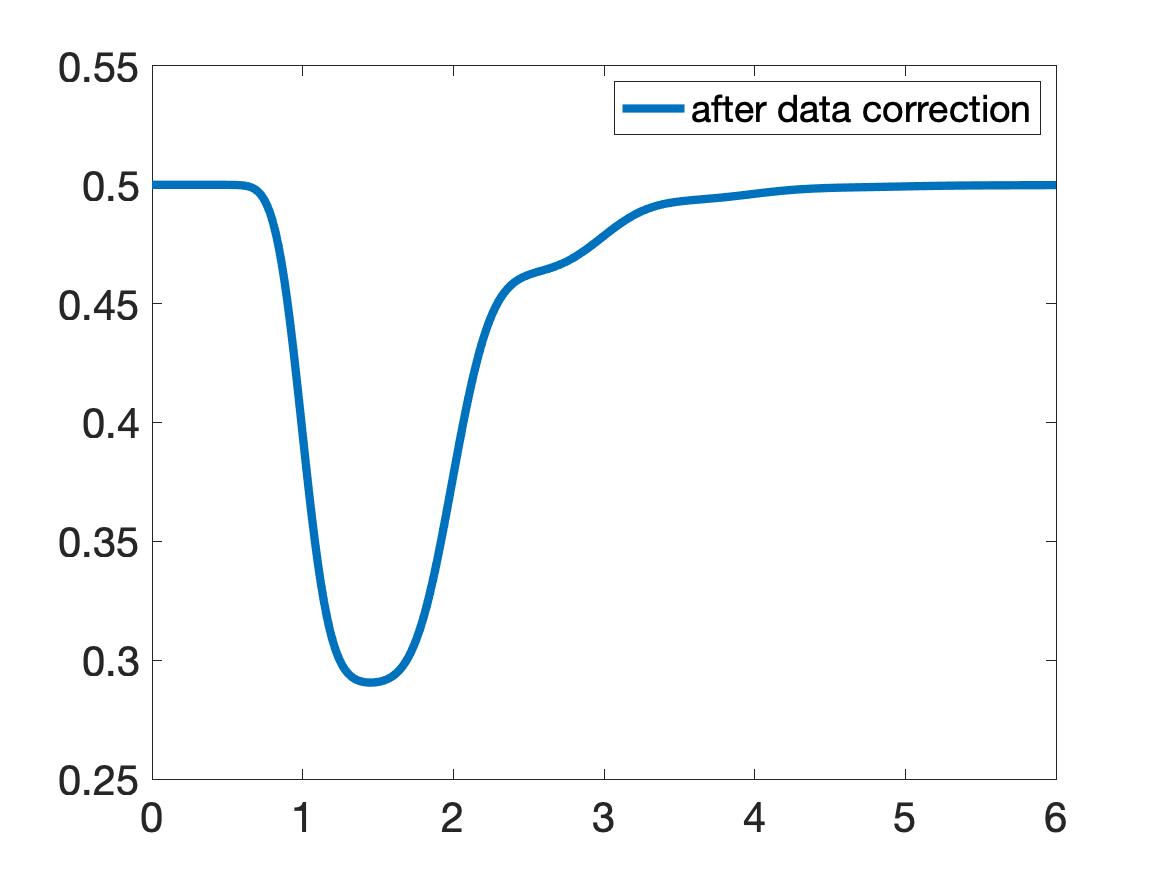}}
\end{center}
\caption{\textit{Illustration of the process of correcting the data near $%
(x=0,t=0).$ We know that when $x$ is small, $c(x)=1$. Therefore by 
(\ref{2.3}), $u(x,t)=\frac{1}{2}$ for $t<|\protect\tau (x)|.$ We, therefore,
set $u(x,t)=\frac{1}{2}$ for $x$ and $t$ small. Figure \protect\ref{fig 1a}
and figure \protect\ref{fig 1b} are the graphs of the function $u(0,t)$
before and after, respectively, this reassignment. These functions are taken
from Test 3 in subsection \protect\ref{sec illus}.}}
\label{fig data correction}
\end{figure}

Having the function $u$ in hands, we compute the functions $g_{0}(t)=u(0,t)$
and $g_{1}(t)=u_{x}(0,t)$ easily. These are the computationally simulated
data for the inverse problem. In the next section, we present the
implementation of the convexification method to solve problem (\ref{q sys}).

\begin{remark}
Let $\delta >0$ be the noise level. We arrange the noisy data by 
\begin{equation*}
g_{0}^{\delta }=g_{0}(1+\delta \mbox{rand})\quad \mbox{and}\quad
g_{1}^{\delta }=g_{1}(1+\delta \mbox{rand}),
\end{equation*}%
where rand is the function that generates uniformly distributed random
numbers in the range $[-1,1].$ The boundary constraints in (\ref{q sys})
involve the derivative of $g_{0}^{\delta }$. We compute $(g_{0}^{\delta
})^{\prime }$ by the Tikhonov regularization method. The Tikhonov
regularization method is well-known. We, therefore, do not describe this
step here. In all numerical tests of subsection \ref{sec illus}, the noise
level is $\delta =0.05,$ i.e. $5\%.$
\end{remark}

\subsection{The numerical implementation of the convexification method}

\label{sec 6.2}

We first present our way to compute the function $q^{(0)}\in H$ in step \ref%
{choose q0} in Algorithm \ref{alg} to initiate the process of minimizing the
objective functional $J_{\lambda ,\alpha ,\beta }.$ In the case when $%
c\equiv 1$, due to (\ref{2.3}), the function $q(x,0)=1/2$ for all $x\in
\lbrack \epsilon ,M].$ Hence, it is natural to set $q^{(0)}(x,0)=1/2$ for
all $x\in \lbrack \epsilon ,M].$ We next find $q^{(0)}(x,t)$ for $t>0$ by
solving the linear partial differential equation obtained by removing the
third term in the left hand side of the equation in (\ref{q sys}). More precisely, we set $q^{(0)}$ as the solution of: 
\begin{equation}
\left\{ 
\begin{array}{ll}
q_{xx}^{(0)}(x,t)-2{q_{xt}^{(0)}(x,t)}=0, & (x,t)\in (\epsilon ,M)\times
(0,T), \\ 
q^{(0)}(\epsilon ,t)=g_{0}(t+\epsilon ), & t\in \lbrack 0,T], \\ 
q_{x}^{(0)}(\epsilon ,t)=g_{1}(t+\epsilon )+g_{0}^{\prime }(t+\epsilon ), & 
t\in \lbrack 0,T], \\ 
q_{x}^{(0)}(M,t)=0, & t\in \lbrack 0,T].%
\end{array}%
\right.  \label{6.2}
\end{equation}%
Denote $Q^{(0)}(x,t)=q_{x}^{(0)}(x,t)$. It follows from (\ref{6.2}) that 
\begin{equation}
\left\{ 
\begin{array}{ll}
Q_{x}^{(0)}(x,t)-2Q_{t}^{(0)}(x,t)=0, & (x,t)\in (\epsilon ,M)\times (0,T),
\\ 
Q^{(0)}(\epsilon ,t)=g_{1}(t+\epsilon )+g_{0}^{\prime }(t+\epsilon ), & t\in
\lbrack 0,T], \\ 
Q^{(0)}(M,t)=0, & t\in \lbrack 0,T].%
\end{array}%
\right.  \label{6.3}
\end{equation}%
The equation in (\ref{6.3}) is a linear transport equation for $Q^{(0)}$
with constant coefficients. Since boundary value problem (\ref{6.3}) is
over-determined, we solve it by the quasi-reversibility method, which was
first introduced in \cite{LattesLions:e1969}. More precisely, we minimize
the functional $I(Q),$ 
\begin{equation*}
I(Q)=\int_{\Omega }|Q_{x}(x,t)-2Q_{t}(x,t)|^{2}dxdt+\int_{0}^{T}|Q(\epsilon
,t)-g_{1}(t+\epsilon )+g_{0}^{\prime }(t+\epsilon )|^{2}dt
\end{equation*}%
\begin{equation}
+\int_{0}^{T}|Q(M,t)|^{2}dt+\eta \Vert Q\Vert _{H^{2}(\Omega _{T})},
\label{702}
\end{equation}%
where $\eta $ is a small number for $Q^{(0)}$. In our computations, $\eta
=10^{-11}.$

We draw the reader's attention to the survey \cite{Klibanov:anm2015} about
the quasi-reversibility method for the existence and uniqueness of the
minimizers of similar functionals as well as for the convergence theorems of
the minimizers to the exact solutions. This method was considered in \cite%
{Klibanov:anm2015} for a variety of ill-posed problems, including
overdetermined ones and for a variety of PDEs. Considerations for (\ref{702})
are quite similar. Thus, we do not discuss these questions here for brevity.

Having $Q^{(0)}(x,t)$ at hands, we find the function $q^{(0)}(x,t)$ as: 
\begin{equation*}
q^{(0)}(x,t)=q^{(0)}(\epsilon ,t)+\int_{\epsilon
}^{x}Q^{(0)}(y,t)dy=g_{0}(t+\epsilon )+\int_{\epsilon }^{x}Q^{(0)}(y,t)dy,
\end{equation*}%
for all $(x,t)\in \lbrack \epsilon ,M]$. Following (\ref{q2c}), we set the
corresponding approximation for the unknown coefficient $c\left( x\right) $
as: 
\begin{equation}
c_{\mathrm{init}}(x)=\frac{1}{16(q^{0})^{4}(x,0)},\quad \mbox{for all }x\in
\lbrack \epsilon ,M].  \label{cinit}
\end{equation}

We now describe our implementation for Step \ref{step min} in Algorithm \ref%
{alg}. This is to minimize the functional $J_{\lambda ,\alpha ,\beta }$. To
work without the boundary constraints in (\ref{q sys}) and to speed up
computations the process, we add to the functional $J_{\lambda ,\alpha
,\beta }$ the boundary terms and minimize the resulting function without
boundary constraints. In addition, although in our theory we use the $%
H^{5}\left( \Omega _{T}\right) -$norm for the regularization term, in
computations we use the simpler to implement $H^{2}\left( \Omega _{T}\right)
-$norm. The resulting functional is still named as $J_{\lambda ,\alpha
,\beta }$, and it is given by%
\begin{equation*}
J_{\lambda ,\alpha ,\beta }(q)=\int_{\Omega _{T}}e^{-2\lambda (x+\alpha t)}%
\Big|q_{xx}(x,t)-q_{xt}(x,t)\frac{1}{2q^{2}(x,0)}+q_{t}(x,t)\frac{q_{x}(x,0)%
}{2q^{3}(x,0)}\Big|^{2}dxdt
\end{equation*}%
\begin{equation*}
+\int_{0}^{T}e^{-2\lambda (x+\alpha t)}|q(\epsilon ,t)-g_{0}(t+\epsilon
)|^{2}dt
\end{equation*}%
\begin{equation}
+\int_{0}^{T}e^{-2\lambda (x+\alpha t)}\left\vert q_{x}(\epsilon
,t)-g_{1}(t+\epsilon )-g_{0}^{\prime }(t+\epsilon )\right\vert ^{2}dt+\beta
\Vert q\Vert _{H^{2}(\Omega _{T})}^{2}.  \label{250}
\end{equation}%
Functional (\ref{250}) can be minimized via a number of optimization
packages. We do so using the ready-to-use optimization toolbox of Matlab.
More precisely, we use the command \textquotedblleft fminunc" of Matlab to
find the minimizer of $J_{\lambda ,\alpha ,\beta }$. The command
\textquotedblleft fminunc" has its own stopping criteria. In our experience
that this command stops when either

\begin{enumerate}
\item \label{case 1} Either a minimizer is found (Matlab lets us know if a
minimizer is found).

\item \label{case 2} Or the number of times Matlab computes the objective
function reaches a default maximum number determined by Matlab.
\end{enumerate}

In the case \ref{case 1}, we take the output of \textquotedblleft fminunc"
as the function $q_{\mathrm{comp}}$ and compute $c_{\mathrm{comp}}$ as in
Step \ref{step 4} of Algorithm \ref{alg}. If \textquotedblleft fminunc"
stops due to the reason of case \ref{case 2}, we understand that the
minimizer is not yet reached. Then, we apply an additional step to speed up
the process. Let $\tilde{q}\left( x,t\right) $ denotes the output of
\textquotedblleft fminunc". We set%
\begin{equation}
\tilde{c}\left( x\right) =\frac{1}{16\tilde{q}^{4}\left( x,0\right) },\quad %
\mbox{for all }x\in \lbrack \epsilon ,M].  \label{6.5555}
\end{equation}%
Next, we solve the following linear boundary value problem with
over-determined boundary data 
\begin{equation}
\left\{ 
\begin{array}{ll}
\tilde{q}_{xx}^{(1)}(x,t)-\frac{\tilde{q}_{xt}^{(1)}(x,t)}{2\tilde{q}%
^{2}(x,0)}+\frac{\tilde{q}_{t}(x,t)\tilde{q}_{x}(x,0)}{2\tilde{q}^{3}(x,0)}%
=0, & (x,t)\in (\epsilon ,M)\times (0,T), \\ 
\tilde{q}^{(1)}(\epsilon ,t)=g_{0}(t+\epsilon ), & t\in \lbrack 0,T], \\ 
\tilde{q}_{x}^{(1)}(\epsilon ,t)=g_{1}(t+\epsilon )+g_{0}^{\prime
}(t+\epsilon ), & t\in \lbrack 0,T], \\ 
\tilde{q}_{x}^{(1)}(M,t)=0, & t\in \lbrack 0,T].%
\end{array}%
\right.  \label{6.4}
\end{equation}%
for a function $\tilde{q}^{(1)}$. Again, we use the quasi-reversibility
method via minimizing the obvious analog of the functional $I(Q)$ in (\ref%
{702}). Next, we set 
\begin{equation}
\tilde{c}_{1}(x)=\frac{1}{16(\tilde{q}^{(1)})^{4}(x,0)},\quad \mbox{for all }%
x\in \lbrack \epsilon ,M].  \label{6.7777}
\end{equation}%
We next minimize the functional $J_{\lambda ,\alpha ,\beta }$ in (\ref{250})
again by \textquotedblleft fminunc" with the initial input $q^{(0)}=\tilde{q}%
^{(1)}$. If the minimizer is found, then we stop. If, however, it is not
found, then we compute a new function $\tilde{c}(x)$ in (\ref{6.5555}) and
proceed as above. This process stops when $\Vert \tilde{c}-\tilde{c}%
_{1}\Vert _{L^{\infty }(\epsilon ,M)}<10^{-3}.$ The final reconstruction of
the function $c$ is $c_{\mathrm{comp}}\left( x\right) =\tilde{c}\left(
x\right) .$ By our computational experience, we need no more than one (1)
correction (\ref{6.4}) for the initial input.

\textbf{Remark 6.3}.\emph{\ In our computations, the parameters for the
Carleman Weight Function are }$\lambda =2$\emph{, }$\alpha =0.3$\emph{, the
regularization parameter }$\beta =10^{-9}$\emph{. }$T=6$\emph{, }$\epsilon
=0 $\emph{\ and }$M=3$\emph{. In (\ref{606})\ }$\alpha T>2M.$\emph{\ But
this condition is not necessary to impose in our numerical experiments.
These numbers are chosen by a trial-and-error procedure. This means that we
try many sets of these parameters to get the best numerical result for one
test, which we call \textquotedblleft reference test" (test 1 in subsection %
\ref{sec illus}). Then we use the same parameters for all other tests,
including the tests with experimental data. In theory, }$\lambda $\emph{\
should be a large number. Here, we choose }$\lambda =2$\emph{\ because this
value is sufficient to obtain satisfactory numerical results. If }$\lambda $%
\emph{\ is too large, then the Carleman Weight Function decays too rapidly.
This causes many difficulties in numerics; especially, on the computing
time. In fact, a similar issue takes place in any asymptotic theory when it
is applied to real computations. Indeed, such a theory basically says that
\textquotedblleft If a certain parameter X is sufficiently large, then a
certain \textquotedblleft good thing"} \emph{takes place". However, when
computing, one needs to estimate X computationally since theoretical
estimates usually more pessimistic than numerical ones.}

\subsection{Numerical results for computationally simulated data}

\label{sec illus}

\begin{figure}[h!]
\subfloat[\label{fig 2a} Numerical result of test 1 for simulated data with
$5\%$ noise]{\includegraphics[width=.3\textwidth]{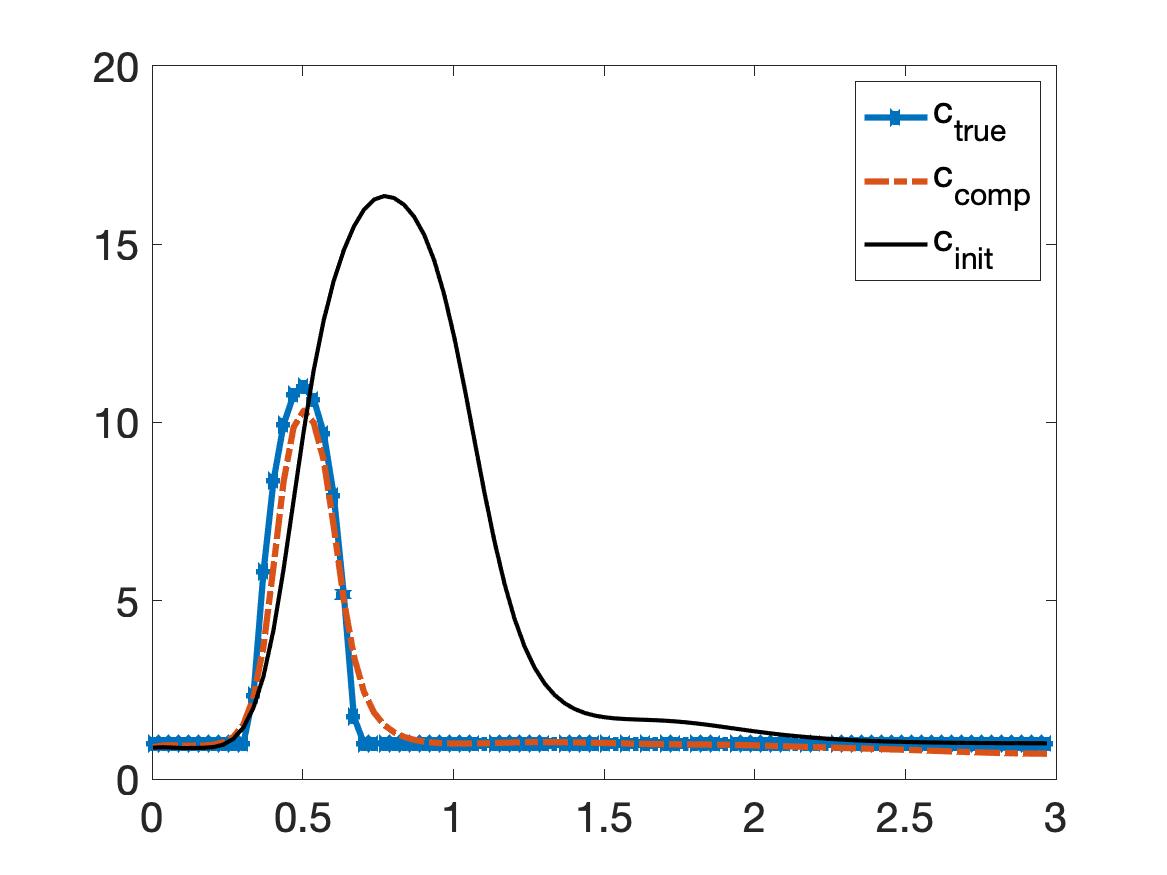}} \quad 
\subfloat[\label{fig 2b} Numerical result for test 2 from simulated data with
$5\%$ noise]{\includegraphics[width=.3\textwidth]{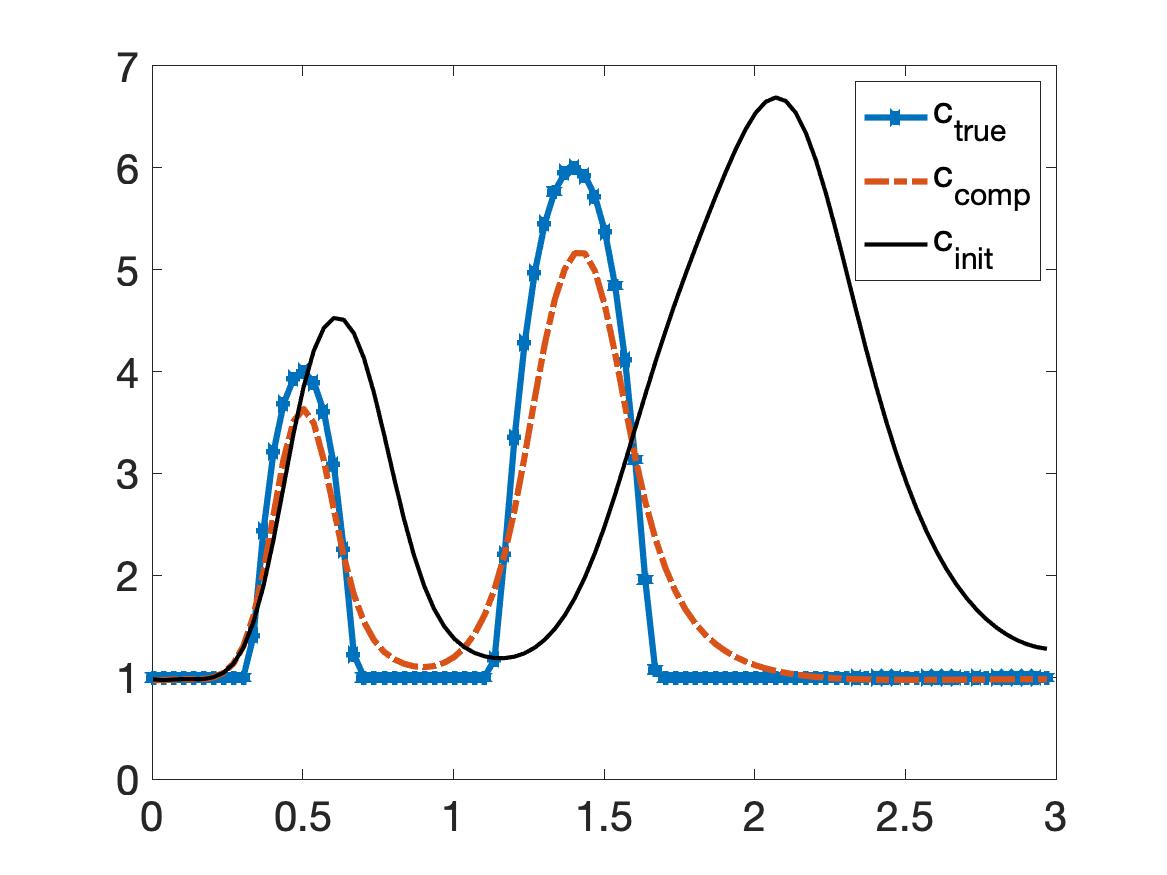}} \quad 
\subfloat[\label{fig 2c} Numerical result for test 3 from simulated data with
$5\%$ noise]{\includegraphics[width=.3\textwidth]{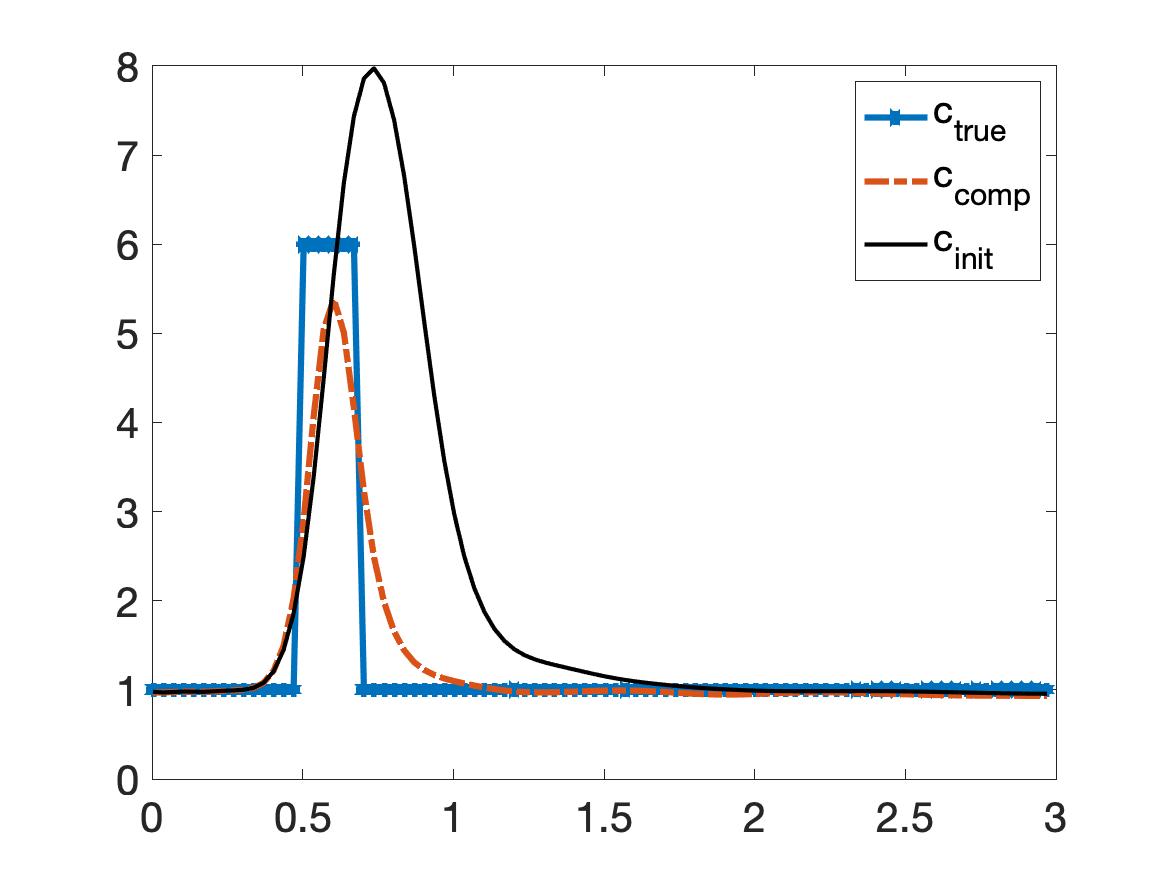}}
\par
\subfloat[\label{fig 2d} Numerical result for test 4 from simulated data with
$5\%$ noise]{\includegraphics[width=.3\textwidth]{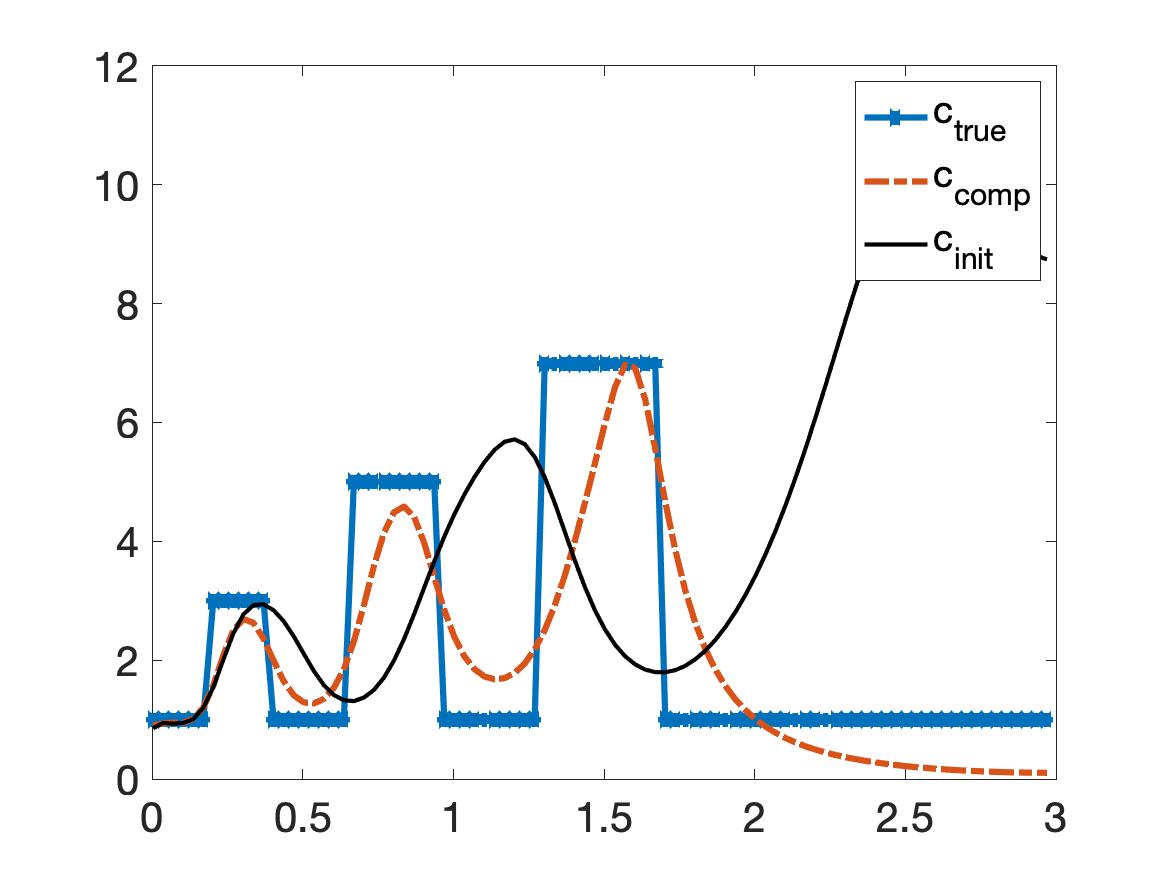}} \quad 
\subfloat[\label{fig 2e} Numerical result for test 5 from simulated data with
$5\%$ noise]{\includegraphics[width=.3\textwidth]{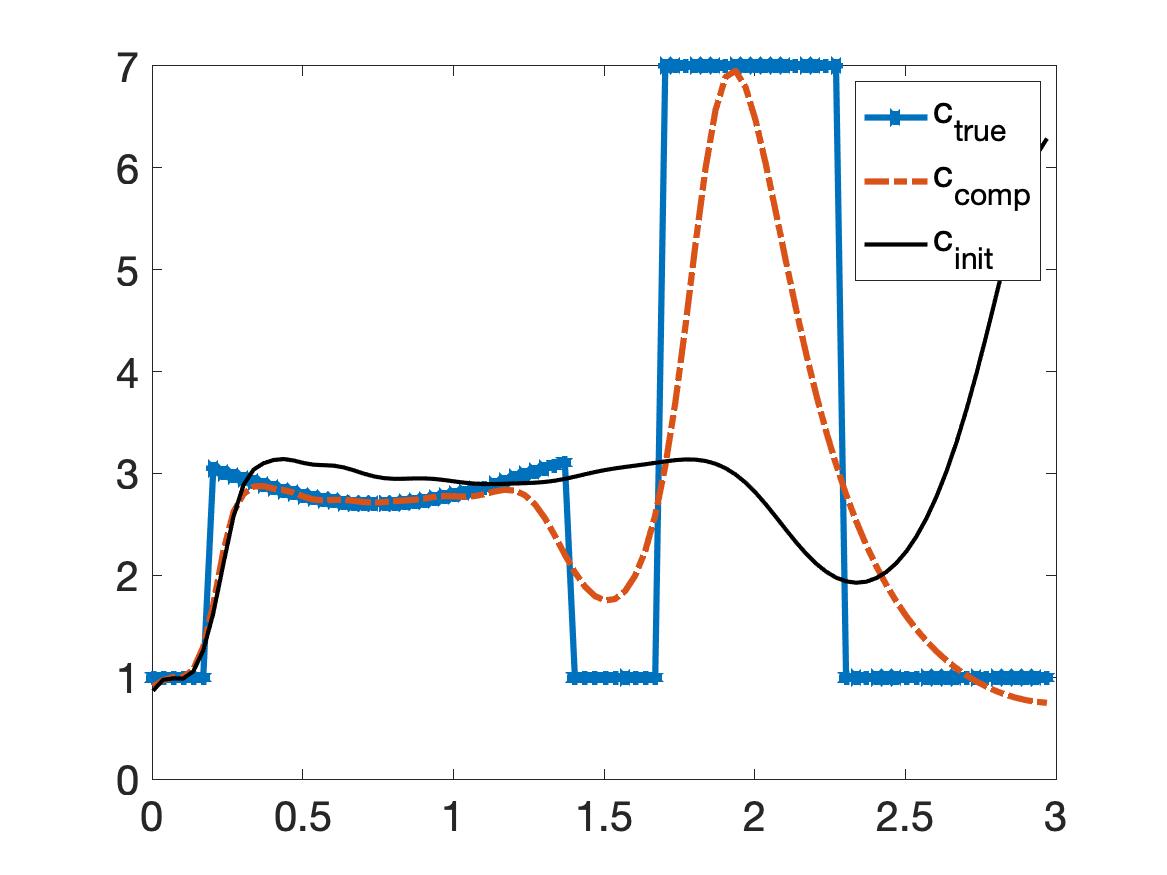}}
\caption{\textit{The true spatially distributed dielectric constant function 
$c_{\mathrm{true}}$, its initial version $c_{\mathrm{init}}$ computed by 
(\ref{cinit}) and its final reconstruction $c_{\mathrm{comp}}$ by our
convexification method. It is evident that in all cases, the initial
solution $c_{\mathrm{init}}$ computed by (\ref{cinit}) already carries some
information of $c_{\mathrm{true}}$. The following iterative steps
significantly improve the positions of ``inclusions" and their values.
Especially, in test 5 (Figure \protect\ref{fig 2e}), the convexification
method successfully reconstructs the curves in the inclusion in the left.}}
\label{fig 3}
\end{figure}

We present five (5) numerical examples to test our convexification method.
The obtained results are displayed in Figure \ref{fig 3}.

\noindent \textit{Test 1.} In this test, we consider the case of one
inclusion with a high inclusion/background contrast. The true function $%
c\left( x\right) $ is given by 
\begin{equation*}
c_{\mathrm{true}}(x)=\left\{ 
\begin{array}{ll}
1+10e^{\frac{(x-0.5)^{2}}{(x-0.5)^{2}-0.2^{2}}} & \mbox{if }|x-0.5|<0.2, \\ 
1 & \mbox{otherwise}.%
\end{array}%
\right.
\end{equation*}%
In this test, we detect one object with a high dielectric constant with the
size 0.4 and the center located at $0.5$. Although the inclusion/background
contrast here is $11/1=11,$ which is high, our method provides good
numerical results without any knowledge of $c_{\mathrm{true}}$ inside of $%
[\epsilon ,M]$. The numerical solution of this test is displayed in Figure %
\ref{fig 2a}. In this test, the function $c_{\mathrm{init}}$ obtained by %
(\ref{cinit}) somewhat provides the information about $c_{\mathrm{true}},$
but it is still far away from $c_{\mathrm{true}}.$ The final reconstruction
quite exactly indicates the position of the \textquotedblleft inclusion".
The maximal value of the computed function $c\left( x\right) $ in the
inclusion is 10.43 (relative error 5.2\%). This value is accurate since we
have the the noise level $\delta =5\%$.

\noindent \textit{Test 2.} We test a more complicated function $c_{\mathrm{%
true}}.$ In this test, the dielectric constant is a smooth function $c\left(
x\right) $ with two (2) inclusions. The function $c_{\mathrm{true}}$ is
given by 
\begin{equation*}
c_{\mathrm{true}}(x)=\left\{ 
\begin{array}{ll}
1+3e^{\frac{(x-0.5)^{2}}{(x-0.5)^{2}-0.2^{2}}} & \mbox{if }|x-0.5|<0.2, \\ 
1+5e^{\frac{(x-1.4)^{2}}{(x-1.4)^{2}-0.3^{2}}} & \mbox{if }|x-1.4|<0.3, \\ 
1 & \mbox{otherwise}.%
\end{array}%
\right.
\end{equation*}%
This test is challenging since the maximal value of the function $c_{\mathrm{%
true}}(x)$ in each inclusion is high (4 and 6). The left inclusion is
blocked by the right inclusion in the in the view of the source and the
detector, both of which are located at $\left\{ x=0\right\} $. The graphs of
the true function $c(x),$ initial and computed solutions are displayed in
Figure \ref{fig 2b}. The function $c_{\mathrm{init}}$ computed by %
(\ref{cinit}) somewhat provides a guess about the shape of $c_{\mathrm{true}%
} $ but is still far away from $c_{\mathrm{true}}.$ The final reconstruction
is good. The computed locations of both inclusions are satisfactory. The
maximal value of the computed function $c\left( x\right) $ in the left
inclusion is 3.40 (relative error 15\%). The maximal value of the computed
function $c(x)$ in the right inclusion is 5.16 (relative error 14\%).

\noindent \textit{Test 3.} We now consider the case when $c_{\mathrm{true}%
}(x)$ is a discontinuous step function, 
\begin{equation*}
c_{\mathrm{true}}(x)=\left\{ 
\begin{array}{ll}
6 & \mbox{if }|x-0.6|<0.1, \\ 
1 & \mbox{otherwise.}%
\end{array}%
\right.
\end{equation*}%
This test is an interesting one. It shows that the convexification method is
stronger than what we can prove in the theory in the sense that the
smoothness condition of the function $c\left( x\right) $ can be relaxed in
numerical studies, although this condition is used in the theoretical part.
The true, initial and the computed solutions of Problem \ref{cip} are
displayed in Figure \ref{fig 2c}. The initial solution obtained by %
(\ref{cinit}) somewhat indicates the inclusion but both the location and the
value of the dielectric constant inside the inclusion are far from correct
ones. However, both the location and the computed dielectric constant meet
the expectation in the final reconstruction of the function $c$. The maximal
value of the computed function $c\left( x\right) $ is 5.6 (relative error is
6.7\%).

\noindent \textit{Test 4.} We consider the case of three inclusions. As in
the previous example, the dielectric constant function $c$ in this case is a
discontinuous one. It is given by 
\begin{equation*}
c_{\mathrm{true}}(x)=\left\{ 
\begin{array}{ll}
3 & \mbox{if }|x-0.3|<0.1, \\ 
5 & \mbox{if }|x-0.8|<0.15, \\ 
7 & \mbox{if }|x-1.5|<0.2, \\ 
1 & \mbox{otherwise}.%
\end{array}%
\right.
\end{equation*}%
Reconstructing this function $c\left( x\right) $ is challenging. In fact,
since we only measure the data at $x=\epsilon $, in the view of the
detector, the second and third inclusions are blocked by the first one.
Nevertheless, our method works well. The numerical solutions are displayed
in Figure \ref{fig 2d}. As in the previous examples, the initial solution $%
c_{\mathrm{init}}\left( x\right) $ provides some information about the
function $c_{\text{true}}\left( x\right) ,$ but the error is large. This
error is corrected by our convexification method. The final reconstruction
successfully shows locations of all three inclusions. The maximal values of
the computed function $c\left( x\right) $ in each inclusion are good. The
computed maximal value of $c\left( x\right) $ in the left inclusion is 2.8
(relative error 6.7\%). The computed maximal value of $c\left( x\right) $ in
the middle inclusion is 4.6 (relative error 8.0\%). The computed maximal
value of $c\left( x\right) $ in the right inclusion is 6.9 (relative error
1.4\%).

\noindent \textit{Test 5.} We now test another interesting case, in which
the true dielectric constant includes two inclusions. The function $c\left(
x\right) $ in the first one is a smooth function, and in the second one it
has a constant value. The true dielectric constant is given by 
\begin{equation*}
c_{\mathrm{true}}=\left\{ 
\begin{array}{ll}
3+0.3\sin (\pi (x-1.25)) & \mbox{if }|x-.08|<0.6, \\ 
7 & \mbox{if }|x-2|<0.3, \\ 
1 & \mbox{otherwise.}%
\end{array}%
\right.
\end{equation*}%
The numerical solution of this test is given in Figure \ref{fig 2e}. The
initial solution $c_{\mathrm{init}}$ obtained by (\ref{cinit}) is far away
from the true function $c\left( x\right) $. It might not contain any
valuable information of the true function $c_{\text{true}}\left( x\right) $.
In the next step, after applying the convexification method, we get a good
reconstruction of $c_{\mathrm{true}}\left( x\right) $. The curve in the
first inclusion locally coincides with the true one. The position and the
computed function $c\left( x\right) $ of the second inclusion are also
accurate. The computed maximal value of $c\left( x\right) $ in the second
inclusion is 6.9 (relative error 1.4\%).

\textbf{Remark 6.4}. \emph{In this section, we have tested our
convexification method for multiple cases. The numerical results show that
our method is robust, since it can be used for the cases when the dielectric
constant has high contrasts, single or multiple inclusions, and a
complicated form. More importantly, we obtain those satisfactory results
without requiring any initial guess.}

\section{Numerical Studies of Experimental Data}

\label{sec exp}

We use the data collected by the Forward Looking Radar built in the US Army
Research Laboratory \cite{NguyenWong:pspie2007}. The goal of this radar is
to detect and identify flash explosive-like targets, such as antipersonnel
land mines and improvised explosive devices. These targets can be both
buried on a few centimeters depth in the ground and located in air, i.e.
above the ground.

The device has an emitter and sixteen (16) detectors. The emitter sends out
only one component of the electric field forward the area that covers the
object and the detectors collect the back scattering electric signal
(voltage) in the time domain. The same component of the electric field is
measured as the one which is generated by the emitter, see Figure \ref{fig
sch} for the schematic diagram of data collection. See the comment about the
validity of the data at the beginning of \cite[Section 2]%
{KlibanovLoc:ipi2016}. The step size in time is 0.133 nanosecond. The
backscattering data in the time domain are collected when the distance
between the radar and the target varies from 8 to 20 meters. We then take
the average of these data with respect to both the position of the radar and
those 16 detectors and use as the 1D data to test our convexification
method. Due to this \textquotedblleft average" of the data, we are unable to
find the location of the target. The location can be found by using the
Ground Positioning System (GPS). The error in each of horizontal coordinates
does not exceed a few centimeters, which is sufficient for practical
purposes. When the target is under the ground, the GPS provides the distance
between the radar and a point on the ground located above the target. As to
the depth of a buried target, it is not of a significant interest, since
horizontal coordinates are known and it is also known that the depth does
not exceed 10 centimeters. We refer to \cite{NguyenWong:pspie2007} for more
details about the data collection process. We refer to previous works of our
group in \cite%
{Karchevskyetal:2013,KlibanovLoc:ipi2016,KlibanovAlex:SIAMjam:2017,KlibanovKolesov:ip2018,Kuzhuget:IP2012, Kuzhuget:IEEE2013,Smirnov:ip2020}
where these experimental data were treated by different inversion algorithms
for Coefficient Inverse Problems.

Hence, the interest here is to compute the values of the dielectric
constants of the targets using these data. Indeed, we hope that in the
future knowledge of dielectric constants, being combined with the knowledge
of other parameters of targets, might help to reduce the false alarm rate.
An interesting feature of our data is that they were collected in the field,
rather than in a simpler case of a laboratory. Besides, all targets were
surrounded by clutter.

As in all previous our above cited works on these data, we have calculated
the relative spatially distributed dielectric constant $c_{\mathrm{rel}}(x)$
of the medium including the background (air or ground) and the target. The
function $c_{\mathrm{rel}}(x)$ is given by 
\begin{equation}
c_{\mathrm{rel}}(x)=\left\{ 
\begin{array}{ll}
\frac{c_{\mathrm{target}}}{c_{\mathrm{bckgr}}} & \mbox{if }x\in D, \\ 
1 & \mbox{otherwise}%
\end{array}%
\right.  \label{7.1}
\end{equation}%
where $D$ is a sub interval of $[\epsilon ,M]$ which is occupied by the
target. Here, $c_{\mathrm{target}}$ is the dielectric constant of the target
and $c_{\mathrm{bckgr}}$ is the dielectric constant of the background. If
the background is air, then $c_{\mathrm{bckgr}}=1$. If the background is dry
sand, then $c_{\mathrm{bckgr}}\in (3,5)$ (see table of dielectric constants
listed on a website of Honeywell, https://goo.gl/kAxtzB). Our inverse solver
in this paper is suitable to compute $c_{\mathrm{rel}}$ given the
backscattering data. The computed $c_{\mathrm{target}}$ follows.

\subsection{Data preprocessing}

It was observed in previous above cited publications of this group about
inversion of these experimental data that there is a significant discrepancy
between the computationally simulated data and experimentally collected
data. Hence, the first step to invert these data is to preprocess them. So
that the preprocessed data and the simulated data would look similarly. We
are doing this by scaling and truncating. We consider two cases.

\begin{figure}[h!]
\begin{center}
\subfloat[\label{fig 4a}The time-dependent raw
data]{\includegraphics[width=.3\textwidth]{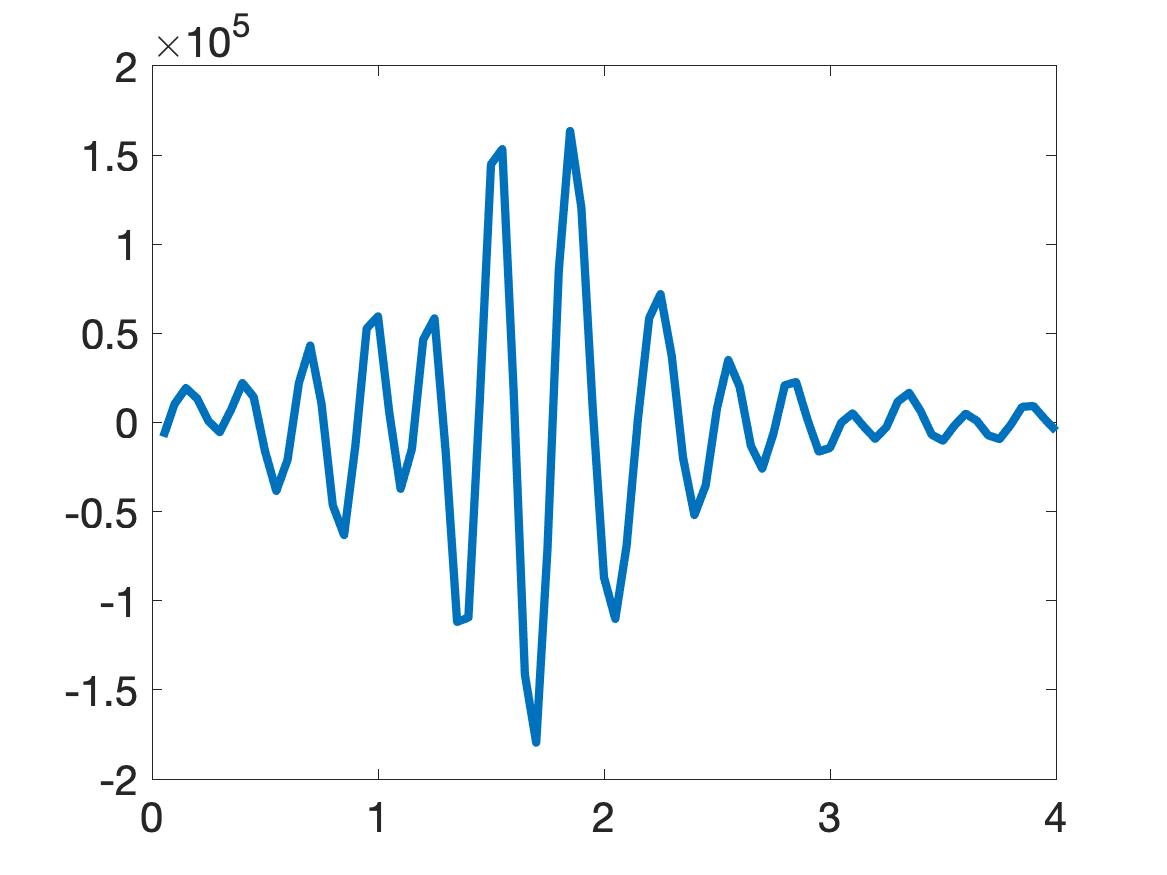}} \quad 
\subfloat[\label{fig 4b}The time-dependent backscattering wave after
preprocessing]{\includegraphics[width=.3\textwidth]{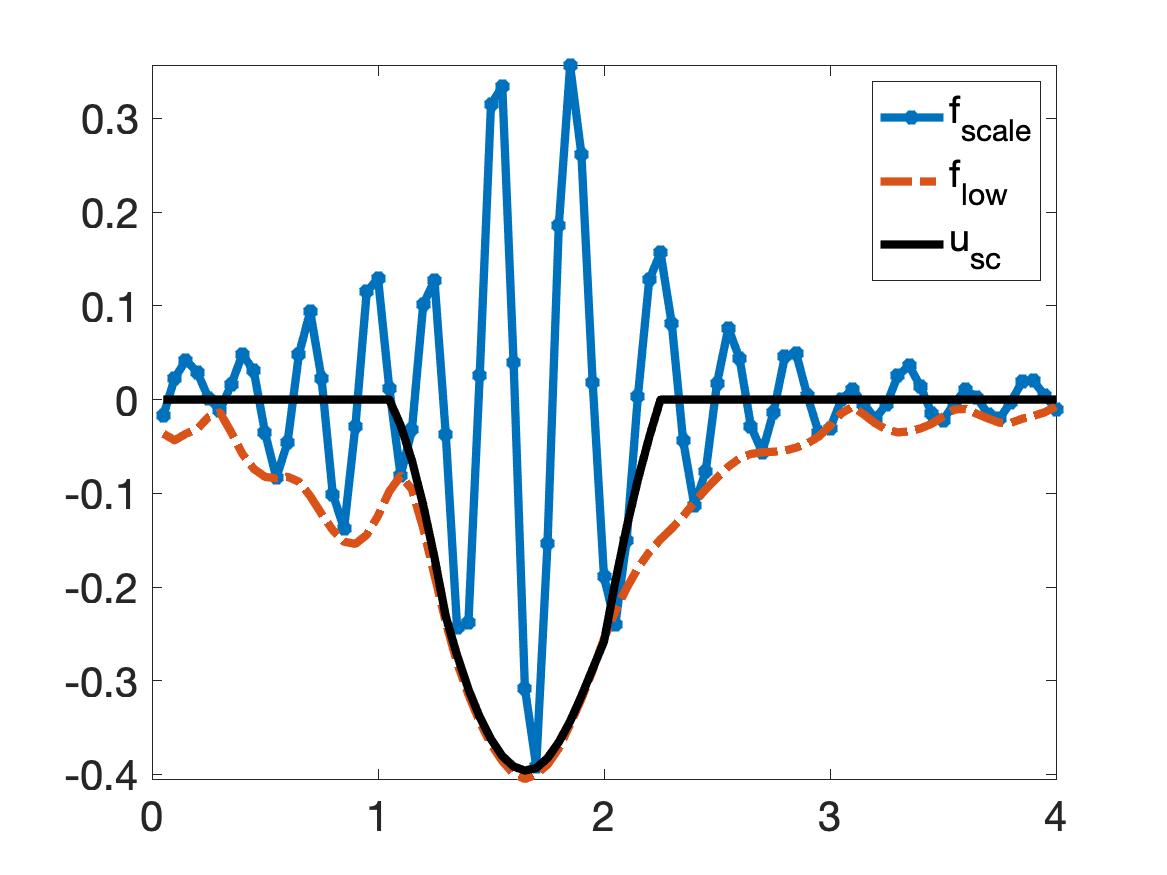}} \quad 
\subfloat[\label{fig 4c}Computed dielectric constant. Its maximal value is
6.76.]{\includegraphics[width=.3\textwidth]{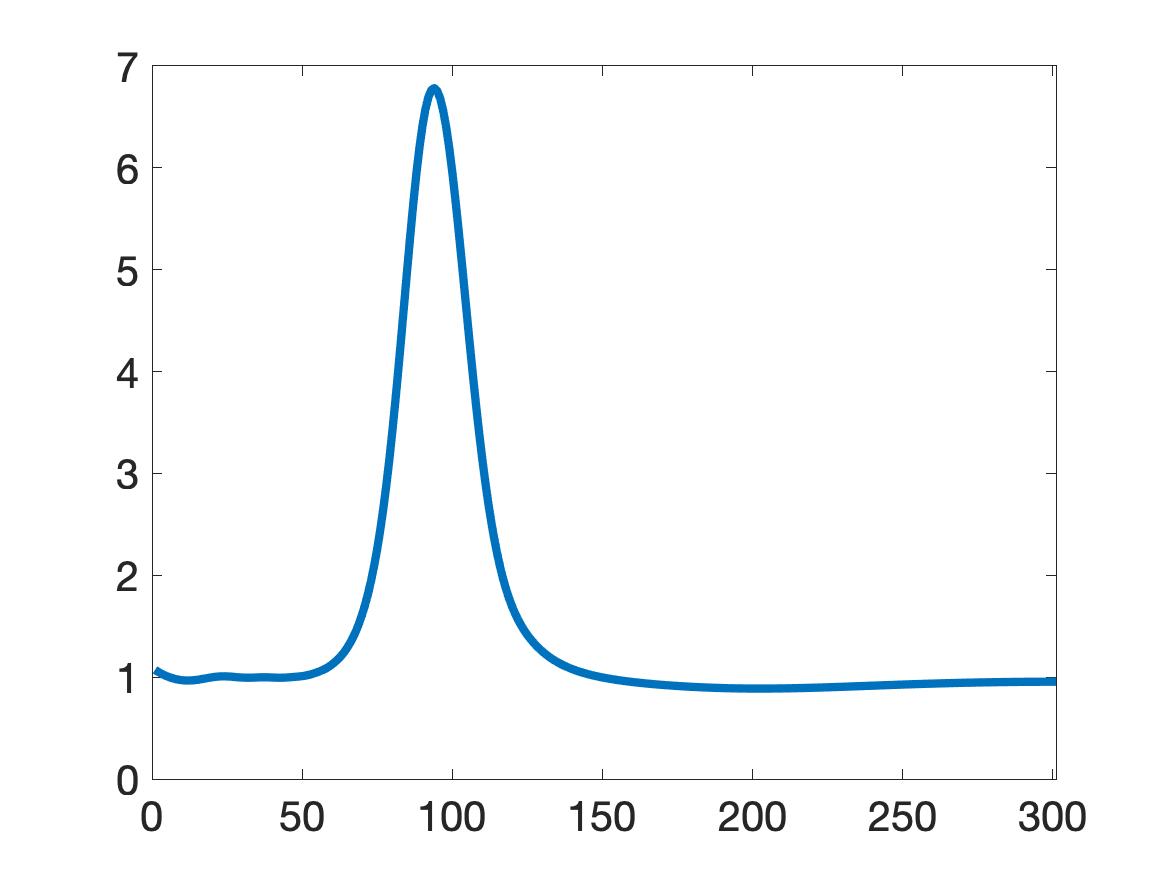} }
\par
\subfloat[\label{fig 4d}The time-dependent raw
data]{\includegraphics[width=.3\textwidth]{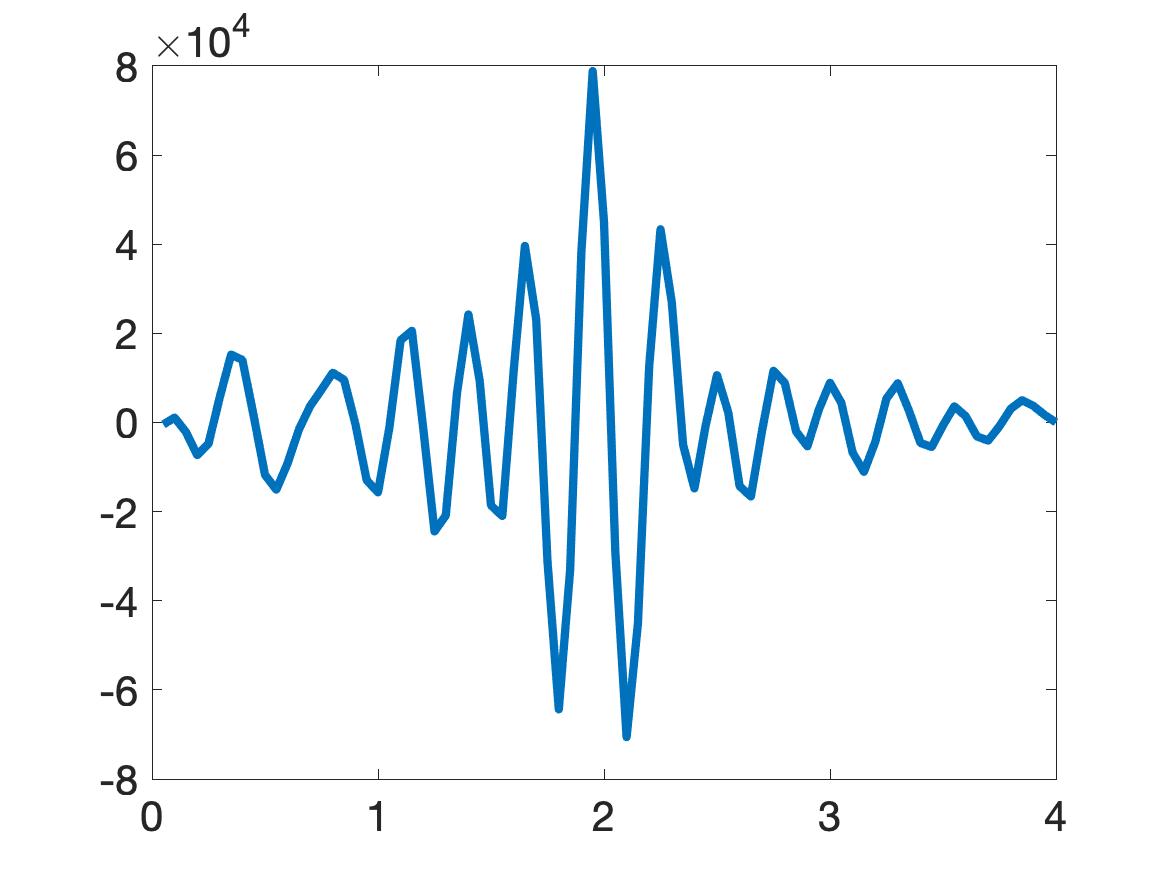}} \quad 
\subfloat[\label{fig 4e}The time-dependent backscattering wave after
preprocessing]{\includegraphics[width=.3\textwidth]{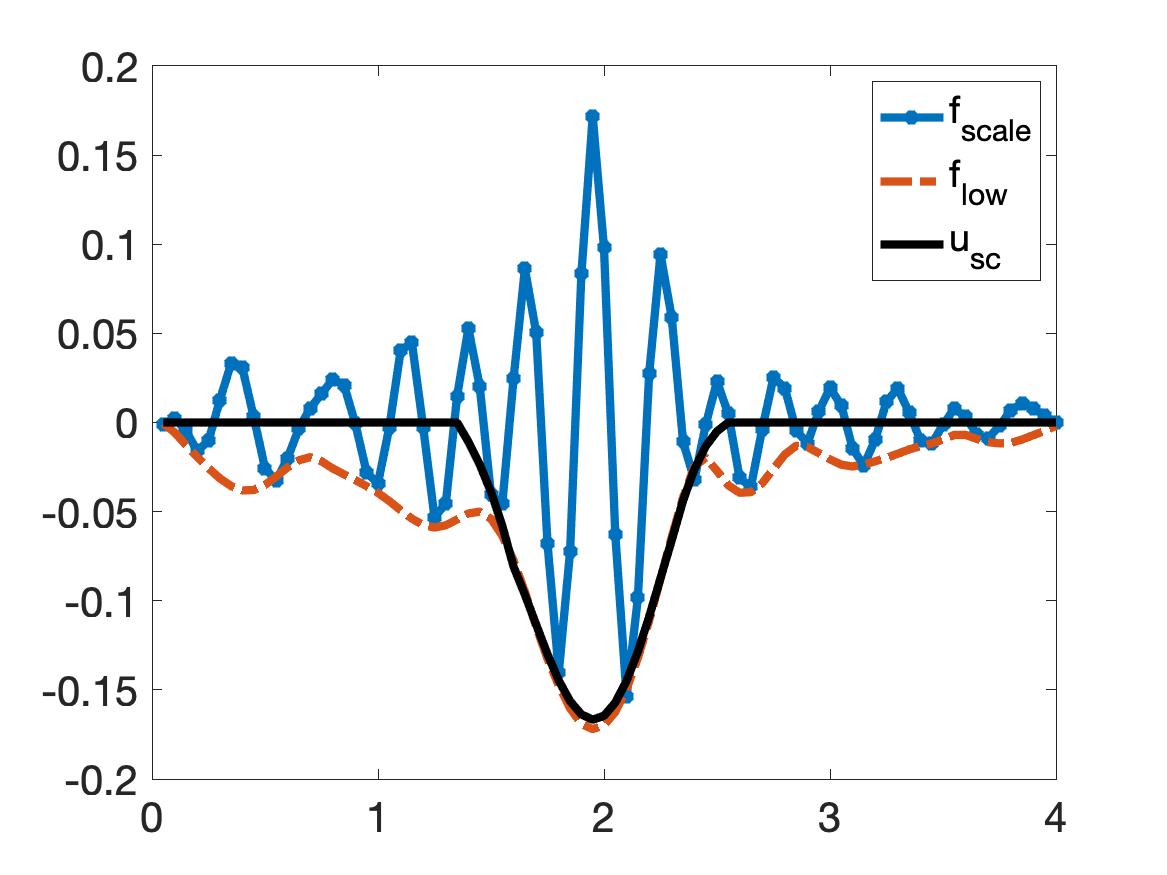}} \quad 
\subfloat[\label{fig 4f}Computed dielectric constant. Its maximal value is
2.2.]{\includegraphics[width=.3\textwidth]{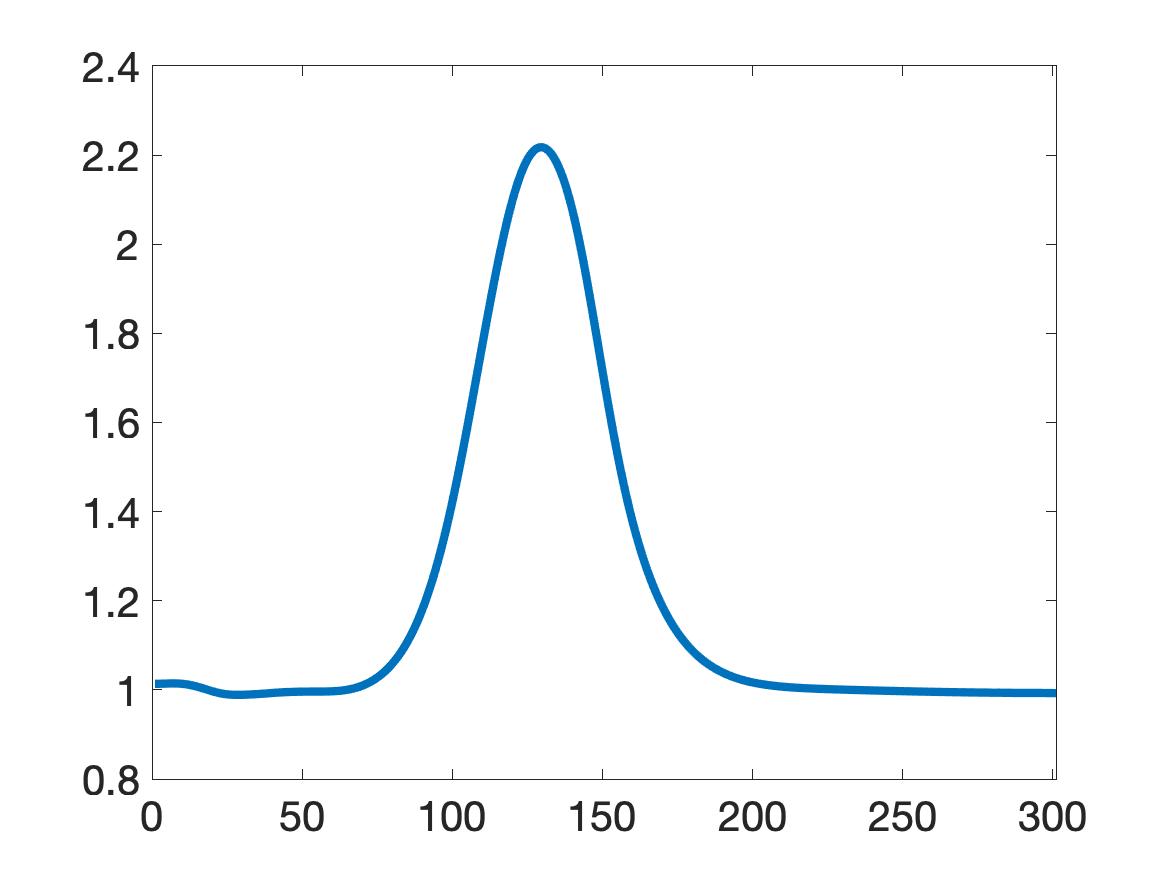} }
\end{center}
\caption{\textit{The case when the target is in the air. The raw and
preprocessed data in the first row correspond to the wave scattered from a bush. The
raw and preprocessed data in the second row correspond to the wave scattered from a
wood stake. The computed dielectric constants for these two tests meet the
expectation since they belong to intervals of their true value, see the last
three rows of Table \protect\ref{tab 1}. }}
\label{fig 4}
\end{figure}

\textit{1. The case when the targets are in air. } We first notice that the
magnitude of the raw experimental data $f_{\mathrm{raw}}$ is large while
that of the simulated data is small. This difference is due to the fact that
we scale the speed of light in the air to be 1. Thus, we compute a
\textquotedblleft scaling factor". To do so, we have to know the true
solution of one set of data generated by a known target. This target is
called the reference object. We choose the reference object as a bush with
its dielectric constant about 6.5 \cite{KlibanovLoc:ipi2016}. We then
generate a corresponding simulated data, named as $f_{\mathrm{sim}}$. The
scaling factor $\mu $ is determined as $\mu \Vert f_{\mathrm{raw}}\Vert
_{L^{\infty }}=\Vert f_{\mathrm{sim}}\Vert _{L^{\infty }}$. The computed
scaling factor is $\mu =459420$. We use the same scaling factor for other
tests. The scaled data $f_{\mathrm{scale}}=\mu f_{\mathrm{raw}}.$ We next
truncate the data. Since the object is placed in the air, then $c_{\mathrm{%
bckgr}}=1$ and $c_{\mathrm{target}}>1$. As seen in Figure \ref{fig data
correction}, the value of the value of the total simulated wave is less than 
$0.5:$ 0.5 should be subtracted, see the first term in the right hand side of %
(\ref{2.3}). This term is responsible for the incident wave. Therefore, the
back scattering wave is non-positive. We thus cut off all positive values of 
$f_{\mathrm{scale}}$ by bounding it by a lower envelop for the scaled data,
see Figures \ref{fig 4b} and \ref{fig 4e} for illustrations of the envelops.
This lower envelop is the graph of the function named $f_{\mathrm{low}}(t)$.
We next truncate $f_{\mathrm{low}}$ because we known that before and after
the backscattering wave hits and then passes the detector, the data is $0$.
This truncating step is as follows. Let $t_{\mathrm{min}}$ be the absolute
minimizer of $f_{\mathrm{low}}(t)$. We keep the value of $f_{\mathrm{low}%
}(t) $ in a neighborhood of $t_{\mathrm{min}}$, say $(t_{\mathrm{min}%
}-10\delta _{t},t_{\mathrm{min}}+10\delta _{t})$ where $\delta _{t}$ is the
step size in time, and re-assign the value of $f_{\mathrm{low}}(t)=0$
outside this neighborhood. The obtained function is the backscattering wave $%
u_{\mathrm{sc}}$. Due to Lemma \ref{lem 2.1}, the total wave at the detector
is $u_{\mathrm{sc}}+0.5.$ See Figures \ref{fig 4a}, \ref{fig 4b}, \ref{fig
4d} and \ref{fig 4e} for the results of data preprocessing.

\begin{figure}[h!]
\begin{center}
\subfloat[\label{fig 5a}The time-dependent raw
data]{\includegraphics[width=.3\textwidth]{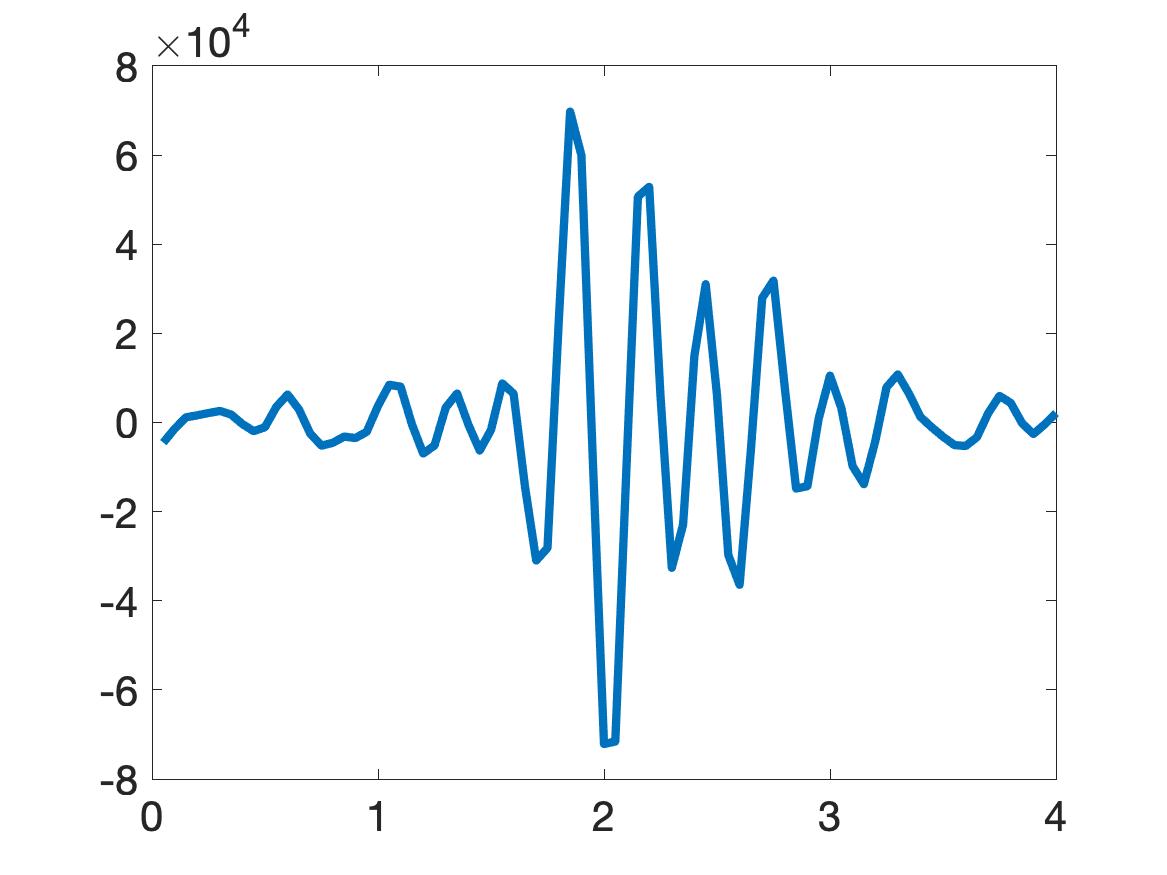}} \quad 
\subfloat[\label{fig 5b}The time-dependent backscattering wave after
preprocessing]{\includegraphics[width=.3\textwidth]{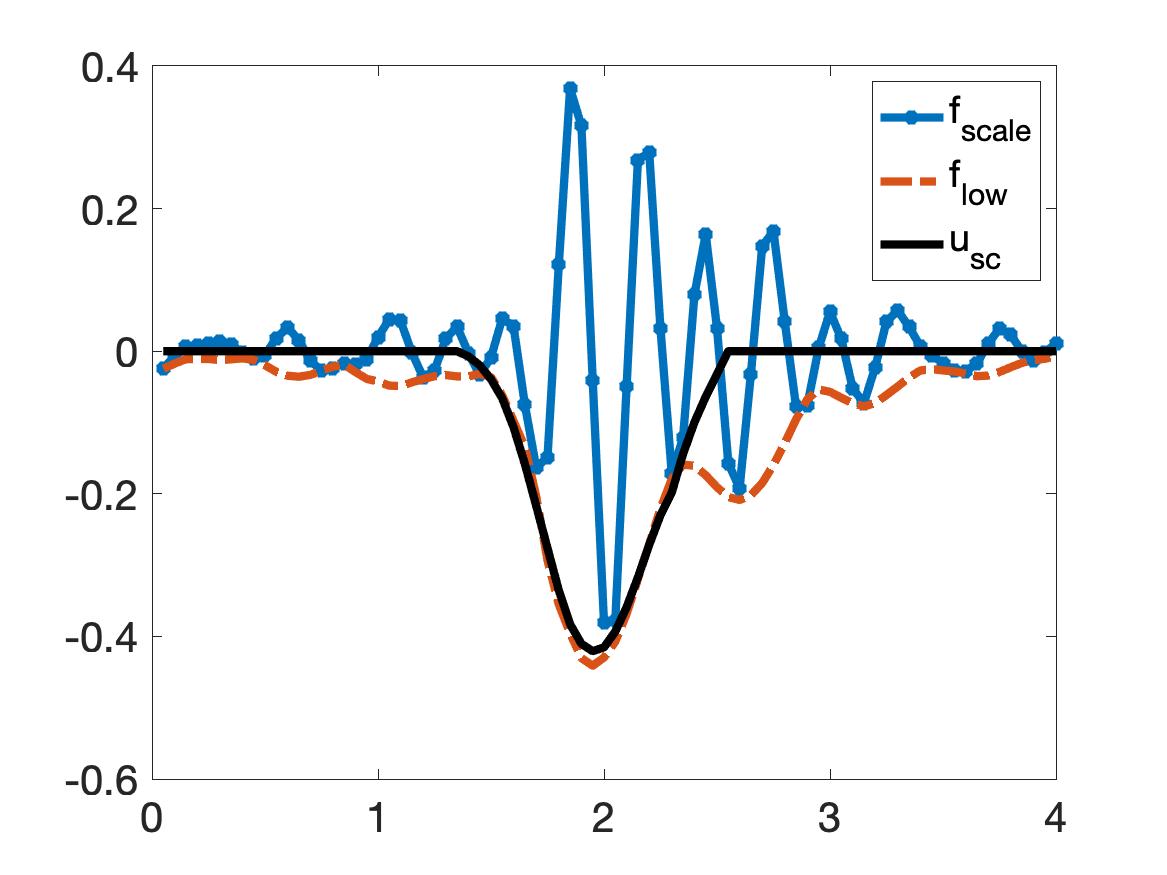}} \quad 
\subfloat[\label{fig 5c}Computed dielectric constant. Its maximal value is
5.2.]{\includegraphics[width=.3\textwidth]{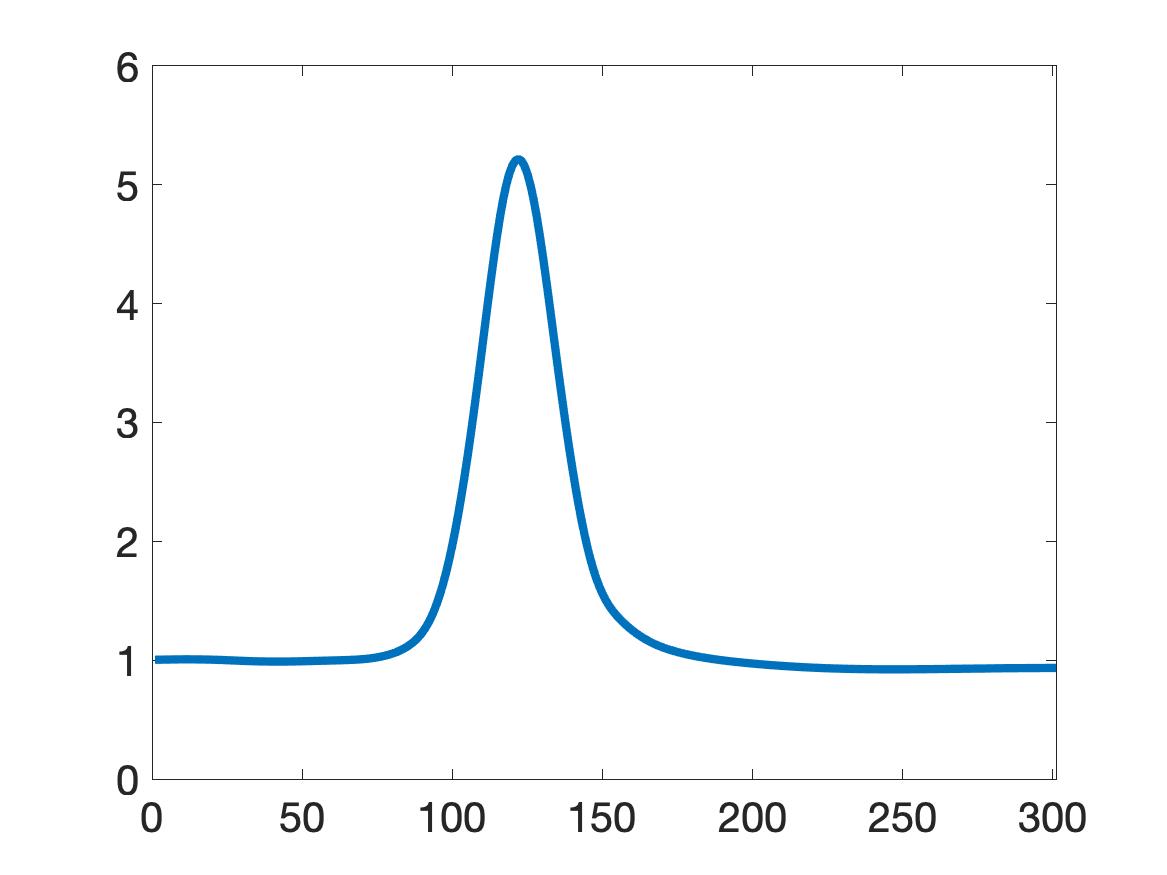} }
\par
\subfloat[\label{fig 5d}The time-dependent raw
data]{\includegraphics[width=.3\textwidth]{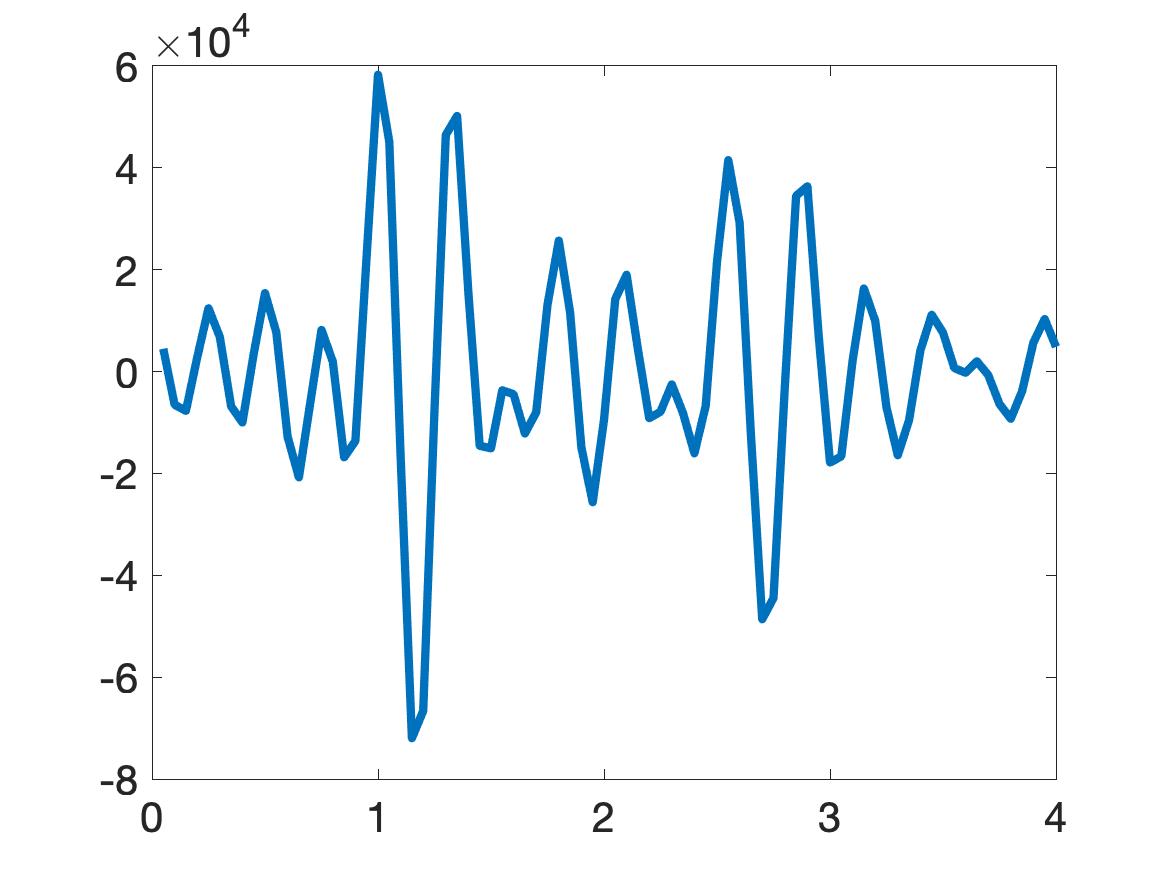}} \quad 
\subfloat[\label{fig 5e}The time-dependent backscattering wave after
preprocessing]{\includegraphics[width=.3\textwidth]{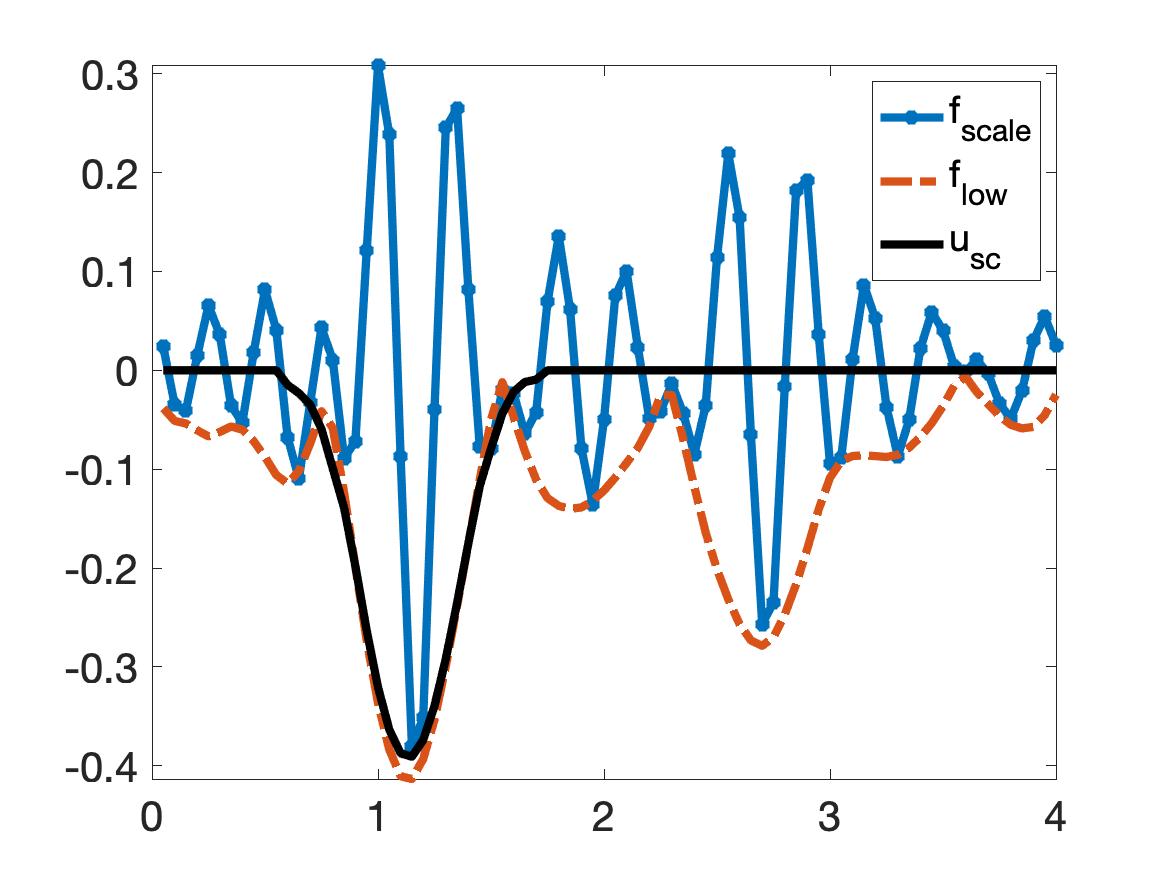}} \quad 
\subfloat[\label{fig 5f}Computed dielectric constant. Its maximal value is
4.7.]{\includegraphics[width=.3\textwidth]{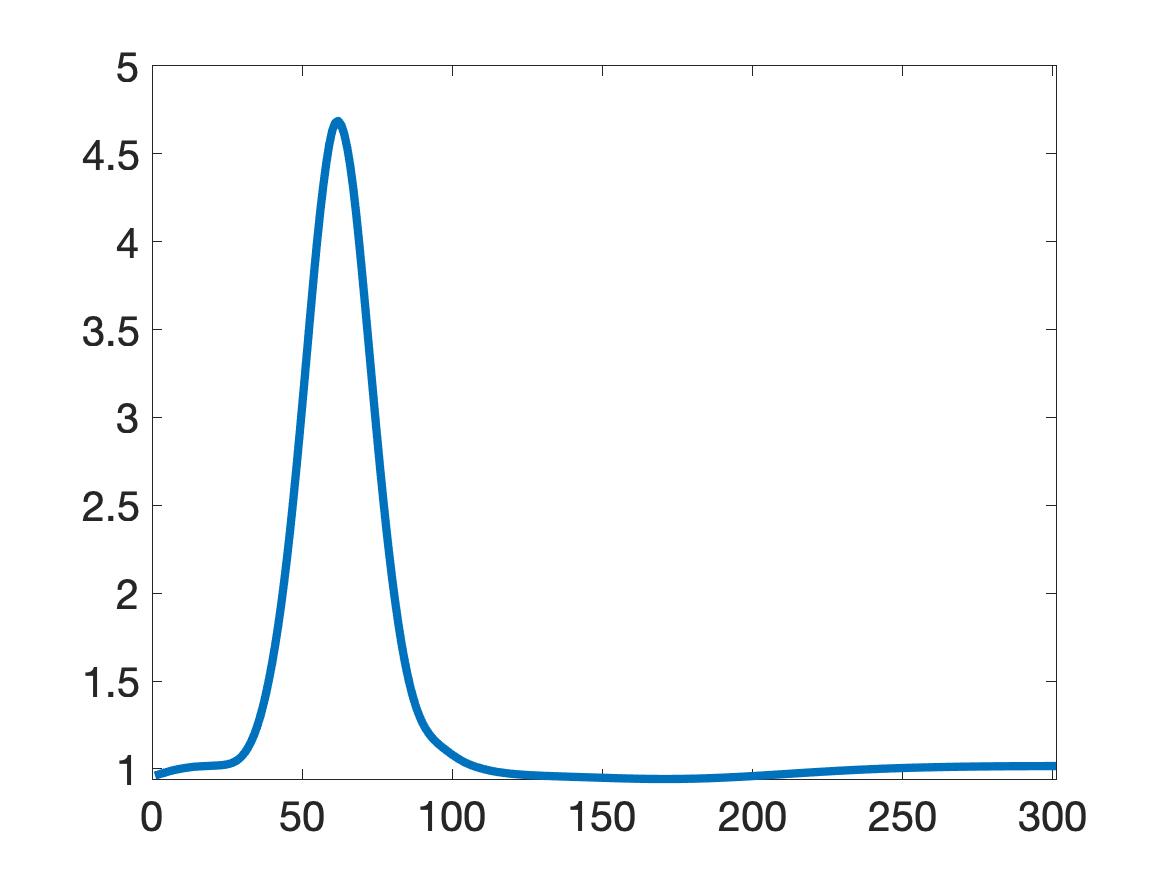} }
\par
\subfloat[\label{fig 5g}The time-dependent raw
data]{\includegraphics[width=.3\textwidth]{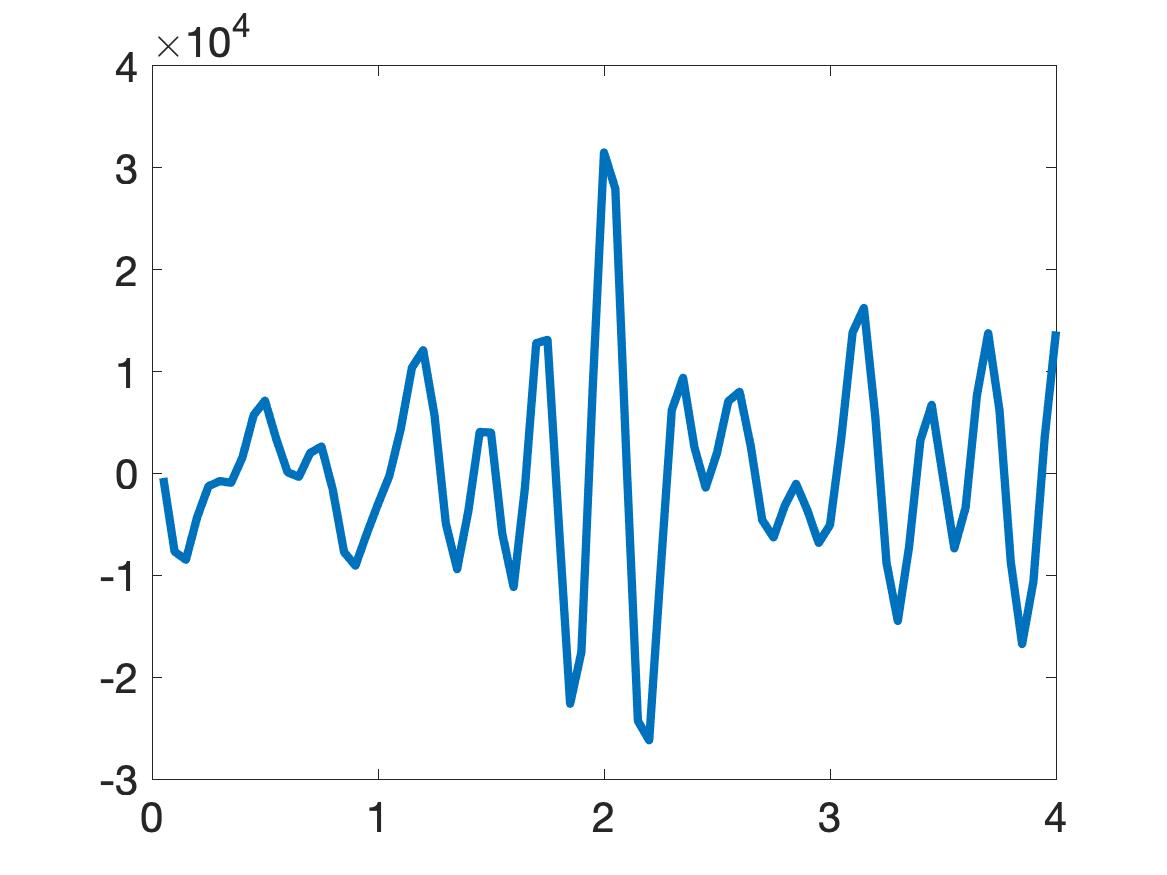}} \quad 
\subfloat[\label{fig 5h}The time-dependent backscattering wave after
preprocessing]{\includegraphics[width=.3\textwidth]{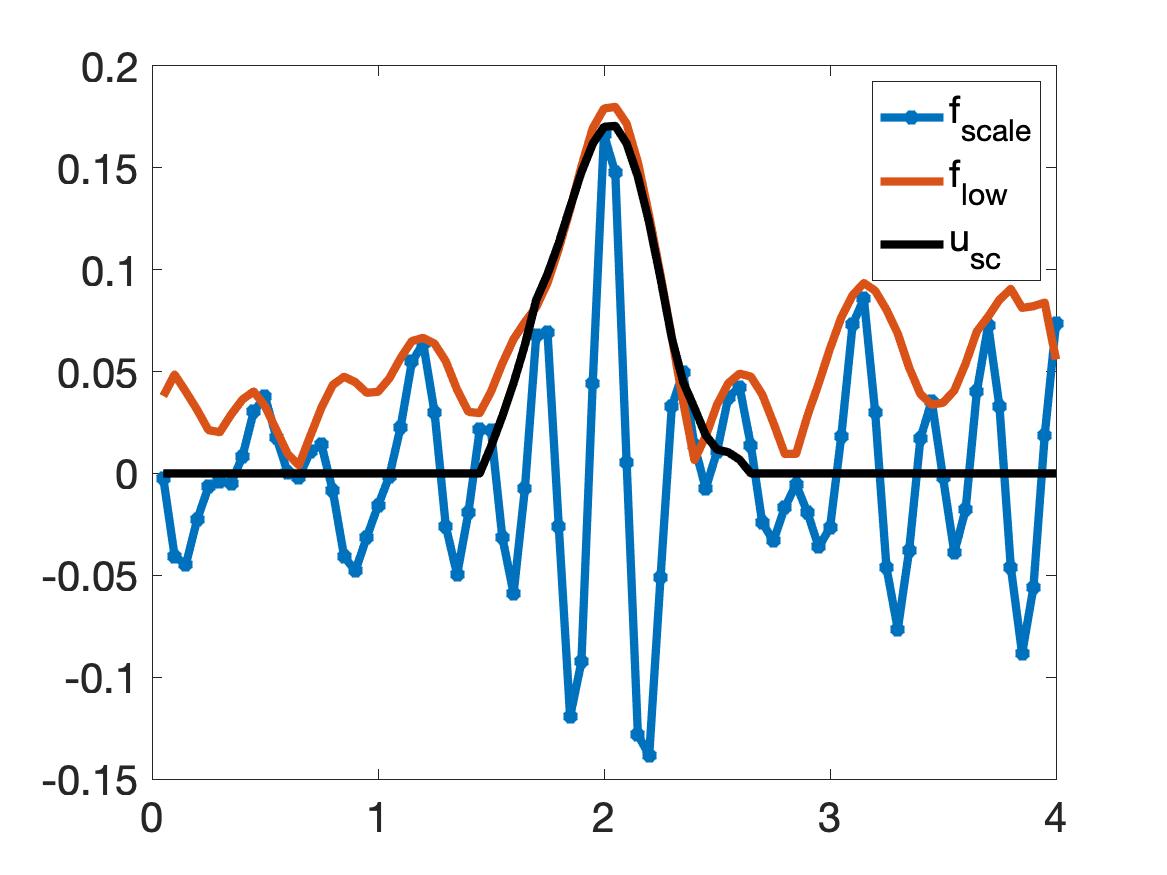}} \quad 
\subfloat[\label{fig 5i}Computed dielectric constant. Its minimal value is
0.37.]{\includegraphics[width=.3\textwidth]{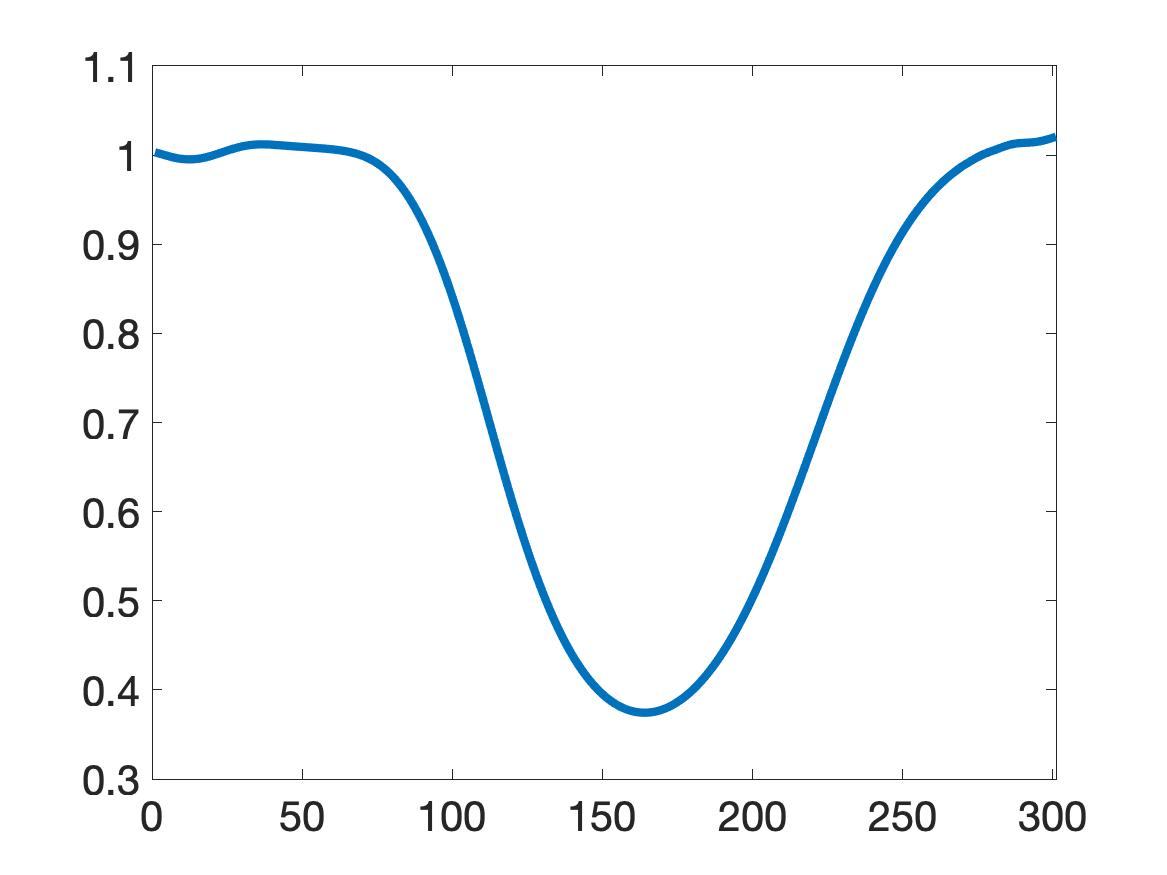} }
\end{center}
\caption{\textit{The case when the target is buried under the ground. The
raw and preprocessed data in the first row correspond to the wave scattered from a
metal box. The raw and preprocessed data in the second row correspond to the wave
scattered from a metal cylinder. The raw and preprocessed data in the third
row correspond to the wave scattered from a plastic cylinder. Unlike the tests in the
first two rows, we choose the upper envelop in this case when preprocessing
the data because $c_{\mathrm{target}} < c_{\mathrm{bckgr}}.$ The computed
dielectric constants for these three tests meet the expectation since they
belong to intervals of their true value, see the last three rows of Table 
\protect\ref{tab 1}. 
}}
\label{fig 5}
\end{figure}

\textit{2. Consider the case when the targets are buried under the ground.}
The background in this case is dry sand. Its dielectric constant is in the
interval $\left[ 3,5\right] .$ We take the average and choose $c_{\mathrm{%
bckgr}}=4.$ We first scale the raw data in the same manner as in case 1 in
which the target is placed in the air. This means that we must know the true
solution of one set of data, generated by a reference target. We use a metal
box with its dielectric constant about 18.5 \cite{KlibanovLoc:ipi2016} as
the reference target. Then, we find a scaling factor $\mu $ such that $\mu
f_{\mathrm{raw}}$ have the same magnitude as the simulated data. The
relative dielectric constant for this reference object is $c_{\mathrm{rel}%
}(x)=4.6$, see (\ref{7.1}). The computed scaling factor is $\mu =189445.$ As
in case 1, the scaled data is denoted by $f_{\mathrm{scale}},$ which is $\mu
f_{\mathrm{raw}}.$ If the dielectric constant of the target $c_{\mathrm{%
target}}$ is larger than that of the background $c_{\mathrm{bckgr}}$, then
the values of the simulated data are less than 0.5 and, hence, the simulated
backscattering wave is non-positive. In this case, we bound $f_{\mathrm{scale%
}}$ by its lower envelop, called $f_{\mathrm{low}}$. If $c_{\mathrm{target}}$
is smaller than $c_{\mathrm{bckgr}}$, the simulated data is larger than 0.5.
In this case, the simulated backscattering wave is non-negative. We hence
bound $f_{\mathrm{scale}}$ by its upper envelop, called $f_{\mathrm{up}}$. A
question arising immediately: how can we know if the dielectric constant of
the target is smaller or larger than that of the background. We answer this
question by experimental observation we got when working with these
experimental data in the past \cite{KlibanovLoc:ipi2016}. We look at the
data and find the three extrema with largest absolute values. If the middle
extremal value among these three is a minimum, then $c_{\mathrm{target}}>c_{%
\mathrm{bckgr}}$. If the middle extreme value is a maximum, then $c_{\mathrm{%
target}}<c_{\mathrm{bckgr}}$. The reader can compare the raw data in Figures %
\ref{fig 5a}, \ref{fig 5d} vs. Figure \ref{fig 5g} for this phenomenon. We
use $f_{\mathrm{envelop}}$ as a common notation for $f_{\mathrm{low}}$ and $%
f_{\mathrm{up}}$. The last step of data preprocessing is the truncation
being applied to $f_{\mathrm{envelop}}$. It is the same as in the truncation
step in case 1. We do not repeat this step here.

The result of the data preprocessing step is the function $g_{0}(t)$ for the
solution of the inverse problem. In comparison with the problem statement in
Problem \ref{cip}, we are missing the knowledge of $g_{1}(t).$ This function
is approximated as follows. Using (\ref{abso1}), we have 
\begin{equation}
u_{x}(x,t)=u_{t}(x,t)\quad \mbox{for all }x<0.  \label{8.2}
\end{equation}%
We accept an error by assuming that (\ref{8.2}) is valid at $0$ in the sense
that we take the limit as $x\rightarrow 0^{-}$. Hence, we can approximate $%
g_{1}(t)=u_{x}(0,t)=u_{t}(0,t)=g_{0}^{\prime }(t)$ for $t>0$. 
%The derivative of $g_0$ is computed by the finite difference method. The step size of the time is $0.05$. This number is taken from \cite{KlibanovLoc:ipi2016}.

\subsection{Numerical results for experimental data}

In this section, we present the numerical results for five (5) tests. The
first two tests are to detect targets in the air and the last three tests
are to identify targets buried a few centimeters under the ground.
Dielectric constants were not measured in these experiments. Therefore, we
have no choice but to compare our computed dielectric constants with those
listed on the website of Honeywell (Table of dielectric constants,
https://goo.gl/kAxtzB). As to the metallic targets, it was numerically
established in \cite{Kuzhuget:IP2012} that one can treat them as dielectrics
with the so-called \textquotedblleft apparent" dielectric constants whose
range is in the interval $\left[ 10,30\right] .$

The reconstructed dielectric constants of these targets are summarized in
Table \ref{tab 1}. 
%In this table, we show the values of computed relative dielectric constant $c_{\rm rel}$ and the computed dielectric constant $c_{\rm target} = c_{\rm rel} c_{\rm bckgr}$.
It can be seen from Table \ref{tab 1} that the computed dielectric constant
of the target and the true one listed on the website of Honeywell (Table of
dielectric constants, https://goo.gl/kAxtzB) are having constant values. In
the table of dielectric constant of Honeywell, the dielectric constant is
not a number. Rather, each dielectric constant of this table is given within
a certain intervals. This interval is listed in the last column of Table \ref%
{tab 1}. It is evident that our computed dielectric constants for all
targets belong to the intervals of the true dielectric constants.
Furthermore, their values are well in the range of those which our group has
computed in previous above cited publications, which have worked with these
experimental data.

\begin{table}[h!]
\begin{center}
\begin{tabular}{|c|c|c|c|c|c|}
\hline
Target & $c_{\mathrm{bckgr}}$ & computed $c_{\mathrm{rel}}$ & $c_{\mathrm{%
bckgr}}$ & computed $c_{\text{target}}$ & True $c_{\text{target}}$ \\ \hline
Bush & 1 & 6.76 & 1 & 6.76 & $[3,20]$ \\ 
Wood stake & 1 & 2.22 & 1 & 2.22 & $[2,6]$ \\ 
Metal box & 4 & 5.2 & $[3,5]$ & $[15.6,26]$ & $[10,30]$ \\ 
Metal cylinder & 4 & 4.7 & $[3,5]$ & $[14.1,23.5]$ & $[10,30]$ \\ 
Plastic cylinder & 4 & 0.37 & $[3,5]$ & $[1.11,1.85]$ & $\left[ 1.1,3.2%
\right] $ \\ \hline
\end{tabular}%
\\[0pt]
\end{center}
\caption{Computed dielectric constants of five targets}
\label{tab 1}
\end{table}

\section{Summary}

\label{sec rem}

We have proposed a new numerical method to solve a highly nonlinear and
severely ill-posed coefficient inverse problem. This method is called the
convexification. Our technique to prove the convexifying phenomenon heavily
relies on a new Carleman estimate, which is proven in Theorem 3.1. The
convexification method has the global convergence property. In fact,
Theorems 4.1, 5.1-5.3 guarantee that the convexification method delivers a
good approximation to the exact solution of the inverse problem without any
advanced knowledge of a small neighborhood of that solution. These results
are verified numerically for both computationally simulated and experimental
data.

\section*{Acknowledgments} 
The work of Klibanov, Le and Loc H. Nguyen , was supported by the US Army
Research Laboratory and US Army Research Office grant W911NF-19-1-0044. The
authors are grateful to Professors Mikhail Kokurin and Oleg Safronov for
their detailed discussions of results of section 4.

%\bibliographystyle{plain}
%\bibliography{../../../../mybib}

\end{document}